\documentclass[reqno]{amsart}
\usepackage[top=1.0in,bottom=1.0in,left=1.0in,right=1.0in]{geometry}
\usepackage{amscd}
\usepackage{amssymb}
\usepackage{multicol}
\usepackage{makecell}
\usepackage[T1]{fontenc}
\usepackage[latin1]{inputenc}
\usepackage{tikz}
\usetikzlibrary{shapes.geometric, arrows}
\usepackage[all]{xy}
\usepackage{times}
\usepackage{comment}

\usepackage{flafter}
\usepackage{placeins}
\usepackage{enumitem}
\usepackage[OT2,T1]{fontenc}
\DeclareSymbolFont{cyrletters}{OT2}{wncyr}{m}{n}
\DeclareMathSymbol{\Sha}{\mathalpha}{cyrletters}{"58}

\usepackage{centernot}
\usetikzlibrary{positioning}
\usepackage{mathrsfs}
\usepackage[color]{showkeys}
\usepackage{url}
\usepackage{centernot}
\usepackage{mathrsfs}
\usepackage[color]{showkeys}
\usepackage{url}
\usepackage{graphicx} 
\usepackage{float} 
\usepackage{subfigure} 

\definecolor{refkey}{rgb}{1,1,1}
\definecolor{labelkey}{rgb}{1,1,1}
\definecolor{cite}{rgb}{0.9451,0.2706,0.4941}
\definecolor{ruri}{rgb}{0.0078,0.4022,0.8010}

\usepackage[%
bookmarks=true,bookmarksnumbered=true,%
colorlinks=true,linkcolor=blue,citecolor=red%
 ]{hyperref}

\def\Hom{{\rm Hom}}
\def\End{{\rm End}}

\def\Aut{{\rm Aut}}

\theoremstyle{plain}

\newtheorem{theorem}{Theorem}[section]

\newtheorem{proposition/example}[theorem]{Proposition/Example}
\newtheorem{proposition}[theorem]{Proposition}
\newtheorem{corollary}[theorem]{Corollary}
\newtheorem{lemma}[theorem]{Lemma}

\theoremstyle{definition}
\newtheorem{definition}[theorem]{Definition}
\newtheorem{remark}[theorem]{Remark}

\newtheorem{example}[theorem]{Example}

\newtheorem{conjecture/question}[theorem]{Conjecture/Question}

\newtheorem{remark/definition}[theorem]{Remark/Definition}
\newtheorem{definition/notation}[theorem]{Definition/Notation}

\numberwithin{equation}{section}

 \theoremstyle{remark}

\numberwithin{equation}{section}

\newcommand{\Ker}{\operatorname{Ker}}

\usepackage{extarrows}\usepackage{bm}

\begin{document}
\title[\textbf{Deformation theory of parabolic representation pairs}]{\textbf{Deformation theory of parabolic representation pairs}}

\author{Zhi Hu}

\address{\textsc{School of Mathematics and Statistics, Nanjing University of Science and Technology, Nanjing 210094, China}}

\email{halfask@mail.ustc.edu.cn }

\author{Pengfei Huang}

\address{\textsc{School of Mathematics, Nanjing University, Nanjing 210093, China}}

\email{pfhwangmath@gmail.com}

\author{Wanmin Yan}

\address{\textsc{School of Mathematics, Nanjing University, Nanjing 210093, China}}

\email{wanminyan@smail.nju.edu.cn}

\author{Runhong Zong}

\address{\textsc{School of Mathematics, Nanjing University, Nanjing 210093, China}}

\email{rzong@nju.edu.cn}


\begin{abstract}
In this paper, we introduce the notions of parabolic representation pairs and the parabolic representation pair variety. We investigate the deformation theory  of parabolic representation pairs. The Zariski tangent space and the tangent quadratic cone of the parabolic representation pair variety are described. By the Riemann--Hilbert--Deligne correspondence, we pro-represent the analytic germs of parabolic representation pair variety  by functors related to certain groupoids of parabolic logarithmic flat bundles. Under suitable assumptions, we prove that the differential graded Lie algebra  (DGLA) controlling the deformation of a parabolic logarithmic flat bundle is partially formal. This leads to the quadraticity of  the subvariety of parabolic representation pair variety, which consists of parabolic representation pairs with fixed eigenvalues of monodromies, at certain generic points. Finally, we construct the moduli space of weighted parabolic representation pairs, and, by means of quiver representation theory, we establish the Kobayashi--Hitchin-type theorem for polystable parabolic representation pairs.
\end{abstract}

\maketitle

\maketitle

\tableofcontents
\section{Introduction}
Let $ \bar X$ be a compact connected oriented Riemann surface of genus $g$, and let $G$ be a complex algebraic group. The representation variety of $\bar X$ is the set $R(\bar \Gamma,G)$ of representations of the fundamental group  $\bar\Gamma:=\pi_1(\bar X)$ of $\bar X$ into $G$. This set inherits a natural structure of an algebraic variety from that on $G$.
 Moreover, there is an action of $G$ on $R(\bar\Gamma,G)$ by conjugation, which corresponds to identifying isomorphic representations. This action is far from being free; indeed, it has varying stabilizers on the subvariety of reducible representations.  For this reason, the orbit space $O(\bar\Gamma,G) = R(\bar\Gamma,G)/G$ is no longer an algebraic variety. To overcome this issue, if $G$ is reductive, we can consider the Geometric Invariant Theory (GIT) quotient $M(\bar\Gamma,G)=R(\bar\Gamma,G)/\!\!/G$, which is usually referred to as  the character variety of $X$ into $G$, and is also known as the moduli space of semisimple representations of $\bar\Gamma$ into $G$. From nonabelian Hodge theory, $M(\bar\Gamma,G)$ can be endowed with a hyperK\"{a}hler structure.

In the case of a punctured Riemann surface, the situation becomes more subtle, requiring the introduction of additional data. Let $D=\{x_1,\cdots,x_n\}$ be a set of $n$ distinct marked points on $\bar X$.
 Denote by $X$ the punctured surface $\bar X\setminus D$, and assume that the Euler number $\chi(X)$ of $X$ is negative. Let $\Gamma$  be the fundamental group $\pi_1(X,b)$ of $X$ with base point $b\in X$. Choose a standard set $\{\alpha_1,\beta_1,\cdots,\alpha_g,\beta_g,\gamma_1,\cdots,\gamma_n\}$ of generators of $\Gamma$, which
 satisfy the single defining  relation
 $$
  [\alpha_1,\beta_1]\cdots[\alpha_g,\beta_g]\gamma_1\cdots\gamma_n=\mathrm{Id}.
 $$
  In particular, $\gamma_i$ is freely homotopic to a small loop that simply winds about
the puncture $x_i$ counterclockwise. The natural embedding $\iota: X\rightarrow \bar X$ induces a homomorphism $\iota_*: \Gamma\rightarrow \bar\Gamma$. The kernel of  $\iota_*$ is the normal closure generated by $\gamma_1,\cdots,\gamma_n$. Let $\Gamma_i\simeq \mathbb{Z}$ be the infinite cyclic subgroup of $\Gamma$ generated by $\gamma_i$.

From now on, we take $G=\mathrm{GL}(r,\mathbb{C})$. Studying the representations of $\Gamma$ into $G$, the monodromy information around punctures must be taken into account. Fix a tuple $\mathbf{C}=(C_1,\cdots,C_n)$ of conjugacy classes  in $G$. The parabolic representation variety $ R(\Gamma, G; \mathbf{C})$ is defined as the collection of representations $\rho: \Gamma\rightarrow G$ such that $\rho(\gamma_i)\in C_i$ for all $i=1,\cdots,n$ \cite{gl}. We can also define the reduced parabolic representation variety $ R(\Gamma, G; \overline{\mathbf{C}})$, as the set of representations $\rho: \Gamma\rightarrow G$ such that $\rho(\gamma_i)\in \overline{C_i}$ for all $i=1,\cdots,n$, where $\overline{C_i}$ denotes the Zariski closure of $C_i$ in $G$. Clearly, $ R(\Gamma, G; \overline{\mathbf{C}})$  is a closed subvariety of $R(\Gamma,G)$, and hence an affine variety. Similarly, the reduced parabolic character variety is defined as $M(\Gamma, G; \overline{\mathbf{C}}):=R(\Gamma, G; \overline{\mathbf{C}})/\!\!/G$, and we define the parabolic character variety $\mathcal{M}(\Gamma, G; {\mathbf{C}})$ as the image of $R(\Gamma, G; {\mathbf{C}})\subset R(\Gamma, G; \overline{\mathbf{C}})$ in $M(\Gamma, G; \overline{\mathbf{C}})$\footnote{If one takes $G=\mathrm{GL}(r,\mathbb{F}_q)$ for a finite field $\mathbb{F}_q$, $ R(\Gamma, G; \mathbf{C})$ is an affine scheme of finite type, then the parabolic character variety can be defined as the corresponding  affine GIT quotient \cite{ka}.} \cite{ll}.  In the case where  all conjugacy classes $C_1,\cdots,C_n$ are semisimple,  $R(\Gamma, G; \mathbf{C})$ coincides with $R(\Gamma, G; \overline{\mathbf{C}})$ and is diffeomorphic to a certain moduli space of semistable parabolic Higgs bundles on $\bar X$ if each $C_i$ has $r$ distinct eigenvalues \cite{hhh}. When $\mathbf{C}$ is generic in the sense of \cite[definition 3.1]{ll}, if $M(\Gamma, G; \overline{\mathbf{C}})$ is non-empty, then $\mathcal{M}(\Gamma, G; {\mathbf{C}})$ is a dense, non-singular, open subset of $M(\Gamma, G; \overline{\mathbf{C}})$, hence is non-empty (cf. \cite[Theorem 1.1]{ll}).

\vspace{5pt}

In this paper, unlike the  considerations above, we not only need the representations $\rho: \Gamma\rightarrow G$, but also keep track of the parabolic groups in which the monodromies around punctures lie. Consequently, a representation with the monodromy around a puncture that lies in the intersection of two different conjugate parabolic subgroups will be recognized as non-equivalent objects once we attach these different parabolic subgroups respectively. More precisely, we consider parabolic representation pairs $(\rho,\mathbf{P})$, where the data $\mathbf{P}$ record the parabolic subgroups to which the $\rho(\gamma_i)$ belong. We then define parabolic representation pair variety $\mathfrak{PRP}(r,\mathbf{d})$ which parametrizes parabolic representation pairs of given parabolic types $\mathbf{d}$. A more suitable definition in nonabelian Hodge theory is to assign the conjugacy classes in Levi parts of the parabolic subgroups in $\mathbf{P}$ (cf. \cite{boa,bo}), which is used in Theorem \ref{tf}.

From the perspective of nonabelian Hodge theory, including the choices of parabolic subgroups associated to the monodromies around  punctures is natural since it corresponds to producing filtered objects \cite{si} via the  choices of parabolic structures on local systems,  flat bundles, or Higgs bundles. In \cite{mk}, Inaba, Iwasaki and Saito studied the Riemann--Hilbert correspondence in the parabolic context. However, they still use the usual (reduced) parabolic character variety, which leads to certain pathological phenomena. In particular, the Riemann--Hilbert correspondence may contract some families of compact subvarieties in the de Rham moduli space onto the singular locus of the character variety if the monodromies are non-generic (especially, they have geometrically  multiple eigenvalues). It is worth pointing out that if one replaces the traditional  parabolic character variety by our parabolic representation pair variety,  these issues are resolved, since varying the parabolic subgroups restores the lost degrees of freedom \cite{hhh}.

The local properties of the classical representation variety and character variety (for a variety of arbitrary dimension and a real algebraic Lie group) have been investigated by many authors. Among the most notable results are those of Goldman and Millson \cite{gm}, who studied the complete local ring $\widehat{\mathbb{O}}_\rho$ at a point $\rho$ in the representation variety for a compact  K\"{a}hler manifold, which pro-represents the functor of deformations of $\rho$.  In particular, they showed that it admits a quadratic presentation, and the similar results were obtained by Nadel for the moduli space of flat bundles \cite{n}. Later,  Eyssidieux and Simpson constructed a functorial mixed Hodge structure on $\widehat{\mathbb{O}}_\rho$ when $\rho$ corresponds to a variation of Hodge structures \cite{es}. Lef\`{e}vre further developed the Eyssidieux--Simpson framework using mixed Hodge diagrams and methods of $L_\infty$-algebras, and extended  Eyssidieux--Simpson's result to the case of smooth quasi-projective variety when $\rho$ has  finite image \cite{lef}. Pridham studied the deformation theory of representations of the algebraic fundamental group \cite{p1,p2}, and construct mixed Hodge structures on real schematic homotopy types of complex quasi-projective varieties, giving mixed Hodge structures on their
homotopy groups and pro-algebraic fundamental groups \cite{p3}.

In this paper, we focus on the local behavior of the parabolic representation pair variety. First, we compute their Zariski tangent spaces and tangent quadratic cones. This leads to the following theorem.

\begin{theorem}
For a representation $\rho: \Xi\rightarrow G$ of a group $\Xi$ into $G$, we introduce the following notations:
\begin{itemize}
  \item  $C^p(\Xi,\mathfrak{g})$ denotes  the abelian group of maps from $\Xi^p$ to $\mathfrak{g}$;
  \item $\delta: C^p(\Xi,\mathfrak{g})\rightarrow C^{p+1}(\Xi,\mathfrak{g})$ is the differential defined via the representation $\rho$;
  \item  $Z^p(\Xi,\mathfrak{g})$ denotes the group of $p$-cocycles;
  \item  $ C^p(\Gamma_i,\mathfrak{g}, \mathfrak{p}_{\mathcal{F}^{(i)}})$ denotes the quotient space  $ C^p(\Gamma_i,\mathfrak{g})/C^p(\Gamma_i,\mathfrak{p}_{\mathcal{F}^{(i)}})$.
\end{itemize}
 Fix a parabolic representation pair $(\rho, \mathbf{P})\in\mathfrak{PRP}(r,\mathbf{d})$, then we have
\begin{enumerate}
  \item  the Zariski tangent space of $\mathfrak{PRP}(r,\mathbf{d})$ at $(\rho, \mathbf{P})$ is given by
      \begin{align*}
        \mathrm{ZT}_{(\rho, \mathbf{P})}\mathfrak{PRP}(r,\mathbf{d})=&\ \{(X, Y_1,\cdots, Y_n)\in Z^1(\Gamma,\mathfrak{g})\bigoplus (\bigoplus\limits_{i=1}^nC^0(\Gamma_i,\mathfrak{g},\mathfrak{p}_{\mathcal{F}^{(i)}}): \\
        &\ \ \ \ \ \ \ \ \ \ \ \ \ \ \ \ \ \ \ \ \ \ \ \ \ \ \ \ \ \ \ \ \ \  \ \ \ \ \ \ \ \ \  [X|_{\Gamma_i}]_i+ \delta  Y_i=0, i=1,\cdots,n\},
      \end{align*}
      and it has complex dimension  $r^2(2g+n-1)$, where $
      [\bullet]_i:C^p(\Gamma_i,\mathfrak{g})\rightarrow  C^p(\Gamma_i,\mathfrak{g},\mathfrak{p}_{\mathcal{F}^{(i)}})
      $
      denotes the natural projection;

  \item the tangent quadratic cone of $\mathfrak{PRP}(r,\mathbf{d})$ at $(\rho, \mathbf{P})$ is given by
  \begin{align*}
   &\mathrm{QC}_{(\rho, \mathbf{P})}\mathfrak{PRP}(r,\mathbf{d})\\
   =&\ \{ (X, Y_1,\cdots, Y_n)\in \mathrm{ZT}_{(\rho, \mathbf{P})}\mathfrak{PRP}(\mathbf{d},r):
    [X,X]_\rho=\delta V\in B^2(\Gamma,\mathfrak{g})\textrm{ for } V\in C^1(\Gamma,\mathfrak{g}),\\
  & \ \ \ \ \ \ \ \ \ \ \ \ \ \ \ \ \ \ \ \ \ \ \ \ \ \ \ \ \ \ \ \ \ \ \     [Y_i, Y_i]_\rho-[V|_{\Gamma_i}]=\delta Z_i\textrm{ for   } Z_i\in \mathfrak{g}/(\mathfrak{p}_{\mathcal{F}^{(i)}}+[ Y_i,  \mathfrak{p}_{\mathcal{F}^{(i)}} ]), i=1,\cdots,n\},
  \end{align*}
  where
  \begin{align*}
   [X,X]_\rho(\omega_1,\omega_2)&=[X(\omega_1),\mathrm{Ad}_{\rho(\omega_1)}X(\omega_2)],\\
   [Y_i, Y_i]_\rho(\gamma_i)&=[ Y_i,\mathrm{Ad}_{\rho(\gamma_i)}Y_i].
  \end{align*}
\end{enumerate}
\end{theorem}

 Next we establish some groupoid versions of the Riemann--Hilbert--Deligne correspondence between certain groupoids of parabolic representation pairs and groupoids of parabolic logarithmic flat bundles.
 This correspondence is achieved via Deligne extension, which can be uniquely lifted to the base of  local Artin
$\mathbb{C}$-algebra by means of  Hensel lemma. However, Deligne extension depends on the choices of argument angles of the eigenvalues of monodromies around punctures\footnote{If one assigns the argument angles to lie in $[0,2\pi)$, it is called the canonical Deligne extension. }. To eliminate this ambiguity, we introduce the notion of RHD-relatedness in the category parabolic logarithmic flat bundles, which identifies two non-isomorphic parabolic logarithmic flat bundles over $(\bar X,D)$ if they give rise to the same filtered local system over $X$. By the Riemann--Hilbert--Deligne correspondence, we can pro-represent the analytic germs of the representation pair variety and the relative representation variety by functors related to certain groupoids of parabolic logarithmic flat bundles. Thus, we obtain the following theorem.

 \begin{theorem}
 Let $(E,\mathbf{F},\nabla)$ be a parabolic logarithmic flat bundle and $(\rho, \mathbf{P})$ the corresponding parabolic representation pair. The groupoids $\mathbb{PFB}'_{( E, {\mathbf{F}}, \nabla)}$ and $\mathbb{HCLC}'_{(C^\infty(E,\mathbf{F}), \bar\partial_E,\nabla)}$ are defined in Sect.~3. Let $\mathbf{Art}_\mathbb{C}$ and $\mathbf{Set}$ denote, respectively, the category of  finite-dimensional local Artin $\mathbb{C}$-algebras that contain unity and have residue field $\mathbb{C}$, and the category of sets.
\begin{enumerate}
  \item The analytic germ of $\mathfrak{PRP}(r,\mathbf{d})$ at $(\rho,\mathbf{P})$ pro-represents the functor $\mathbb{F}_{(E,\mathbf{F},\nabla)}: \mathbf{Art}_\mathbb{C}\rightarrow \mathbf{Set}$ defined by
      \begin{align*}
        A\longmapsto \mathrm{Iso}\mathbb{PFB}'_{( E, {\mathbf{F}}, \nabla)}(A),
      \end{align*}
      where $\mathrm{Iso}\widetilde{\mathbb{PFB}}^{[0,1)'}_{( E, {\mathbf{F}}, \nabla)}(A)$ denotes the set of isomorphism classes in $\mathbb{PFB}'_{( E, {\mathbf{F}}, \nabla)}(A)$.

  \item The analytic germ of $R(X,G;\mathbf{d})$ at $\rho$ pro-represents the functor $\mathbb{G}_{(E,\mathbf{F},\nabla)}: \mathbf{Art}_\mathbb{C}\rightarrow \mathbf{Set}$ defined by
      \begin{align*}
        A\longmapsto \mathrm{Iso}\mathbb{HCLC}'_{(C^\infty(E,\mathbf{F}), \bar\partial_E,\nabla)}(A),
      \end{align*}
      where $\mathrm{Iso}\mathbb{HCLC}'_{(C^\infty(E,\mathbf{F}), \bar\partial_E,\nabla)}(A)$ denotes the set of isomorphism classes in $\mathbb{HCLC}'_{(C^\infty(E,\mathbf{F}), \bar\partial_E,\nabla)}(A)$.
\end{enumerate}
 \end{theorem}

In Sect.~4, we consider the deformations of parabolic logarithmic flat bundles. It is known that Deligne formulated the following philosophy:
in characteristic  0, any deformation problem is controlled by a differential graded Lie algebra  (DGLA) via the Maurer--Cartan equation and gauge equivalence. This philosophy is now a theorem due to Lurie \cite{lur} and Pridham \cite{prid}, using ideas from derived geometry.  We give DGLA descriptions of some deformation functors associated to a parabolic logarithmic flat bundle $(E,\mathbf{F},\nabla)$. Let $\End(E,\mathbf{F})$ denote  the sheaf of germs of endomorphisms of $E$ which preserve the parabolic structure $\mathbf{F}$. The parabolic de Rham complex associated to $(E,\mathbf{F},\nabla)$ is given by
 \begin{align*}
   \mathcal{C}_{(E,\mathbf{F},\nabla)}: \End(E,\mathbf{F})\xlongrightarrow{\nabla}\End(E,\mathbf{F})\otimes_{\mathcal{O}_{\bar X}}\Omega^1_{\bar X},
 \end{align*}
which is regarded as a sheaf of DGLAs. Let $\mathcal{A}_{(E,\mathbf{F},\nabla)}$ be a sheaf of DGLAs providing an acyclic resolution
$$
\mathcal{C}_{(E,\mathbf{F},\nabla)}\longrightarrow{\mathcal{A}_{(E,\mathbf{F},\nabla)}}
$$
of $\mathcal{C}_{(E,\mathbf{F},\nabla)}$ as a morphism of sheaves of DGLA. Then we obtain a DGLA, which is denoted as $\Gamma(\bar X, {\mathcal{A}_{(E,\mathbf{F},\nabla)}})$.
We have the following theorem.

\begin{theorem}
Let $(E,\mathbf{F},\nabla)$ be a parabolic logarithmic flat bundle and $(\rho, \mathbf{P})$ be the corresponding parabolic  representation pair.
   The groupoid $\mathbb{DEF}_{(\rho,\mathbf{P})}(A)$ is defined in Sect. 3, and we define the following  functor  \begin{align*}\overline{\mathrm{ Def}}_{(\rho,\mathbf{P})}: \mathbf{Art}_\mathbb{C}&\longrightarrow \mathbf{Set}\\
 A&\longmapsto\mathrm{Iso} \mathbb{DEF}_{(\rho,\mathbf{P})}(A),
\end{align*}
where $\mathrm{Iso} \mathbb{DEF}_{(\rho,\mathbf{P})}(A)$ denotes the set of isomorphism classes in $\mathbb{DEF}_{(\rho,\mathbf{P})}(A)$. Then the deformation functor $ \mathrm{Def}_{\Gamma(\bar X, {\mathcal{A}_{(E,\mathbf{F},\nabla)}})}$ associated to the DGLA $\Gamma(\bar X, {\mathcal{A}_{(E,\mathbf{F},\nabla)}})$ is isomorphic to the functor $\overline{\mathrm{ Def}}_{(\rho,\mathbf{P})}$.
\end{theorem}

 The basic strategy for proving the above theorem is to apply Goldman--Millson's fundamental theorem, which asserts that quasi-isomorphic DGLAs have isomorphic deformation functors \cite[Theorem 2.4]{gm}. For this purpose, we will construct an intermediate DGLA, the Thom--Whitney--Sullivan DGLA (TWS DGLA), which plays the role of a bridge. The logical chain of the proof is briefly outlined below:
$$\begin{tikzpicture}
\draw [-, thick, color=orange] (-3.5,0.5) -- (-0.4,0.5) -- (-0.4,1.5) -- (-3.5,1.5) -- (-3.5,0.5);
\draw [-, thick, color=orange] (-3.5,-0.5) -- (-0.4,-0.5) -- (-0.4,-1.5) -- (-3.5,-1.5) -- (-3.5,-0.5);
\draw [-, thick, color=orange] (4.3,0.5) -- (1,0.5) -- (1,1.5) -- (4.3,1.5) -- (4.3,0.5);
\draw [-, thick, color=orange] (4.3,-0.5) -- (1,-0.5) -- (1,-1.5) -- (4.3,-1.5) -- (4.3,-0.5);
\draw [-, thick, color=orange] (8.2,-0.5) -- (5.7,-0.5)--(5.7,-1.5)--(8.2,-1.5)--(8.2,-0.5);
\draw [-, thick, color=orange] (12.1,-0.5) -- (9.6,-0.5)--(9.6,-1.5)--(12.1,-1.5)--(12.1,-0.5);
\draw [->, thick] (-0.4,1) -- (1,1);
\draw [->, thick] (-0.4,-1) -- (1,-1);
\draw [<->, thick] (4.3,-1) -- (5.7,-1);
\draw [<->, thick] (8.2,-1) -- (9.6,-1);
\draw [<->, thick] (-2.5,0.45) -- (-2.5,-0.45);
\draw [<->, thick] (3,0.45) -- (3,-0.45);
\node at (-1.9,1.2) {\text{DGLA}};
\node at (-1.9,0.8) {\text{$\Gamma(\bar X, {\mathcal{A}_{(E,\mathbf{F},\nabla)}})$}};
\node at (2.7,1.2) {\text{deformation functor}};
\node at (2.7,0.8) {\text{$ \mathrm{Def}_{\Gamma(\bar X, {\mathcal{A}_{(E,\mathbf{F},\nabla)}})}$}};
\node at (-1.9,-0.8) {\text{TWS DGLA}};
\node at (-1.9,-1.2) {$\mathrm{TWS}_\mathcal{U}(\mathcal{C}_{(E,\mathbf{F},\nabla)})$};
\node at (2.7,-0.8) {\text{functor }};
\node at (2.7,-1.2) {\text{$\mathrm{Def}_{\mathrm{TWS}_\mathcal{U}(\mathcal{C}_{(E,\mathbf{F},\nabla)})}$}};
\node at (7,-0.8) {\text{functor }};
\node at (7,-1.2) {\text{$\mathfrak{H}^1_{\exp(\mathcal{C}_{(E,\mathbf{F},\nabla)})}$}};
\node at (10.9,-0.8) {\text{functor }};
\node at (10.9,-1.2) {\text{$\overline{\mathrm{ Def}}_{(\rho,\mathbf{P})}$}};
\node [right] at (-2.5,0) {};
\node [above] at (0,-1) {};
\end{tikzpicture}$$

 As a consequence, the deformation of a parabolic logarithmic flat bundle $(E,\mathbf{F},\nabla)$  is controlled by the DGLA $M_{(E,\mathbf{F},\nabla)}$, defined as follows:
\begin{itemize}
  \item $M_{(E,\mathbf{F},\nabla)}^i=\bigoplus\limits_{\bullet=p+q}\Gamma(\bar X, \mathcal{A}^{p,q}_{(\bar X,D)}(\End(E,\mathbf{F})))$, $i=0,1,2$;
  \item the bracket is naturally induced from the Lie bracket on $\End(E)$;
  \item the differential is given by $d_\nabla=\bar\partial_E+\nabla$.
\end{itemize}
Similarly, the deformation of a  logarithmic flat bundle $(E,\nabla)$ with the fixed compatible parabolic structure   is controlled by the DGLA $\tilde M_{(E,\mathbf{F},\nabla)}$, where
\begin{itemize}
  \item $\tilde M_{(E,\mathbf{F},\nabla)}^0=M_{(E,\mathbf{F},\nabla)}^0$;
  \item $\tilde M_{(E,\mathbf{F},\nabla)}^1=\Gamma(\bar X,\mathcal{A}^{0,1}_{(\bar X,D)}(\End(E)))\bigoplus\Gamma(\bar X,\mathcal{A}^{1,0}_{(\bar X,D)}(\End(E,\mathbf{F})))$;
  \item $\tilde M_{(E,\mathbf{F},\nabla)}^2=\Gamma(\bar X,\mathcal{A}^{1,1}_{(\bar X,D)}(\End(E)))$.
\end{itemize}
 When all the residues of $\nabla$ at punctures are regular semisimple elements in $\mathfrak{g}$, the deformation functors associated to the above two DGLAs are isomorphic, which implies that the parabolic representation pair variety $\mathfrak{PRP}(r,\mathbf{d})$is non-singular at the corresponding parabolic representation pair $(\rho,\mathbf{P})$.

In parallel with the  famous work of Goldman--Millson \cite{gm}, one may ask whether the DGLA $M_{(E,\mathbf{F},\nabla)}$ becomes formal when the stability condition is imposed on the parabolic logarithmic flat bundle. Unfortunately, in the quasi-projective case,  we lack a sufficiently good Hodge theory to guarantee formality. This is not too discouraging, since we still have the $L_2$-Hodge theory for harmonic bundles at hand. But this leads to a new problem: the $L_2$-condition is not preserved under Lie brackets. This  motivates us  to introduce the notion of partial formality. A DGLA $L$ is said to be partially formal if it  contains a formal sub-DGLA satisfying some extra conditions. If the DGLA $L$ is used to describe the local behavior  of a variety associated to  some moduli problem, then the  partial formality of $L$  is related to the local quadraticity  of a certain subvariety.  Under some assumptions (especially Jordan stability of logarithmic flat bundles), we show the partial formality of $M_{(E,\mathbf{J},\nabla)}$, where $\mathbf{J}$ is a parabolic structure arising from the Jordan decompositions of the residues of $\nabla$.

\begin{theorem}
Let  $(E,\nabla)$ be a logarithmic flat bundle  over $(\bar X,D)$. Assume that
\begin{itemize}
  \item  the residues of $\nabla$ at the punctures are semisimple and all their eigenvalues lie in $[0,1)$;
    \item  $(E,\nabla)$ is Jordan stable,
    \end{itemize}
then DGLA $M_{(E,\mathbf{J},\nabla)}$ is partially formal.
\end{theorem}
From this theorem, we conclude  the quadraticity of a certain subvariety of $\mathfrak{PRP}(r,\mathbf{d})$.
\begin{theorem}\label{tf}
Fixing $\mathbf{e}=\{\overrightarrow{e^{(i)}}\}_{i=1,\cdots,n}=\{(e^{(i)}_1,\cdots,e^{(i)}_{r})\}_{i=1,\cdots,n}$ for  $e^{(i)}_\mu=e^{2\pi\sqrt{-1}a^{(i)}_\mu}$ with $a^{(i)}_\mu\in[0,1)$, let $\mathfrak{PRP}(r,\mathbf{d},\mathbf{e})$ be a subvariety of $\mathfrak{PRP}(r,\mathbf{d})$ consisting of parabolic representation pairs $(\rho,\mathbf{P})$ such that the eigenvalues of $\rho(\gamma_i)$ are exactly $e^{(i)}_1,\cdots,e^{(i)}_{r}$ for each $i$.
 If $\mathfrak{PRP}(r,\mathbf{d},\mathbf{e})$ contains a parabolic representation pair $(\rho,\mathbf{P})$ corresponding to a  logarithmic flat bundle as in the above theorem, then
$\mathfrak{PRP}(r,\mathbf{d},\mathbf{e})$ is quadratic at $(\rho,\mathbf{P})$.
\end{theorem}

Finally, in Sect. 5, since affine GIT is no longer applicable to the parabolic representation pair variety, we construct the moduli space of weighted parabolic representation pairs by translating these objects into representations of a certain star-shaped quiver with loops. It is worth emphasizing that the monotonicity of the weight vector at each puncture plays a crucial role, as it ensures that the stability of parabolic representation pairs coincides with the stability of representations of the corresponding quiver. As an application, we have the following Kobayashi--Hitchin-type theorem,  which can be regarded as a generalization of Totaro's result \cite[Theorem 2]{to} on polystable filtered vector spaces introduced by Faltings \cite{fa}.

\begin{theorem}
Given a rational weight system $\mathbf{w}=\{\overrightarrow{w^{(i)}}\}_{i=1,\cdots,n}$, where each $\overrightarrow{w^{(i)}}=(w^{(i)}_0,\cdots, w^{(i)}_{s_i})\in \mathbb{Q}^{s_i+1}$ satisfies $w^{(i)}_0<\cdots <w^{(i)}_{s_i}$. Then a  parabolic representation pair $(\rho,\mathbf{P})\in \mathfrak{PRP}(r,\mathbf{d})$ is polystable with respect to $\mathbf{w}$ if and only if there exists a $\mathbf{w}$-good Hermitian metric $h$ on $V$, namely the following equations are satisfied
\begin{align*}
    \left\{
    \begin{aligned}
        &\sum_{j=1}^g ([\rho(\alpha_j),\rho^{*_h}(\alpha_j)]+[\rho(\beta_j),\rho^{*_h}(\beta_j)])+\sum_{i=1}^n[\rho(\gamma_i),\rho^{*_h}(\gamma_i)]+\sum\limits_{i=1}^n{p_h}^{(i)}_1\\
 &\hspace{138pt}=\ (\frac{1}{r} \sum\limits_{i=1}^n\sum\limits_{\ell=0}^{s_i}w_{\ell}^{(i)}d_{\ell}^{(i)}-\sum\limits_{i=1}^nw_{0}^{(i)})\mathrm{Id},\\[3pt]
 &[\rho(\gamma_i)|_{V^{(i)}_\ell},(\rho(\gamma_i)|_{V^{(i)}_\ell})^{*_h}]+{p_h}^{(i)}_\ell
  =\ (w^{(i)}_{\ell-1}-w^{(i)}_\ell+1)\mathrm{Id},\quad 1\leq \ell\leq s_i, 1\leq i\leq n,
    \end{aligned}
    \right.
\end{align*}
where $\rho^{*_h}(\omega)$ and ${p_h}^{(i)}_\ell$ denote, respectively, the conjugate of $\rho(\omega)$ and the orthogonal projection  of $V^{(i)}_{\ell}$  onto $V^{(i)}_{\ell+1}$ with respect to the metric $h$. Moreover, if $h'$ is another $\mathbf{w}$-good Hermitian metric  on $V$, then there exists $g\in G$ such that
\begin{align*}
  h(u,v)=h'(g\cdot u,g\cdot v)
\end{align*}
for any $u,v\in V$.
\end{theorem}


\section{Deformations of parabolic representation pairs}

\subsection{Parabolic representation pairs}

Identifying $G$ with $\mathrm{GL}(V)$ for an $r$-dimensional complex vector space $V$ and writing  $\mathfrak{g}=\mathfrak{gl}(V)$, a subgroup $P$ of $G$ is called parabolic if one of the following equivalent conditions holds:
\begin{itemize}
  \item $P$ contains a maximal connected solvable subgroup
of $G$;
  \item $P$ consists of all automorphisms of $V$ that preserve a fixed flag $ \mathcal{F}=\{V_\ell\}_{\ell=0,\cdots,s}$, i.e. as a nested sequence of nonzero proper subspaces of $V$
  \begin{align*}
\{0\}=V_{s+1}\subsetneq V_s\subsetneq\cdots\subsetneq V_{1}\subsetneq V_{0}=V
  \end{align*}
   of type $\overrightarrow{d}=(d_0,\cdots,d_{s})$, where $d_\ell=\dim_\mathbb{C}V_{\ell}/V_{\ell+1}$ with $\sum\limits_{\ell=0}^{s}d_\ell=r$;
  \item the quotient $G/P$ has the structure of a projective variety, which is a flag variety consisting of all flags in  $V$ of some type $\overrightarrow{d}$.
\end{itemize}
Hence, any parabolic subgroup $P$ of $G$ is uniquely associated to a flag  $\mathcal{F}$; we also denote $P$ by $P_\mathcal{F}$, and the type of $\mathcal{F}$ is also called the type of $P$. Up to conjugate by $G$, a parabolic subgroup of $G$ is uniquely determined by the numerical invariant: type. Denote by $\mathrm{FL}(\overrightarrow{d},r)$ the flag variety consisting of flags of type $\overrightarrow{d}$ in an $r$-dimensional vector space $V$, whose elements correspond parabolic subgroups of $G$ of type $\overrightarrow{d}$.

Let $X=\bar X\setminus D$ be a punctured surface. We define the representation variety
\begin{align*}
 R(\Gamma,G)=\textrm{Hom}(\Gamma, G)
\end{align*}
which consists of the representations of $\Gamma$ on $V$. It is isomorphic to the non-singular  affine variety $G^{1-\chi(X)}$. We also define the multiple flag variety \cite{sm}
\begin{align*}
 \mathrm{FL}(r,\mathbf{d})=\mathrm{FL}(r,\overrightarrow{d^{(1)}})\times \cdots\times \mathrm{FL}(r,\overrightarrow{d^{(n)}})
\end{align*}
for $\mathbf{d}=\{\overrightarrow{d^{(i)}}\}_{i=1,\cdots,n}=\{(d^{(i)}_0,\cdots,d^{(i)}_{s_i})\}_{i=1,\cdots,n}$ with $\sum\limits_{\ell=0}^{s_i}d^{(i)}_{\ell}=r$, which is a non-singular Fano projective variety.

\begin{definition}\
\begin{enumerate}
  \item A \emph{parabolic representation pair} of rank $r$ and type $\mathbf{d}=\{\overrightarrow{d^{(i)}}\}_{i=1,\cdots,n}$ is a pair
 \begin{align*}
  (\rho,\mathbf{P})\in R(\Gamma,G)\times\mathrm{FL}(r,\mathbf{d}),
 \end{align*}
 where $\rho\in R(\Gamma,G)$, and  $\mathbf{P}\in \mathrm{FL}(r,\mathbf{d})$ is given by $\mathbf{P}=(P_{\mathcal{F}^{(1)}},\cdots,P_{\mathcal{F}^{(n)}})\in\mathrm{FL}(r,\mathbf{d}) $ for parabolic subgroups $P_{\mathcal{F}^{(i)}}\in \mathrm{FL}(\overrightarrow{d^{(i)}},r)$ of type $\overrightarrow{d^{(i)}}$, such that $\rho(\gamma_i)\in P_{\mathcal{F}^{(i)}}$ holds for each $i$.
The subset of $R(\Gamma,G)\times\mathrm{FL}(r,\mathbf{d})$ consisting of parabolic representation pairs is denoted by $\mathfrak{PRP}(r,\mathbf{d})$\footnote{We will employ the same notation in different contexts, such as sets, varieties, and categories.}.
  \item Two parabolic representation pairs $ (\rho,\mathbf{P})$ and $(\rho',\mathbf{P}')$ are \emph{equivalent} if there exists $g\in G$ such that $\rho'(\omega)=g\rho(\omega)g^{-1}$ for any $\omega\in \Gamma$ and $P'_{\mathcal{F}^{(i)}}=gP_{\mathcal{F}^{(i)}}g^{-1}$ for $i=1,\cdots, n$.
  \item  The \emph{transformation groupoid} corresponding to $(\mathfrak{PRP}(r,\mathbf{d}),G)$ is denoted by $\mathbb{PRP}(r,\mathbf{d})$.
\end{enumerate}

\end{definition}

\begin{proposition}
$\mathfrak{PRP}(r,\mathbf{d})$ is a closed subvariety of $R(\Gamma,G)\times\mathrm{FL}(r,\mathbf{d})$, called the \emph{parabolic representation pair variety}.
\end{proposition}
\begin{proof}
  We have a surjective forgetful map
  \begin{align*}
  f_2: \mathfrak{PRP}(r,\mathbf{d})&\longrightarrow \mathrm{FL}(r,\mathbf{d})\\
    (\rho,\mathbf{P})&\longmapsto \mathbf{P}.
  \end{align*}
  It is clear that $f_2$ is a fibration. Given  $\mathbf{P}=(P_{\mathcal{F}^{(1)}},\cdots,P_{\mathcal{F}^{(n)}})$, it suffices to show  the preimage $f_2^{-1}(\mathbf{P})$ is a closed subvariety of $R(\Gamma,G)$.This follows thta fact that  every parabolic subgroup of $G$ is Zariski-closed in $G$.
\end{proof}
\begin{remark}We call  $R(\Gamma,G;\mathbf{P}):=f_2^{-1}(\mathbf{P})$ the relative representation variety of type $\mathbf{P}$. For this object, one can ask whether there exists a moduli space with respect to the transformation group $\bigcap\mathbf{P}$, which is the intersection of all parabolic groups in $\mathbf{P}$. This question is very complicated  because, in general, $\bigcap\mathbf{P}$ is non-reductive, the classical GIT does not apply to the construction of the corresponding moduli spaces.  A natural approach to resolving this issue is to apply non-reductive GIT, as developed in \cite{bdh,bd,bk}.
\end{remark}
In \cite[Section 3]{hs}, the authors define a subset $R(\Gamma,G;\mathbf{d})$ of $R(\Gamma,G)$, which  consists of representations $\rho:\Gamma\rightarrow G$ satisfying that  there exists a flag $\mathcal{F}^{(i)}\in \mathrm{FL}(r,\overrightarrow{d}^{(i)})$ with $\rho(\gamma_i)\in P_{\mathcal{F}^{(i)}}$ for each $i$. Obviously,  $R(\Gamma,G;\mathbf{d})$ is a  closed subvariety of $R(\Gamma,G)$, and in general, it is very singular (cf. \ref{co}-(2)).  There is another surjective forgetful map given by
  \begin{align*}
  f_1: \mathfrak{PRP}(r,\mathbf{d})&\longrightarrow R(\Gamma,G;\mathbf{d})\\
    (\rho,\mathbf{P})&\longmapsto \rho.
  \end{align*}
  The preimage $f_1^{-1}(\rho)$ is a closed subvariety of $\mathrm{FL}(r,\mathbf{d})$ consisting of the fixed points of the diagonal $\rho(\gamma_1)\times\cdots\times \rho(\gamma_n)$-action on $\mathrm{FL}(r,\mathbf{d})$, and in particular, for a generic representation $\rho$,   $f_1^{-1}(\rho)$ is exactly a single point. We call $R(\Gamma,G;\mathbf{d})$ the \emph{extended representation variety} of type $\mathbf{d}$,
  and correspondingly the \emph{extended character variety} of type $\mathbf{d}$ is defined as
  $$
  M(\Gamma, G; \mathbf{d}):=R(\Gamma, G; \mathbf{d})/\!\!/G.
  $$
The transformation groupoid corresponding to $(R(\Gamma,G;\mathbf{d}),G)$ is denoted by $\mathbb{R}(\Gamma,G;\mathbf{d})$.

In summary, we have the following commutative diagram
  \begin{align*}
   \CD
  \mathrm{FL}(r,\mathbf{d})@=\mathrm{FL}(r,\mathbf{d})\\
   @A f_2 AA @A p_2 AA \\
      \mathfrak{PRP}(r,\mathbf{d})@>\iota>>R(\Gamma,G)\times\mathrm{FL}(r,\mathbf{d}) \\
     @V f_1 VV @V p_1 VV  \\
     R(\Gamma,G;\mathbf{d}) @>\iota>> R(\Gamma,G)
   \endCD,
  \end{align*}
where $\iota$ denotes the natural embedding, $p_1$ and $p_2$ denote the natural projections.

\begin{example}\
\begin{enumerate}
     \item Assume $r\geq n$.  We call $\mathbf{P}$  a \emph{braid} \cite{c} if
 \begin{itemize}
   \item all $P_{\mathcal{F}^{(i)}}$ are  Borel subgroups of $G$, i.e. all $\mathcal{F}^{(i)}$ are complete flags on $V$;
   \item there is an order $(x_{i_1},\cdots,x_{i_n})$ of $n$ punctures such that the following conditions are satisfied
\begin{itemize}
\item the two flags $\mathcal{F}^{(i_l)}$ and $\mathcal{F}^{(i_{l+1})}$ are in $\jmath_{i_l}$-position for certain $\jmath_{i_l}\in\{1,\cdots, r-1\}$, i.e.
$$
V^{(i_l)}_\ell\left\{
               \begin{array}{ll}
                      =V^{(i_{l+1})}_\ell, & \hbox{$\ell\neq \jmath_{i_l}$,} \\
                      \neq V^{(i_{l+1})}_\ell, & \hbox{$\ell=\jmath_{i_l}$,}
               \end{array}
               \right.
$$
 \item the Borel subgroups $P_{\mathcal{F}^{(i_1)}}$ and $P_{\mathcal{F}^{(i_n)}}$ are opposite.
 \end{itemize}
   \end{itemize}
 One says a parabolic representation pair $ (\rho,\mathbf{P})$ is a \emph{braid representation pair} if $\mathbf{P}$ is a braid.

   \item  Assume $g=0$ and fix $\mathbf{P}\in \mathrm{FL}(r,\mathbf{d})$. The representations $\rho\in R(\Gamma,G;\mathbf{P})$ are closely related to the \emph{parabolic Deligne--Simpson problem}\footnote{For the classical Deligne--Simpson problem, please see \cite{cr}.}: do there exist matrices  $A_i\in P_{\mathcal{F}^{(i)}}, i=1,\cdots,n$, which do not admit common nontrivial proper invariant subspaces, satisfy $A_1\cdots A_n=\mathrm{Id}$? For example, take $r=2, n=3$. Consider the following flags on the vector space $V=\mathrm{Span}_{\mathbb{C}}\{e_1,e_2\}$:
                   \begin{align*}
                    \mathcal{F}_1: &\ 0\subset\mathrm{Span}_{\mathbb{C}}\{e_1\} \subset V,\\
                   \mathcal{F}_2: &\ 0\subset\mathrm{Span}_{\mathbb{C}}\{e_2\} \subset V,\\
                   \mathcal{F}_3: &\ 0\subset\mathrm{Span}_{\mathbb{C}}\{\frac{3+\sqrt{5}}{2}e_1+\frac{1+\sqrt{5}}{2}e_2\} \subset V,\\
                   \mathcal{F}_4: &\ 0\subset\mathrm{Span}_{\mathbb{C}}\{\frac{1+\sqrt{5}}{2}e_1+\frac{3+\sqrt{5}}{2}e_2\} \subset V.
                   \end{align*}
If one takes $P_{\mathcal{F}^{(1)}}=P_{\mathcal{F}_1}$ and $P_{\mathcal{F}^{(2)}}=P_{\mathcal{F}^{(3)}}=P_{\mathcal{F}_2}$, then the parabolic Deligne--Simpson problem admits no solution. In contrast, if one takes $P_{\mathcal{F}^{(1)}}=P_{\mathcal{F}_1}$, $P_{\mathcal{F}^{(2)}}=P_{\mathcal{F}_3}$, and $P_{\mathcal{F}^{(3)}}=P_{\mathcal{F}_4}$, then the parabolic Deligne--Simpson problem has solutions. For instance, one can take
$$
  A_1=\left(\begin{array}{cc}
           1 & -2\\
           0 & 1 \\
        \end{array}
        \right),\quad
    A_2=\left(\begin{array}{cc}
              1 & 1\\
              1 & 0 \\
          \end{array}
         \right),\quad
    A_3=\left(\begin{array}{cc}
              0 & 1\\
              1 & 1 \\
          \end{array}
          \right).
$$
\item Assume $g\geq 1$, $n$ and $r$ are coprime. For any $\mathbf{d}$,  the extended character variety $M(\Gamma, G; \mathbf{d})$ contains a smooth subvariety of dimension $r^2(2g-2)+2$, given by
$$T=\{\rho\in R(\Gamma,G):\rho(\gamma_i)=e^{\frac{2\pi\sqrt{-1}}{r}}\mathrm{Id}, i=1\cdots, n\}/\!\!/G.$$
$T$ is recognized as a certain twisted character variety \cite{hv}, hence by nonabelian Hodge correspondence, it is diffeomorphic to the moduli space of stable Higgs bundles of rank $r$ and degree $n$ over $\bar X$.
\end{enumerate}
\end{example}

\begin{definition}[\cite{si,tm2}]
Let $\mathbb{L}$ be a local system of rank $r$ over $X$.
\begin{enumerate}
    \item A \emph{filtered structure} of type $\mathbf{d}$ on $\mathbb{L}$ a collection $\{\mathbb{F}_{U_i}\}_{i=1,\cdots,n}$, where each $\mathbb{F}_{U_i}$ is a decreasing filtration of sub-local systems of $\mathbb{L}|_{U_i}$ on an appropriate punctured neighborhood $U_i$ of $x_i$:
    \begin{align*}
  \mathbb{F}_{U_i}:
        \mathbb{L}|_{U_i}= \mathbb{L}^{(U_i)}_{0}\supsetneq \mathbb{L}^{(U_i)}_{1}\supsetneq\cdots\supsetneq\mathbb{L}^{(U_i)}_{s_i}\supsetneq\mathbb{L}^{(U_i)}_{s_i+1}=0
  \end{align*}
  such that the rank of $\mathbb{L}^{(U_i)}_{\ell}/\mathbb{L}^{(U_i)}_{\ell+1}$ is $d^{(i)}_\ell$;
\item Two filtered structures $\{\mathbb{F}_{U_i}\}_{i=1,\cdots,n}$ and $\{\mathbb{F}_{U_i'}\}_{i=1,\cdots,n}$ are said to be \emph{equivalent} if for each $i$, there exists a punctured neighborhood $U_i''$ of $x_i$ with $U_i''\subseteq U_i\cap U_i'$ such that
\begin{align*}
    \mathbb{L}_j^{(U_i)}|_{U_i''}=\mathbb{L}_j^{(U_i')}|_{U_i''}\ \ \text{for any}\ j;
\end{align*}
\item A \emph{filtered local system} of rank $r$ and type $\mathbf{d}$ over $X$ is pair $(\mathbb{L},\mathbf{F})$ consisting of a local system $\mathbb{L}$ of rank $r$ together with an equivalence class $\mathbf{F}=\{\mathbb{F}_{U_i}\}_{i=1,\cdots,n}$ of filtered structures of type $\mathbf{d}$ on $\mathbb{L}$.
\end{enumerate}
\end{definition}

Let $\mathfrak{ FLS}(r,\mathbf{d})$ be the category of filtered local systems of rank $r$ and type $\mathbf{d}$ over $X$, then it is equivalent to the category $\mathfrak{PRP}(r,\mathbf{d})$ of parabolic representation pairs of the same rank and type.
Indeed, given a parabolic representation pair $(\rho, \mathbf{P})\in\mathfrak{ PRP}(r,\mathbf{d})$, one defines
  \begin{align*}
   \mathbb{L}^\rho=\tilde X\times V/{\sim_{(\Gamma,\rho)}}
  \end{align*}
  where $\pi: \tilde X\rightarrow X$ is the universal cover of $X$, and the equivalence relation is defined by
$(x,v)\sim_{(\Gamma,\rho)}(x',v')$ if there exists $\omega\in \Gamma$ such that $x'=\omega\cdot x$ and $v'=\rho(\omega)\cdot v$. Thus, $\mathbb{L}^\rho$ is a local system over $X$.
Similarly, one defines a local system over $U_i$ as follows
\begin{align*}
   {\mathbb{L}^\rho}^{(i)}_\ell=\pi^{-1}( U_i)\times V_\ell/\sim_{\gamma_i},
  \end{align*}
  where $(x,v)\sim_{\gamma_i}(x',v')$ if there exists an integer $k$ such that $x'=\gamma_i^k\cdot x$ and  $v'=\rho(\gamma_i^k)\cdot v$. Thus, $\{ {\mathbb{L}^\rho}^{(i)}_\ell\}$ forms a filtration ${\mathbb{F}^\rho}_{U_i}$ of sub-local systems of $\mathbb{L}^\rho|_{U_i}$. One easily checks that the functor given by
\begin{align*}
(\rho,\mathbf{P})\mapsto (\mathbb{L}^\rho, \mathbf{F}^\rho=\{{\mathbb{F}^\rho}_{U_i}\}_{i=1,\cdots,n})
\end{align*}
provides an equivalence between  $\mathfrak{PRP}(r,\mathbf{d})$  and  $\mathfrak{ FLS}(r,\mathbf{d})$.

\begin{definition}\
\begin{enumerate}
  \item A \emph{weighted parabolic representation pair} $((\rho,\mathbf{P}),\mathbf{w})$ is a parabolic representation pair $(\rho,\mathbf{P})$ endowed with a weight system $\mathbf{w}=\{\overrightarrow{w^{(i)}}\}_{i=1,\cdots,n}$, where $\overrightarrow{w^{(i)}}=(w^{(i)}_0,\cdots, w^{(i)}_{s_i})\in \mathbb{R}^{s_i+1}$ with $w^{(i)}_0<\cdots <w^{(i)}_{s_i}$. The \emph{degree} and \emph{slope} of  $((\rho,\mathbf{P}),\mathbf{w})$ are defined, respectively, by
\begin{align*}
 \deg((\rho,\mathbf{P}),\mathbf{w})&=\sum_{i=1}^n\overrightarrow{d^{(i)}}\cdot\overrightarrow{w^{(i)}}=\sum_{i=1}^n\sum_{\ell=0}^{s_i}{d^{(i)}_\ell}{w^{(i)}_\ell},\\
 \mu((\rho,\mathbf{P}),\mathbf{w})&=\frac{\deg((\rho,\mathbf{P}),\mathbf{w})}{r}.
\end{align*}
  \item Let $\rho'$ be a sub-representation of $\rho$ on an $r'$-dimensional subspace $V'\subset V$.
   Then
   $$
   \mathcal{F}^{(i)}\bigcap V'=\{V_\ell\bigcap V'\}_{\ell=1,\cdots,s_i}
   $$
   can be made into a flag ${\mathcal{F}^{(i)}}'=\{V'_{\ell'}\}_{\ell'=0,\cdots,s_i'}$ in $V'$ of type $\overrightarrow{{d^{(i)}}'}$ by removing duplicate copies. We define
   $\mathbf{P}'=(P_{{\mathcal{F}^{(1)}}'},\cdots,P_{{\mathcal{F}^{(n)}}'})\in\mathrm{FL}(\{\overrightarrow{{d^{(i)}}'}\}_{i=1,\cdots,n},r') $. For each $i$ and $\ell'$, we define
   $$
   S^{(i)}(\ell')=\{\ell\in{0,\cdots, s_i}: \dim_\mathbb{C}V_\ell\bigcap V'=\dim_\mathbb{C}V'_{\ell'}\},
   $$
   and the weight system $\mathbf{w}'=\{\overrightarrow{{w^{(i)}}'}\}_{i=1,\cdots,n}=\{(w^{{(i)}'}_0,\cdots, w^{{(i)}'}_{s_i'})\}_{i=1,\cdots,n}$ by
   $
    w^{{(i)}'}_{\ell'}= w^{(i)}_{\max\{S^{(i)}(\ell')\}}$.
  Therefore, $((\rho',\mathbf{P}'),\mathbf{w}')$ is a weighted parabolic representation pair, called the \emph{induced weighted parabolic sub-representation pair}. A weighted parabolic representation pair $((\rho,\mathbf{P}),\mathbf{w})$ is called \emph{stable} (respectively, \emph{semistable}) if for any induced weighted parabolic sub-representation pair $((\rho',\mathbf{P}'),\mathbf{w}')$ we have the inequality
   \begin{align*}
      \mu((\rho',\mathbf{P}'),\mathbf{w}')<(\textrm{respectively}, \leq ) \ \mu((\rho,\mathbf{P}),\mathbf{w}).
   \end{align*}
    A weighted parabolic representation pair $((\rho,\mathbf{P}),\mathbf{w})$ is called \emph{polystable}, if $\rho$ decomposes as a direct sum of subrepresentations $\rho=\bigoplus\limits_{\mu}\rho_\mu$ such that the induced weighted parabolic sub-representation pair associated to each summand $\rho_\mu$ is  of the same slope and stable.
\end{enumerate}
\end{definition}

\begin{definition}\label{n}
A pair $(\rho:\Gamma\rightarrow G,\mathbf{P}=(P_1,\cdots,P_n))$ is called a stable/semistable/polystable \emph{unweighted parabolic representation pair of degree zero} if
 \begin{itemize}\item
 $\rho(\gamma_i)\in P_i$ for each $i$;
   \item there is a weight system $\mathbf{w}$ such that  \begin{itemize}
     \item  $\deg((\rho,\mathbf{P}),\mathbf{w})$=0,
     \item the weighted parabolic representation pair $((\rho,\mathbf{P}),\mathbf{w})$ is stable/semistable/polystable.
   \end{itemize}
 \end{itemize}
  \end{definition}

The remarkable tame nonabelian Hodge correspondence due to Simpson \cite{si} (also cf. \cite{tm}) predicts the equivalence of the following categories:
\begin{itemize}
  \item the category of stable unweighted parabolic representation pairs of degree zero;
  \item the category of stable filtered local system of degree zero over $X$;
  \item the category of stable filtered logarithmic Higgs bundles of degree zero over $\bar X$;
  \item the category of stable filtered logarithmic flat bundles of degree zero over $\bar X$.
\end{itemize}

 At the end of this section, we briefly indicate how to generalize the above notions to the principal objects. In particular, to establish Simpson's tame nonabelian Hodge correspondence \cite{bo,hksz}, the definition of stability plays a crucial role. Let $H$ be a complex reductive group with Lie algebra $\mathfrak{h}$.
Fix a maximal torus $T$ of $H$ with Lie algebra $\mathfrak{t}$, and let $C_T=\mathrm{Hom}(\mathbb{C}^*,T)$ be the group of cocharacters of $T$, regarded as a lattice in $\mathfrak{t}$. An element of $\mathfrak{t}_\mathbb{R}=C_T\otimes_{\mathbb{Z}}\mathbb{R}$ is called a \emph{weight}. Given a weight $\theta\in \mathfrak{t}_\mathbb{R}$, the associated parabolic subgroup $P_\theta$ is defined by
\begin{align*}
{Q}_\theta=\{ h\in {H}: z^{\theta}  hz^{-\theta}\textrm{ has a limit as}\ z\rightarrow0\textrm{ along any ray}\}.
\end{align*}

\begin{definition}\
\begin{enumerate}
  \item An \emph{$H$-parabolic representation pair of type $\overrightarrow{\theta}$}\footnote{The local system version is defined in \cite[Definition 6.1]{bo} and \cite[Definition 4.1]{hs}.} is a pair $(\rho:\Gamma\rightarrow H,\mathbf{Q}=(Q_1,\cdots,Q_n))$, where each $Q_i$ is a parabolic subgroup of $H$, and $\overrightarrow{\theta}=(\theta_1,\cdots,\theta_n)\in (\mathfrak{t}_{\mathbb{R}})^n$, satisfying
      \begin{itemize}
        \item $\rho(\gamma_i)\in Q_i$ for each $i$,
        \item $Q_{i}$ is conjugate to $ Q_{\theta_i}$.
      \end{itemize}
  \item An $H$-parabolic representation pair $(\rho:\Gamma\rightarrow H,\mathbf{Q}=(Q_1,\cdots,Q_n))$ of type $\overrightarrow{\theta}$ is called \emph{semistable} (resp. \emph{stable}), if for each proper parabolic subgroup $H'$ of $H$ satisfying $\rho(\Gamma)\subset H'$ and for each nontrivial antidominant character $\chi'\in \mathrm{Hom}(\mathfrak{h}',\mathbb{C})$ of the Lie algebra $\mathfrak{h}'$ of $H'$, which is trivial on the centralizer of $\mathfrak{h}'$, we have
         \begin{align*}
         \deg_{\theta,\chi'}(\rho,\mathbf{Q},H')=\sum\limits_{i=1}^n \hat{\chi}'(\tilde\theta_i)\geq (\textrm{resp}. >) \ 0,
          \end{align*}
          where, via the root decomposition with respect to the maximal torus contained in $H'$, we extend the character $\chi'$ of $\mathfrak{h}'$ to a character $\hat \chi'\in \mathrm{Hom}(\mathfrak{h},\mathbb{C})$ of $\mathfrak{h}$ that is invariant under the action of the Weyl group $W(T)$ associated to $T$ when restricted to $\mathfrak{t}$, and $\tilde\theta_i=\mathrm{Ad}_{h_i}\theta_i$ for a unique element $h_i\in H$, up to the Weyl group $W(T)$, such that $Q_{\tilde\theta_i}=Q_i$.
\end{enumerate}
  \end{definition}

  The following proposition is the analogue of \cite[Lemma 2.1]{r} (see also \cite[Proposition 4.2]{la}).

\begin{proposition}
An $H$-parabolic representation pair $(\rho:\Gamma\rightarrow H,\mathbf{Q}=(Q_1,\cdots,Q_n))$ of type $\overrightarrow{\theta}$ is semistable if and only if
   \begin{align*}
     \sum\limits_{i}m(e^{i})e^{(i)}\leq 0,
   \end{align*}
where $e^{(i)}$ denotes an eigenvalue of $\tilde\theta_i$ under the adjoint action on $\mathfrak{h}$, with eigenvectors lying in $\mathfrak{h}'$, and $m(e^{i})$ is the corresponding multiplicity of $e^{(i)}$ in $\mathfrak{h}'$.
\end{proposition}

\subsection{Zariski tangent spaces and tangent quadratic cones}

Let $\mathbf{Art}_\mathbb{C}$ be the category of finite-dimensional local Artin $\mathbb{C}$-algebras that contain unity and have residue field $\mathbb{C}$, and let $\mathbf{Set}$ be the category of sets.
Fix a local Artin $\mathbb{C}$-algebra $A\in \mathbf{Art}_\mathbb{C}$. Then we have $\mathbb{C}$-algebra homomorphisms
\begin{align*}
 \mathbb{C}\xrightarrow{\iota}A\xrightarrow{q}\mathbb{C}
\end{align*}
where $q$ denotes the residue map. Let $G(A)$ be the Lie group of $A$-points whose Lie algebra is $\mathfrak{g}(A)=\mathfrak{g}\otimes_\mathbb{C} A$. Then we obtain Lie-group homomorphisms
\begin{align*}
G\xrightarrow{\iota}G(A)\xrightarrow{q}G.
\end{align*}
The kernel of $q$ is the unipotent group $G^0(A):=\exp{(\mathfrak{g}\otimes_{\mathbb{C}} \mathfrak{m}_A)}$.
We define
\begin{align*}
  \mathfrak{PRP}(r,\mathbf{d})(A)=\{( \rho_A: \Gamma\rightarrow G(A), {\mathbf{P}_A}=(P_{\mathcal{F}^{(1)}_A},\cdots,P_{\mathcal{F}^{(n)}_A})):  \rho_{A}(\gamma_i)\in P_{\mathcal{F}^{(i)}_A}, i=1,\cdots,n\},
\end{align*}
where $\{\mathcal{F}^{(1)}_A,\cdots,\mathcal{F}^{(n)}_A\}$ are  filtrations of $A$-modules on $V(A)=V\otimes_{\mathbb{C}}A$ of type $\mathbf{d}$ (that is, the $\ell$th-graded piece has dimension ${d^{(i)}}_\ell$ as an $A$-module), and $P_{\mathcal{F}^{(i)}_A}$ is a subgroup of $G(A)$ consisting of elements that preserve $\mathcal{F}^{(i)}_A$.
Given a flag $\mathcal{F}^{(i)}$ on $V$, let $\mathcal{F}^{(i)}(A)$ be the filtration on $V(A)$ given by $V^{(i)}_\ell(A)=V^{(i)}_\ell\otimes_\mathbb{C}A$. Then each $\mathcal{F}^{(i)}_A$ can be written as $\mathcal{F}^{(i)}_A=\tilde g \mathcal{F}^{(i)}(A)$ for some $\tilde g\in G(A)$. In particular, if $q(\mathcal{F}^{(i)}_A)=\mathcal{F}^{(i)}$, then one can take $\tilde g\in G^0(A)$.

Take $(\rho, \mathbf{P})\in \mathfrak{PRP}(r,\mathbf{d})$, and let $\widehat{\mathcal{O}}_{(\rho, \mathbf{P})}$ be the formal completion of the local ring to $\mathfrak{PRP}(r,\mathbf{d})$ at $(\rho, \mathbf{P})$, then we have the deformation functor
$$
\mathrm{ Def}_{(\rho, \mathbf{P})}: \mathbf{Art}_\mathbb{C}\rightarrow \mathbf{Set}
$$
defined by
\begin{align*}
\mathrm{ Def}_{(\rho, \mathbf{P})}(A)&:=\Hom(\widehat{\mathcal{O}}_{(\rho, \mathbf{P})},A)\\
&\ =\{( \rho_A, \mathbf{P}_A)\in \mathfrak{PRP}(r,\mathbf{d})(A): q\circ \rho_A=\rho,\ q(P_{\mathcal{F}^{(i)}_A})=P_{\mathcal{F}^{(i)}}, i=1,\cdots, n\}.
\end{align*}

We can also define \begin{align*}
                     &R(\Gamma,G;\mathbf{d})(A)\\
                     =&\ \{ \rho_A: \Gamma\rightarrow G(A):\textrm{ there are filtrations }\{\mathcal{F}^{(1)}_A,\cdots,\mathcal{F}^{(n)}_A\} \textrm{   of type  } \mathbf{d}\textrm{ such that  }   \rho_{A}(\gamma_i)\in P_{\mathcal{F}^{(i)}_A},i=1,\cdots,n\}
                   \end{align*}
and the deformation functor
$$
\mathrm{Def}^R_{\rho}: \mathbf{Art}_\mathbb{C}\rightarrow \mathbf{Set}
$$
by
\begin{align*}
 \mathrm{Def}^R_{\rho}(A)=\{ \rho_{A}\in R(\Gamma,G;\mathbf{d})(A):  q\circ \rho_{A}=\rho\}.
\end{align*}for a representation $\rho\in R(\Gamma,G;\mathbf{d})$.
Similarly, fixing $\mathbf{P}\in \mathrm{FL}(r,\mathbf{d})$,  we  define the deformation functor
$$
\mathrm{Def}^{\mathbf{P}}_{\rho}: \mathbf{Art}_\mathbb{C}\rightarrow \mathbf{Set}
$$
by
\begin{align*}
  \mathrm{Def}^{\mathbf{P}}_{\rho}(A)=\{ \rho_{A}: \Gamma\rightarrow G(A):  q\circ \rho_{A}=\rho,  \rho_{A}(\gamma_i)\in P_{\mathcal{F}^{(i)}}(A), i=1,\cdots, n\}.
\end{align*}

We can consider the following groupoids:
\begin{itemize}
  \item the transformation groupoid corresponding to $(\mathfrak{PRP}(r,\mathbf{d})(A), G(A))$ is denoted by $\mathbb{PRP}(r,\mathbf{d})(A)$;
  \item the transformation groupoid corresponding to $(\mathrm{ Def}_{(\rho, \mathbf{P})}(A), G^0(A))$ is denoted by $\mathbb{ DEF}_{(\rho, \mathbf{P})}(A)$;
  \item the transformation  groupoid corresponding to $(R(\Gamma,G;\mathbf{d})(A), G(A))$ is denoted by $\mathbb{R}(\Gamma,G;\mathbf{d})(A)$;
  \item the transformation  groupoid corresponding to $ (\mathrm{Def}^R_{\rho}(A), G^0(A))$ is denoted by $\mathbb{DEF}^R_{\rho}(A)$.
\end{itemize}

\begin{definition}
For a representation $\rho: \Xi\rightarrow G$ of a group $\Xi$ into $G$, we treat the Lie algebra $\mathfrak{g}$ of $G$ as a $\Xi$-module via the composition $\mathrm{Ad}_\rho$, obtained by composing $\rho$ with the adjoint representation $\mathrm{Ad}: G\to\mathrm{GL}(\mathfrak{g})$ of $G$. Let $C^p(\Xi,\mathfrak{g})$ ($p\geq 0$\footnote{One assigns $C^0(\Xi,\mathfrak{g})=\mathfrak{g}$.}) denote the abelian group of maps from $\Xi^p$  to $\mathfrak{g}$, and define the \emph{differential}  $\delta: C^p(\Xi,\mathfrak{g})\rightarrow C^{p+1}(\Xi,\mathfrak{g})$ by
\begin{align*}
  \delta X(\omega_1,\cdots,\omega_{p+1})= &\ \mathrm{Ad}_{\rho(\omega_1)}X(\omega_2,\cdots,\omega_{p+1})\\
  &+\sum\limits_{j=1}^p(-1)^jX(\omega_1,\cdots, \omega_{j-1},\omega_j\omega_{j+1},\omega_{j+2},\cdots,\omega_{p+1})\\
  &+(-1)^{p+1}X(\omega_1,\cdots,\omega_p)
\end{align*}
for $\omega_1,\cdots,\omega_{p+1}\in\Xi$. The groups of \emph{$p$-cocycles} and \emph{$p$-coboundaries} are defined by
\begin{align*}
  Z^p(\Xi,\mathfrak{g})&=\Ker (\delta: C^p(\Xi,\mathfrak{g})\rightarrow C^{p+1}(\Xi,\mathfrak{g}) ),\\
  B^p(\Xi,\mathfrak{g})&=\mathrm{Im }(\delta: C^{p-1}(\Xi,\mathfrak{g})\rightarrow C^{p}(\Xi,\mathfrak{g}) ),
\end{align*}
respectively. The \emph{$p$-th cohomology} associated to  $\rho: \Xi\rightarrow G$ is then defined as
\begin{align*}
  H^p(\Xi,\mathfrak{g})=\frac{ Z^p(\Xi,\mathfrak{g})}{B^p(\Xi,\mathfrak{g})}.
\end{align*}
\end{definition}

\begin{example}For our purpose of  studying a parabolic representation pair $(\rho, \mathbf{P})\in\mathfrak{PRP}(r,\mathbf{d})$, we consider
 $\Xi=\Gamma_i$, regarding $\mathfrak{g}$ as a $\Gamma_i$-module via the action of $\mathrm{Ad}_\rho$, then we have
\begin{align*}
  Z^1(\Gamma_i,\mathfrak{g})&\simeq\mathfrak{g},\\
  B^1(\Gamma_i,\mathfrak{g})&\simeq\mathrm{Im}(\mathrm{Ad}_{\rho(\gamma_i)}-\mathrm{Id}: \mathfrak{g\rightarrow\mathfrak{g}}).
\end{align*}
Hence, $H^1(\Gamma_i,\mathfrak{g})\simeq\mathrm{Coker}(\mathrm{Ad}_{\rho(\gamma_i)}-\mathrm{Id}: \mathfrak{g\rightarrow\mathfrak{g}})$. We also introduce the following notations:
\begin{itemize}
  \item  $ C^p(\Gamma_i,\mathfrak{g}, \mathfrak{p}_{\mathcal{F}^{(i)}})$ denotes the quotient space  $ C^p(\Gamma_i,\mathfrak{g})/C^p(\Gamma_i,\mathfrak{p}_{\mathcal{F}^{(i)}})$;
        \item $
        [\bullet]_i:C^p(\Gamma_i,\mathfrak{g})\rightarrow  C^p(\Gamma_i,\mathfrak{g},\mathfrak{p}_{\mathcal{F}^{(i)}})
        $
        denotes the natural projection;
        \item the linear map
       $ \delta: C^p(\Gamma_i,\mathfrak{g},\mathfrak{p}_{\mathcal{F}^{(i)}})\rightarrow C^{p+1}(\Gamma_i,\mathfrak{g},\mathfrak{p}_{\mathcal{F}^{(i)}})$
        is  induced from the differential $\delta: C^p(\Gamma_i,\mathfrak{g})\rightarrow  C^{p+1}(\Gamma_i,\mathfrak{g})$;
        \item  $H^p(\Gamma_i,\mathfrak{g},\mathfrak{p}_{\mathcal{F}^{(i)}})=\frac{Z^p(\Gamma_i,\mathfrak{g},\mathfrak{p}_{\mathcal{F}^{(i)}})}{B^p(\Gamma_i,\mathfrak{g},\mathfrak{p}_{\mathcal{F}^{(i)}})}$, where $Z^p(\Gamma_i,\mathfrak{g},\mathfrak{p}_{\mathcal{F}^{(i)}})$ and $B^p(\Gamma_i,\mathfrak{g},\mathfrak{p}_{\mathcal{F}^{(i)}})$ denote the kernels and images of $\delta$, respectively.
\end{itemize}
\end{example}

\begin{theorem}\label{o}
Fix a parabolic representation pair $(\rho, \mathbf{P})\in\mathfrak{PRP}(r,\mathbf{d})$.
\begin{enumerate}
  \item  The Zariski tangent space of $\mathfrak{PRP}(r,\mathbf{d})$ at $(\rho, \mathbf{P})$ is given by
      \begin{align*}
        \mathrm{ZT}_{(\rho, \mathbf{P})}\mathfrak{PRP}(r,\mathbf{d})=&\ \{(X, Y_1,\cdots, Y_n)\in Z^1(\Gamma,\mathfrak{g})\bigoplus (\bigoplus\limits_{i=1}^nC^0(\Gamma_i,\mathfrak{g},\mathfrak{p}_{\mathcal{F}^{(i)}}): \\
        &\ \ \ \ \ \ \ \ \ \ \ \ \ \ \ \ \ \ \ \ \ \ \ \ \ \ \ \ \ \ \ \ \ \  \ \ \ \ \ \ \ \ \ [X|_{\Gamma_i}]_i+ \delta_i  Y_i=0, i=1,\cdots,n \},
      \end{align*} and it has complex dimension  $r^2(2g+n-1)$.

  \item The tangent quadratic cone of $\mathfrak{PRP}(r,\mathbf{d})$ at $(\rho, \mathbf{P})$ is given by
  \begin{align*}
   &\mathrm{QC}_{(\rho, \mathbf{P})}\mathfrak{PRP}(r,\mathbf{d})\\
   =&\ \{ (X, Y_1,\cdots, Y_n)\in \mathrm{ZT}_{(\rho, \mathbf{P})}\mathfrak{PRP}(\mathbf{d},r):
    [X,X]_\rho=\delta V\textrm{ for } V\in C^1(\Gamma,\mathfrak{g}),\\
  & \ \ \ \ \ \ \ \ \ \ \ \ \ \ \ \ \ \ \ \ \ \ \ \ \ \ \ \ \ \ \ \ \ \ \      [Y_i, Y_i]_\rho-[V|_{\Gamma_i}]_i=\delta Z_i\textrm{ for   } Z_i\in \mathfrak{g}/(\mathfrak{p}_{\mathcal{F}^{(i)}}+[ Y_i,  \mathfrak{p}_{\mathcal{F}^{(i)}} ]), i=1,\cdots,n\},
  \end{align*}
  where
  \begin{align*}
   [X,X]_\rho(\omega_1,\omega_2)&=[X(\omega_1),\mathrm{Ad}_{\rho(\omega_1)}X(\omega_2)],\\
   [ Y_i, Y_i]_\rho(\gamma_i)&=[ Y_i,\mathrm{Ad}_{\rho(\gamma_i)} Y_i].
  \end{align*}
\end{enumerate}
\end{theorem}

\begin{proof}
(1)  Let $A_k$ be the truncated polynomial ring $\mathbb{C}[t]/(t^{k+1})$, and let $q^k_l:A_k\rightarrow A_l$ be the natural projection for $k\geq l$. The Zariski tangent space of $\mathfrak{PRP}(\mathbf{d},r)$ at $(\rho, \mathbf{P})$ is given by
\begin{align*}
 \mathrm{ZT}_{(\rho, \mathbf{P})}\mathfrak{PRP}(r,\mathbf{d})&= \mathrm{ Def}_{(\rho, \mathbf{P})}({A_1})\\
 &=\{( \rho_{A_1}: \Gamma\rightarrow G(A_1), \mathbf{P}_{A_1}=(P_{{\mathcal{F}^{(1)}_{A_1}}},\cdots,P_{{\mathcal{F}^{(n)}_{A_1}}})):\\
&\ \ \ \ \ \ \ \ \ \ \ \ \ \ \ \ \ \ \   \ \ \ \ \ \ \ \ \  \rho_{A_1}(\gamma_i)\in P_{{\mathcal{F}^{(i)}_{A_1}}}, q^1_0\circ\rho_{A_1}=\rho, q^1_0(P_{{\mathcal{F}^{(i)}_{A_1}}})=P_{\mathcal{F}^{(i)}}, i=1,\cdots, n\}.
\end{align*}
There is a semidirect product decomposition $G(A_1)=G^0(A_1)G$, where $G^0(A_1)=\exp(\mathfrak{g}\otimes_\mathbb{C}\mathfrak{m}_1)$ for the maximal ideal $\mathfrak{m}_1=(t)/(t^2)$  of $A_1$. Hence, there exist
$X\in C^1( \Gamma, \mathfrak{g})$ and $Y_1\in \mathfrak{g}/\mathfrak{p}_{\mathcal{F}^{(1)}},\cdots,Y_n\in\mathfrak{g}/\mathfrak{p}_{\mathcal{F}^{(n)}}$
such that
\begin{itemize}
  \item $\rho_{A_1}(\omega)=\exp(X(\omega)t) \rho(\omega)$ for any $\omega\in \Gamma$;
  \item $P_{{\mathcal{F}^{(i)}_{A_1}}}=\exp(\tilde Y_it)P_{\mathcal{F}^{(i)}(A_1)}\exp(-\tilde Y_it)$,
  \end{itemize}
where $\tilde Y_i$ denotes a lifting of $Y_i$ in  $C^0(\Gamma_i,\mathfrak{g})$.

One can check the following:
\begin{itemize}
  \item Since $\rho_{A_1}$ is a group homomorphism, we have
  \begin{align}\label{1}
   \exp(X(\omega_1)t) \rho(\omega_1)\exp(X(\omega_2)t) \rho^{-1}(\omega_1)=\exp(X(\omega_1\omega_2)t)
  \end{align}
  for any $\omega_1,\omega_2\in\Gamma$. Taking the derivative with respect to $t$ on both sides of \eqref{1} and then evaluating at $t=0$ yields
  \begin{align*}
  X(\omega_1\omega_2)=X(\omega_1)+\mathrm{Ad}_{\rho(\omega_1)}X(\omega_2),
  \end{align*}
thus $X\in Z^1(\Gamma,\mathfrak{g})$  \cite{l,gm}.
  \item The condition  $\rho_{A_1}(\gamma_i)\in P_{\mathcal{F}^{(i)}_{A_1}}$, namely
  \begin{align*}
    \exp(-Y_it) \exp(X(\gamma_i)t) \rho(\gamma_i)\exp(Y_it)\in P_{\mathcal{F}^{(i)}(A_1)}=P_{\mathcal{F}^{(i)}}(A_1),
  \end{align*}
 implies that
  \begin{align*}
   \frac{d}{dt}|_{t=0} \exp(-\tilde Y_it) \exp(X(\gamma_i)t) \rho(\gamma_i)\exp(\tilde Y_it)\equiv0\  ( \mathrm{mod}\  \mathfrak{p}_{\mathcal{F}^{(i)}}),
  \end{align*}
which gives rise to
\begin{align*}
  X(\gamma_i)\equiv(\mathrm{Id}-\mathrm{Ad}_{\rho(\gamma_i)})\tilde Y_i\  ( \mathrm{mod}\  \mathfrak{p}_{\mathcal{F}^{(i)}}).
\end{align*}
\end{itemize}

  Since there is an isomorphism
  \begin{align*}
   & \{(\Theta_1,\cdots,\Theta_n;Y_1,\cdots,Y_n)\in(\bigoplus\limits_{i=1}^nZ^1(\Gamma_i,\mathfrak{g}))\bigoplus(\bigoplus\limits_{i=1}^nC^0(\Gamma_i,\mathfrak{g},\mathfrak{p}_{\mathcal{F}^{(i)}})): [\Theta_i]+\delta_i Y_i=0, i=1,\cdots,n \}\\
  \simeq &\ (\bigoplus\limits_{i=1}^n   \mathfrak{p}_{\mathcal{F}^{(i)}})\bigoplus(\bigoplus\limits_{i=1}^n \mathrm{Im}(\delta_i))\bigoplus(\bigoplus\limits_{i=1}^n \mathrm{Ker}(\delta_i))\\
  \simeq &\ \mathfrak{g}^{\oplus n},
  \end{align*}
  we conclude that
  \begin{align*}
     \mathrm{ZT}_{(\rho, \mathbf{P})}\mathfrak{PRP}(r,\mathbf{d})\simeq Z^1(\Gamma,\mathfrak{g}).
  \end{align*}
  It is known that the complex dimension of $Z^1(\Gamma,\mathfrak{g})$ is $r^2(2g+n-1)$.

Thus, we complete the proof of (1).

(2) To determine the tangent quadratic cone of $\mathfrak{PRP}(r,\mathbf{d})$ at $(\rho, \mathbf{P})$, we consider
\begin{align*}
 \mathrm{ Def}_{(\rho, \mathbf{P})}(A_2) &=\{( \rho_{A_2}: \Gamma\rightarrow G(A_2), \mathbf{P}_{A_2}=(P_{{\mathcal{F}^{(1)}_{A_2}}},\cdots,P_{{\mathcal{F}^{(n)}_{A_2}}})):\\
&\ \ \ \ \ \ \ \ \ \ \ \ \ \ \ \ \ \ \   \ \ \ \ \ \ \ \ \   \rho_{A_2}(\gamma_i)\in P_{{\mathcal{F}^{(i)}_{A_2}}}, q^2_0\circ\rho_{A_2}=\rho, q^2_0(P_{{\mathcal{F}^{(i)}_{A_2}}})=P_{\mathcal{F}^{(i)}}, i=1,\cdots, n\}
\end{align*}
Similarly, note that  $G(A_2)=\exp(\mathfrak{g}\otimes_\mathbb{C}\mathfrak{m}_2)\cdot G$, where $\mathfrak{m}_2=(t)/(t^3)$ is the maximal ideal of $A_2$, and
\begin{align*}
 \exp{(at+bt^2)}\exp{(ct+dt^2)}=\exp{((a+c)t+(b+d+\frac{1}{2}[a,c])t^2)}
\end{align*}
for $at+bt^2, ct+dt^2\in \mathfrak{g}\otimes_\mathbb{C}\mathfrak{m}_2$. Hence, there exist
$X,V\in C^1( \Gamma, \mathfrak{g})$ and $Y_1\in \mathfrak{g}/\mathfrak{p}_{\mathcal{F}^{(1)}}, Z_1\in \mathfrak{g}/(\mathfrak{p}_{\mathcal{F}^{(1)}}+[\tilde{Y}_1, \mathfrak{p}_{\mathcal{F}^{(1)}}]),\cdots,Y_n\in \mathfrak{g}/\mathfrak{p}_{\mathcal{F}^{(n)}}, Z_n\in \mathfrak{g}/(\mathfrak{p}_{\mathcal{F}^{(n)}}+[\tilde{Y}_n, \mathfrak{p}_{\mathcal{F}^{(n)}}])$
such that
  \begin{itemize}
  \item $\rho_{A_2}(\omega)=\exp(X(\omega)t+V(\omega)t^2) \rho(\omega)$ for any $\omega\in \Gamma$;
  \item $P_{\mathcal{F}^{(i)}_{A_2}}=\exp(\tilde Y_it+\tilde Z_it^2)P_{\mathcal{F}^{(i)}(A_2)}\exp(-\tilde Y_it-\tilde Z_it^2)$,
  \end{itemize}
  where $\tilde Z_i$ denotes a lifting of $Y_i$ in  $C^0(\Gamma_i,\mathfrak{g})$.

  Also, we can check the following:
  \begin{itemize}
    \item Since $\rho_{A_2}$ is a group homomorphism, we have
    \begin{align}\label{2}
   &\exp(X(\omega_1)t+V(\omega_1)t^2) \rho(\omega_1)\exp(X(\omega_2)t+V(\omega_2)t^2) \rho^{-1}(\omega_1) \nonumber\\
   =&\ \exp(X(\omega_1\omega_2)t+V(\omega_1\omega_2)t^2).
  \end{align}
    Taking values at $t=0$ of the first and second derivatives with respect to $t$ on both sides of \eqref{2} yields, respectively,
    \begin{align*}
  X(\omega_1\omega_2)&=X(\omega_1)+\mathrm{Ad}_{\rho(\omega_1)}X(\omega_2),\\
  V(\omega_1\omega_2)&=V(\omega_1)+\mathrm{Ad}_{\rho(\omega_1)}V(\omega_2)+\frac{1}{2}[X(\omega_1),\mathrm{Ad}_{\rho(\omega_1)}X(\omega_2)],
\end{align*}
where we use the identities
\begin{align*}
\frac{d}{dt}\bigg|_{t=0} \rho(\omega_1)\exp(X(\omega_2)t+V(\omega_2)t^2)\rho^{-1}(\omega_1)&=\mathrm{Ad}_{\rho(\omega_1)}X(\omega_2),\\
 \frac{d^2}{dt^2}\bigg|_{t=0} \rho(\omega_1)\exp(X(\omega_2)t+V(\omega_2)t^2)\rho^{-1}(\omega_1)&=(\mathrm{Ad}_{\rho(\omega_1)}X(\omega_2))^2+2\mathrm{Ad}_{\rho(\omega_1)}V(\omega_2).
\end{align*}
    \item The condition  $\rho_{A_2}(\gamma_i)\in P_{\mathcal{F}^{(i)}_{A_2}}$ implies
    \begin{align*}
     \frac{d}{dt}\bigg|_{t=0}\exp(-\tilde Y_it-\tilde Z_it^2)\exp(X(\gamma_i)t+V(\gamma_i)t^2)\rho(\gamma_i)\exp(\tilde Y_it+\tilde Z_it^2)
     &\equiv0 \ ( \mathrm{mod}\  \mathfrak{p}_{\mathcal{F}^{(i)}} ),\\
      \frac{d^2}{dt^2}\bigg|_{t=0}\exp(-\tilde Y_it-\tilde Z_it^2)\exp(X(\gamma_i)t+V(\gamma_i)t^2)\rho(\gamma_i)\exp(\tilde Y_it+\tilde Z_it^2)&\equiv0\  ( \mathrm{mod}\  \mathfrak{p}_{\mathcal{F}^{(i)}}).
    \end{align*}
    It follows that
    \begin{align*}
   X(\gamma_i)&\equiv(\mathrm{Id}-\mathrm{Ad}_{\rho(\gamma_i)})\tilde Y_i \  ( \mathrm{mod}\  \mathfrak{p}_{\mathcal{F}^{(i)}}),\\
    V(\gamma_i)&\equiv(\mathrm{Id}-\mathrm{Ad}_{\rho(\gamma_i)})\tilde Z_i-\frac{1}{2}[\tilde Y_i,\mathrm{Ad}_{\rho(\gamma_i)}\tilde Y_i] \  ( \mathrm{mod}\   \mathfrak{p}_{\mathcal{F}^{(i)}}).
  \end{align*}
  \end{itemize}

  We thus complete the proof of (2).
 \end{proof}

\begin{corollary}\label{co}\
\begin{enumerate}
  \item  Given $\mathbf{P}\in \mathrm{FL}(r,\mathbf{d})$, we have:
  \begin{enumerate}
    \item The Zariski tangent space of $R(\Gamma,G;\mathbf{P})$ at $\rho$ is given by
\begin{align*}
 ZT_{(\rho, \mathbf{P})}R(\Gamma,G;\mathbf{P})
& = \mathrm{Def}^{\mathbf{P}}_{\rho}(A_1)\\
&=\{ \rho_{A_1}: \Gamma\rightarrow G(A_1):  q\circ \rho_{A_1}=\rho,  \rho_{A_1}(\gamma_i)\in P_{\mathcal{F}^{(i)}}(A_1), i=1,\cdots, n\}\\
 &= \{X\in Z^1(\Gamma,\mathfrak{g}):X(\gamma_i)\in\mathfrak{p}_{\mathcal{F}^{(i)}}, i=1,\cdots,n\},
\end{align*}
and it has complex dimension $r^2(2g-1)+\sum\limits_{i=1}^nf_i$, where
$f_i=\frac{1}{2}(r^2+\sum\limits_{\ell=0}^{s_i}(d^{(i)}_\ell)^2)$.
    \item The tangent quadratic cone of $R(\Gamma,G;\mathbf{P})$ at $\rho$  is
 \begin{align*}
   \mathrm{QC}_{(\rho, \mathbf{P})}R(\Gamma,G;\mathbf{P})&=\{ X\in \mathrm{ZT}_{(\rho, \mathbf{P})}R(\Gamma,G;\mathbf{P}):   [X,X]_\rho=\delta V\\
  & \ \ \ \ \ \ \ \ \ \  \ \ \ \ \ \ \ \ \ \ \ \ \textrm{ for } V\in C^1(\Gamma,\mathfrak{g}) \textrm{ satisfying } V(\gamma_i)\in \mathfrak{p}_{\mathcal{F}^{(i)}}, i=1,\cdots,n\}.
  \end{align*}
  \end{enumerate}

  \item Fix $\rho\in R(\Gamma,G;\mathbf{d})$, and choose $(P_{\mathcal{F}^{(1)}},\cdots,P_{\mathcal{F}^{(n)}})\in\mathrm{FL}(r,\mathbf{d})$ such that $\rho(\gamma_i)\in P_{\mathcal{F}^{(i)}} $, then:
  \begin{enumerate}
    \item The Zariski tangent space of $R(\Pi,G;\mathbf{d})$ at $\rho$ is given by
\begin{align*}
 ZT_{\rho}R(\Gamma,G;\mathbf{d})
\simeq\Ker(\Lambda: Z^1(\Gamma,\mathfrak{g})\rightarrow\bigoplus\limits_{i=1}^nH^1(\Gamma_i,\mathfrak{g},\mathfrak{p}_{\mathcal{F}^{(i)}})),
\end{align*} where  $\Lambda$ is the composition
 \begin{align*}
   Z^1(\Gamma,\mathfrak{g})\rightarrow \bigoplus\limits_{i=1}^nZ^1(\Gamma_i,\mathfrak{g})\rightarrow\bigoplus\limits_{i=1}^nZ^1(\Gamma_i,\mathfrak{g},\mathfrak{p}_{\mathcal{F}^{(i)}})\rightarrow \bigoplus\limits_{i=1}^nH^1(\Gamma_i,\mathfrak{g},\mathfrak{p}_{\mathcal{F}^{(i)}}),
 \end{align*}
 consisting of natural restriction and projections;
    \item the tangent quadratic cone of $R(\Gamma,G;\mathbf{d})$ at $\rho$ is given by
  \begin{align*}
   &\mathrm{QC}_{\rho}R(\Gamma,G;\mathbf{d})\\
   =&\ \{ X\in ZT_{\rho}R(\Gamma,G;\mathbf{d}):
    [X,X]_\rho=\delta V, \textrm{ where } V\in C^1(\Gamma,\mathfrak{g}) \textrm{ satisfying  }  [V|_{\Gamma_i}]_i=[Y_i, Y_i]_\rho+\delta Z_i\\
 &\ \ \ \ \ \ \ \ \ \ \ \ \ \ \ \ \ \textrm{ for } Y_i\in C^0(\Gamma_i,\mathfrak{g},\mathfrak{p}_{\mathcal{F}^{(i)}}), Z_i\in \mathfrak{g}/( \mathfrak{p}_{\mathcal{F}^{(i)}}+[ Y_i,  \mathfrak{p}_{\mathcal{F}^{(i)}} ]) , i=1,\cdots,n\}.
  \end{align*}
  \end{enumerate}
\end{enumerate}

\end{corollary}

\section{Deformations of parabolic logarithmic flat bundles}

\subsection{Parabolic-groupoid version of the Riemann--Hilbert--Deligne correspondence}

Now we turn to parabolic logarithmic flat bundles $(E,\mathbf{F}, \nabla)$ of rank $r$ and type $\mathbf{d}$, where
\begin{itemize}
    \item $E$ is a vector bundle of rank $r$ over $\bar X$;
    \item the parabolic structure $\mathbf{F}=\{\mathcal{F}^{(i)}\}_{i=1,\cdots,n}$ consists of the flags
$$
\mathcal{F}^{(i)}: \{0\}=E^{(i)}_{s_i+1}\subsetneq E^{(i)}_{s_i}\subsetneq\cdots\subsetneq E^{(i)}_{0}=E|_{x_i}
$$
of type $\overrightarrow{d^{(i)}}$;
\item $\nabla: E\rightarrow E\otimes_{\mathcal{O}_{\bar X}}\Omega^1_{\bar X}(\mathcal{D})$ is a logarithmic connection with poles along the divisor $\mathcal{D}=x_1+\cdots+x_n$  compatible with $\mathbf{F}$, that is, the residue  $\mathrm{Res}(\nabla, x_i)$ of $\nabla$ at $x_i$ preserves the flag  $\mathcal{F}^{(i)}$.
\end{itemize}
Let $\End(E,\mathbf{F})$ be the sheaf of germs of
endomorphisms of $E$ which preserve the   flag $\mathcal{F}^{(i)}$ at  each $x_i$, namely it is characterized by the short exact sequence
\begin{align*}
  0\rightarrow\End(E,\mathbf{F})\rightarrow\End(E)\rightarrow\bigoplus\limits_{i=1}^nS_i\rightarrow0,
\end{align*}
where $S_i$ is the skyscraper sheaf supported at $x_i$ with fiber $S_i|_{x_i}=\mathfrak{g}/\mathfrak{p}_{\mathcal{F}^{(i)}}$.  Obviously, $\End(E,\mathbf{F})$ is torsion-free, hence is locally free.  Let  $\Aut(E,\mathbf{F})$ be the group of automorphisms of $E$ that preserve $\mathbf{F}$.

The set/category of parabolic logarithmic flat bundles of rank $r$ and type $\mathbf{d}$ is denoted by $\mathfrak{PFB}(r,\mathbf{d})$. It is well known that there is an equivalence between the categories $\mathfrak{PRP}(r,\mathbf{d})$ and $\mathfrak{PFB}(r,\mathbf{d})$ given by the following functors \begin{align*}
\mathrm{RHD}: \mathfrak{PRP}(r,\mathbf{d})&\longrightarrow \mathfrak{PFB}(r,\mathbf{d}),\\
\mathrm{RHD}^{-1}: \mathfrak{PFB}(r,\mathbf{d})&\longrightarrow \mathfrak{PRP}(r,\mathbf{d}),
\end{align*}
 such equivalence is called the Riemann--Hilbert--Deligne correspondence. We describe it respectively as follows.
\begin{itemize}
  \item Given a parabolic representation pair  $(\rho, \mathbf{P})$, or equivalently a filtered local system $(\mathbb{L}_\rho, \mathbf{F}^\rho)$ of rank $r$ and type $\mathbf{d}$, let $(\mathcal{E}^\rho,\nabla^\rho)$ be the flat bundle over $X$ corresponding to the local system $\mathbb{L}_\rho$, and $(\mathcal{E}^{(i)}_\ell, \nabla^{(i)}_\ell)$ be the flat bundles over a small punctured neighborhood $U_i$ of $x_i $, corresponding to the sub-local systems $ {\mathbb{L}^\rho}^{(i)}_\ell$. One chooses a stalk $V^{(i)}_\ell$ of  $ {\mathbb{L}^\rho}^{(i)}_\ell$ around $x_i$. Let $\{{e^{(i)}_\ell}_\mu\}$ be the distinct eigenvalues of the $\rho(\gamma_i)$-action on $V^{(i)}_\ell$, with arguments chosen in $[0,2\pi)$. Define the generalized eigenspace
      \begin{align*}
        B_{{e^{(i)}_\ell}_\mu}=\{v\in V^{(i)}_\ell: (\rho(\gamma_i)-{e^{(i)}_\ell}_\mu)^{v^{(i)}_\ell}v=0\},
      \end{align*}for $v^{(i)}_\ell=\dim_{\mathbb{C}}V^{(i)}_\ell$
      and the associated nilpotent operator
      \begin{align*}
        N_{{e^{(i)}_\ell}_\mu}=\log\{({e^{(i)}_\ell}_\mu)^{-1}\rho(\gamma_i)|_{ B_{{e^{(i)}_\ell}_\mu}}\}.
      \end{align*}
     For $v\in  B_{{e^{(i)}_\ell}_\mu}$, which can be viewed as a multivalued horizontal section of $(\mathcal{E}^{(i)}_\ell, \nabla^{(i)}_\ell)$, then all the sections of $(\mathcal{E}^{(i)}_\ell, \nabla^{(i)}_\ell)$ are as follows
      \begin{align*}
        \tilde v=\exp\{\frac{1}{2\pi\sqrt{-1}}\log z^{(i)}(\log {e^{(i)}_\ell}_\mu+ N_{{e^{(i)}_\ell}_\mu})\}v,
      \end{align*}
      where $z^{(i)}$ is the local coordinate on $U_i$, generate an extension $\bar {\mathcal{E}}^{(i)}_\ell$ of $\mathcal{E}^{(i)}_\ell$ to the neighborhood $\bar U_i=U_i\bigcup\{x_i\}$ of $x_i$. Meanwhile, $\nabla^{(i)}_t$ extends to a logarithmic connection
      $$
      \bar{\nabla}^{(i)}_\ell: \bar {\mathcal{E}}^{(i)}_\ell\rightarrow
      \bar {\mathcal{E}}^{(i)}_\ell\otimes_{\mathcal{O}_{\bar U_i}}\Omega^1_{\bar U_i}(x_i).
      $$
      In particular, $\{(\bar {\mathcal{E}}^{(i)}_{s_i},\bar{\nabla}^{(i)}_{s_i})\}_{i=1,\cdots,n}$ produces the Deligne canonical extension $(\bar{\mathcal{E}}^\rho, \bar\nabla^\rho )$ of $(\mathcal{E}^\rho,\nabla^\rho)$. The filtration
       \begin{align*}
\mathcal{D}^{(i)}: \{0\}= \bar {\mathcal{E}}^{(i)}_{s_i+1}|_{x_i}\subsetneq \bar {\mathcal{E}}^{(i)}_{s_i}|_{x_i}\subsetneq\cdots\subsetneq\bar {\mathcal{E}}^{(i)}_{1}|_{x_i}\subsetneq \bar {\mathcal{E}}^{(i)}_{0}|_{x_i}=\bar{\mathcal{E}}^\rho|_{x_i}
  \end{align*}
  is preserved by the residue of $\bar \nabla$ at $x_i$. Therefore, from $(\rho, \mathbf{P})\in\mathfrak{PRP}(r,\mathbf{d})$, we obtain a parabolic logarithmic flat bundle
  \begin{align*}
   (\bar{\mathcal{E}}^\rho, \bar\nabla^\rho, \mathbf{F}=\{\mathcal{D}^{(i)}\}_{i=1,\cdots,n})\in\mathfrak{PFB}(r,\mathbf{d})
  \end{align*}
   of rank $r$ and type $\mathbf{d}$ over $\bar X$ \footnote{The parabolic structure given above is generally coarser than that in \cite{si}, since we do not impose the constraint of zero parabolic degree with respect to some weight system on the extended parabolic bundle.}. Note that the real parts  of the eigenvalues of $ \mathrm{Res}(\bar\nabla^\rho,x_i)$ lie in $[0,1)$, and  the degree of $\bar{\mathcal{E}}^\rho$ is  $-\sum a$, where the sum is taken over all such real parts $a$ of  eigenvalues of the residues of  $ \bar\nabla^\rho$ at the punctures.
       \item Conversely, given a parabolic logarithmic flat bundle $(E,\mathbf{F}, \nabla)$ of rank $r$ and type $\mathbf{d}$, let $\mathbb{L}^{(E,\nabla)}$ be the local system over $X$ associated to the flat bundle $(E,\nabla)|_X$. There is a trivialization of $(E,\nabla)|_{U_i}$ over some punctured neighborhood $U_i$ of $x_i$ of the form
       $$
       (E|_{x_i}\times U_i, d+\frac{Q_i}{z^{(i)}}dz^{(i)}), \quad Q_i\in \mathfrak{p}_{\mathcal{F}^{(i)}},
       $$
       so that $(V^{(i)}_\ell\times U_i, d+\frac{Q_i}{z^{(i)}}dz^{(i)})$ for the subspace $V^{(i)}_\ell$ of $E|_{x_i}$ in the flag $\mathcal{F}^{(i)}\in \mathbf{F}$ defines a sub-local system $\mathbb{L}^{(i)}_\ell$ of $\mathbb{L}^{(E,\nabla)}$. Thus we obtain a filtration
           \begin{align*}
             \mathbb{D}_{U_i}:  0=\mathbb{L}^{(i)}_{s_i+1}\subsetneq \mathbb{L}^{(i)}_{s_i}\subsetneq\cdots\subsetneq \mathbb{L}^{(i)}_{1}\subsetneq \mathbb{L}^{(i)}_{0}=\mathbb{L}^{(E,\nabla)}|_{U_i}.
           \end{align*}
           Consequently, from $(E,\mathbf{F}, \nabla)\in \mathfrak{PFB}(r,\mathbf{d})  $, we obtain a filtered local system
           \begin{align*}
             (\mathbb{L}^{(E,\nabla)},\mathbf{D}=\{\mathbb{D}_{U_i}\}_{i=1,\cdots,n} )\in \mathfrak{ FLS}(r,\mathbf{d})
           \end{align*}
          of rank $r$ and type $\mathbf{d}$ over $X$, hence an object in $\mathfrak{ PRP}(r,\mathbf{d})$.
\end{itemize}

 Fix a parabolic bundle $(E,\mathbf{F})$, and let $\mathfrak{LC}_{(E,\mathbf{F})}$ be the space of  logarithmic connections on $(E,\mathbf{F})$ that are compatible with $\mathbf{F}$. The transformation groupoid corresponding to $(\mathfrak{LC}_{(E,\mathbf{F})}, \Aut(E,\mathbf{F}))$ is denoted by $\mathbb{LC}_{(E,\mathbf{F})}$. The evaluation map yields a group homomorphism $\varepsilon:\Aut(E,\mathbf{F})\rightarrow G$ such that the following commutative diagram
\begin{align*}
   \CD
     \mathfrak{LC}_{(E,\mathbf{F})}@>\mathrm{RHD}^{-1}>>\mathfrak{PRP}(r,\mathbf{d}) \\
     @V \Phi VV @V \varepsilon(\Phi) VV  \\
     \mathfrak{LC}_{(E,\mathbf{F})} @>\mathrm{RHD}^{-1}>> \mathfrak{PRP}(r,\mathbf{d})
   \endCD
  \end{align*}
for any $\Phi\in \Aut(E,\mathbf{F})$.  The groupoid  $\mathbb{PFB}(r,\mathbf{d})$ is defined as the disjoint union, over the set $\mathrm{Iso}\mathfrak{PB}(r,\mathbf{d})$ of isomorphism classes of parabolic  bundles of rank $r$ and type $\mathbf{d}$ over $\bar X$, of the groupoids $\mathbb{LC}_{(E,\mathbf{F})}$ for $(E,\mathbf{F})\in \mathfrak{PFB}(r,\mathbf{d})$.

 \begin{definition}\
 \begin{enumerate}
    \item We say that two parabolic logarithmic flat bundles $(E,\mathbf{F},\nabla), (E',\mathbf{F}',\nabla')\in \mathfrak{PRP}(r,\mathbf{d})$ are \emph{RHD-related} if
    \begin{itemize}
   \item  $(E,\mathbf{F},\nabla)$ is not isomorphic to $(E',\mathbf{F}',\nabla')$;
   \item  they produce the same filtered local system, i.e. $\mathrm{RHD}^{-1}(E,\mathbf{F},\nabla)=\mathrm{RHD}^{-1}(E',\mathbf{F}',\nabla')$.
   \end{itemize}

        \item We define the  groupoid  $\mathbb{PFB}^{[0,1)}(r,\mathbf{d})$ to have the same objects as $\mathbb{PFB}(r,\mathbf{d})$, and morphisms consisting of those in $\mathbb{PFB}(r,\mathbf{d})$ together with the RHD-relatedness.
\end{enumerate}
\end{definition}

\begin{proposition} If two parabolic logarithmic flat bundles $(E,\mathbf{F},\nabla), (E',\mathbf{F}',\nabla')\in \mathfrak{PRP}(r,\mathbf{d})$ are RHD-related, then
 for each $x_i$, there is a neighborhood $\bar U_i$ of $x_i$ such that
 \begin{itemize}
   \item $(E,\nabla)|_{\bar U_i}=\bigoplus\limits_{\alpha}(E_\alpha,\nabla_\alpha)$, where $E_\alpha|_{x_i}$ is the generalized eigenspace of $\mathrm{Res}(\nabla,x_i)$, and $\nabla_\alpha=\nabla|_{E_\alpha}$;
   \item $(E',\nabla')|_{\bar U_i}=\bigoplus\limits_{\alpha}(E'_\alpha,\nabla'_\alpha)$, where
   \begin{itemize}
       \item $E'_\alpha=E_\alpha\otimes \mathcal{O}_{\bar X}(n^{(i)}_\alpha x_i)|_{\bar U_i}$ for some integer $n^{(i)}_\alpha$,
       \item $\nabla'_\alpha=\nabla_\alpha\otimes D_\alpha^{(i)}$ with the connection $ D_\alpha^{(i)}=d-\frac{n^{(i)}_\alpha}{z^{(i)}}dz^{(i)}$ on $ \mathcal{O}_{\bar X}(n^{(i)}_\alpha x_i)|_{\bar U_i}$;
   \end{itemize}
   \item $\mathbf{F}'\bigcap E'_\alpha|_{x_i}=(\mathbf{F}\bigcap E_\alpha|_{x_i})\otimes  \mathcal{O}_{\bar X}(n^{(i)}_\alpha x_i)|_{x_i}$.
 \end{itemize}\end{proposition}
\begin{remark}  The  groupoid $\mathbb{PFB}^{[0,1)}(r,\mathbf{d})$ is equivalent to the transformation groupoid of $\mathbb{R}$-filtered logarithmic flat bundles ($\mathbb{R}$-filtered regular $D_X$-modules) of rank $r$ and type $\mathbf{d}$ introduced by Simpson \cite{si}.
 \end{remark}

Given a local Artin $\mathbb{C}$-algebra $A\in\mathbf{Art}_\mathbb{C}$,  let $\bar X_A=\bar X\times \mathrm{Spec}(A)$, then we have the following commutative diagram
\begin{align*}
   \CD
    \bar X@>\iota>>\bar X_A @>q>>\bar X\\
     @V VV @V p VV   @V VV \\
     \mathrm{Spec}(\mathbb{C}) @>\iota>> \mathrm{Spec}(A)@>q>> \mathrm{Spec}(\mathbb{C})
   \endCD.
  \end{align*}
  Hence $\bar X_A$ can be interpreted as the topological space $\bar X$ equipped with the structure sheaf $\mathcal{O}_{\bar X}(A)=\mathcal{O}_{\bar X}\otimes_{\mathbb{C}}A$.

\begin{definition}
A \emph{parabolic logarithmic flat  $A$-bundle} is a triple  $( E_A,  {\mathbf{F}_A},\nabla_A)$, where
 \begin{itemize}
   \item  $E_A$ is a locally free sheaf of  $\mathcal{O}_{\bar X}(A)$-modules over $\bar X_A$;
   \item
  $\mathbf{F}_A=\{\mathcal{F}^{(i)}_A\}_{i=1,\cdots,n}$ consists of  filtrations $\mathcal{F}^{(i)}_A$ of $A$-modules on $(\iota^*E_A)|_{x_i}$;
   \item $\nabla_A: E_A\rightarrow  E_A\otimes_{\mathcal{O}_{\bar X_A}}q^*(\Omega^1_{\bar X}(\mathcal{D}))$ is a logarithmic connection on $E_A$ compatible with $\mathbf{F}_A$, namely, the residue $\mathrm{Res}(\iota^*\nabla_A, x_i)\in \mathfrak{g}(A)$ preserves $\mathcal{F}^{(i)}_A$ for each $i$.
 \end{itemize}
\end{definition}

The set/category of parabolic  $A$-bundles of rank $r$ and type $\mathbf{d}$ is denoted by $\mathfrak{PB}(r,\mathbf{d})(A)$, and the set/category of parabolic logarithmic flat $A$-bundles of rank $r$ and type $\mathbf{d}$ is denoted by $\mathfrak{PFB}(r,\mathbf{d})(A)$.
Given a parabolic logarithmic flat $A$-bundle $( E_A,  {\mathbf{F}_A},\nabla_A)$, we have a parabolic logarithmic flat bundle $(q^\sharp (E_A), q(\mathbf{F}_A), q^\sharp(\nabla_A))$, where $q^\sharp=q\circ\iota^*$. Hence, there are functors
\begin{align*}
   q: \mathfrak{PB}(r,\mathbf{d})(A)&\longrightarrow \mathfrak{PB}(r,\mathbf{d}),\\
  q: \mathfrak{PFB}(r,\mathbf{d})(A)&\longrightarrow \mathfrak{PFB}(r,\mathbf{d}) .
\end{align*}
We also introduce the following sets/categories:
\begin{itemize}
  \item given a parabolic bundle $(E,\mathbf{F})$, the  set/category
$\mathfrak{PB}_{(E,\mathbf{F})}(A)$ consists of parabolic  $A$-bundles $(E_A,\mathbf{F}_A)$  of rank $r$ and type $\mathbf{d}$ such that $q^\sharp(E_A)=E$ and $q({\mathbf{F}_A})=\mathbf{F}$;
  \item given a parabolic  $A$-bundle $( E_A,{\mathbf{F}_A})$, the  set/category
$\mathfrak{LC}_{(E_A,{\mathbf{F}_A})}$ consists  of logarithmic connections $ \nabla_A$ on $ E_A$ compatible with $\mathbf{F}_A$;
\item given  a parabolic logarithmic flat bundle $(E,\mathbf{F},\nabla)$ and a parabolic $A$-bundle $( E_A,{\mathbf{F}_A})\in\mathfrak{PB}_{(E,\mathbf{F})}(A)$, the set/category
$\mathfrak{LC}_{(E, \mathbf{F},\nabla, E_A,{\mathbf{F}_A})}$ consists of  logarithmic  connections $\nabla_A$ on $E_A$ compatible with $\mathbf{F}_A$ such that  $q^\sharp( \nabla_A)=\nabla$;
\item given a parabolic logarithmic flat bundle $( E, {\mathbf{F}}, \nabla)$, the set/category
$\mathfrak{PFB}_{( E, {\mathbf{F}}, \nabla)}(A)$ consists of parabolic logarithmic flat $A$-bundles $( E_A,\mathbf{F}_A, \nabla_A)$  such that $q^\sharp( E_A)=E$, $q({\mathbf{F}_A})= {\mathbf{F}}$, and $q^\sharp( \nabla_A)=\nabla$.
\end{itemize}
\begin{proposition}We have  functors
\begin{align*}
   \mathrm{RHD}_A: \mathfrak{PRP}(r,\mathbf{d})(A)&\longrightarrow \mathfrak{PFB}(r,\mathbf{d})(A), \\\mathrm{RHD}_A^{-1}: \mathfrak{PFB}(r,\mathbf{d})(A)&\longrightarrow \mathfrak{PRP}(r,\mathbf{d})(A),
\end{align*}
satisfying
\begin{enumerate}
  \item  given a parabolic $A$-bundle $(E_A,\mathbf{F}_A)\in \mathfrak{PB}(r,\mathbf{d})(A)$, there is a group homomorphism
\begin{align*}
\varepsilon_A: \Aut(E_A,\mathbf{F}_A)\longrightarrow  G(A)
\end{align*} such that the following diagram
\begin{align*}
  \xymatrix{
 \mathfrak{LC}_{q(E_A,\mathbf{F}_A)}\ar[d]^{q^\sharp(\Phi_A)}\ar@(ur,ul)[rrr]^{\mathrm{RHD}^{-1}}&\mathfrak{LC}_{( E_A,\mathbf{F}_A)} \ar[d]^{\Phi_A}\ar[l]_-{q} \ar[r]^-{\mathrm{RHD}^{-1}_A} & \mathfrak{PRP}(r,\mathbf{d})(A)\ar[d]^{\varepsilon_A(\Phi_A)} \ar[r]^-{q}&\mathfrak{PRP}(r,\mathbf{d})\ar[d]^{q(\varepsilon_A(\Phi_A))=\varepsilon(q^\sharp(\Phi_A))}\\
 \mathfrak{LC}_{q(E_A,\mathbf{F}_A)}\ar@(dr,dl)[rrr]^{\mathrm{RHD}^{-1}}&\mathfrak{LC}_{( E_A,\mathbf{F}_A)} \ar[l]_-{q} \ar[r]^-{\mathrm{RHD}^{-1}_A} & \mathfrak{PRP}(r,\mathbf{d})(A) \ar[r]^-{q}&\mathfrak{PRP}(r,\mathbf{d})}
\end{align*}
commutes for any $\Phi_A\in \Aut(E_A,\mathbf{F}_A)$, where $\Aut(E_A,\mathbf{F}_A)$  is  the group of automorphisms of $E_A$ that preserve $\mathbf{F}_A$;
  \item given a parabolic logarithmic flat bundle $(E,\mathbf{F},\nabla)$ with corresponding parabolic  representation pair $(\rho, \mathbf{P})$   and a parabolic $A$-bundle $( E_A,{\mathbf{F}_A})\in\mathfrak{PB}_{(E,\mathbf{F})}(r,\mathbf{d})(A)$, there is a group homomorphism
\begin{align*}
\varepsilon_A: \Aut_0(E_A,\mathbf{F}_A)\longrightarrow  G^0(A),
\end{align*} such that  the following  diagram
\begin{align*}
   \CD
    \mathfrak{LC}_{(E, \mathbf{F},\nabla, E_A,{\mathbf{F}_A})}@>\mathrm{RHD}^{-1}_A>>\mathrm{ Def}_{(\rho, \mathbf{P})}(A) \\
     @V \Phi_A VV @V \varepsilon_A(\Phi_A) VV  \\
 \mathfrak{LC}_{(E, \mathbf{F},\nabla, E_A,{\mathbf{F}_A})}@>\mathrm{RHD}^{-1}_A>> \mathrm{ Def}_{(\rho, \mathbf{P})}(A)
   \endCD
  \end{align*}
commutes for any $\Phi_A\in \Aut_0(E_A,\mathbf{F}_A)$, where $$\Aut_0(E_A,\mathbf{F}_A)=\{\Phi_A\in  \Aut(E_A,\mathbf{F}_A):q^\sharp(\Phi_A)=\mathrm{Id}\}.$$
\end{enumerate}

\end{proposition}
\begin{proof} We only need to show the canonical Deligne extension can be uniquely lifted to the base of  local Artin
$\mathbb{C}$-algebra $A$. Consider $T_A\in GL(V(A))$ for a free $A$-module $V(A)$  of rank $r$ with $q(T_A)=T\in GL(V)$. The characteristic polynomial $P(T)$ of $T$ is given by $P_T(x)=\det(x\mathrm{Id}-T)=\prod\limits_{l=1}^{t}(x-\lambda_l)^{m_l}$ for distinct eigenvalues $\lambda_l$ of $T$ with multiplicity $m_l$, and the characteristic polynomial $P_{T_A}(x)$ of $T_A$ is given by $P_{T_A}(x)=\det(x\mathrm{Id}-T_A)\in A[x]$. It is clear that $q(P_{T_A}(x))=P_T(x)$. By Hensel lemma, there are unique  monic polynomials $P_1(x),\cdots,P_t(x)\in A[x]$ satisfying that
\begin{itemize}
  \item they are pairwise coprime, namely for any two polynomials $P_a(x), P_b(x)$ one finds polynomials $K_a(x), K_b(x)\in A[x]$ such that $K_a(x)P_a(x)+K_b(x)P_b(x)=1$;
  \item  $q(P_l(x))=(x-\lambda_l)^{m_l}$;
  \item $P_{T_A}(x)=\prod\limits_{l=1}^{t}P_l(x)$.
\end{itemize}
Consequently, $V(A)$ decomposes into $V(A)=\bigoplus\limits_{l=1}^tV_l(A)$ with $V_l(A)=\Ker(P_l(T_A))$. The restriction of $T_A$ on $V_l(A)$ is denoted by ${_l}T_A$, then one writes ${_l}T_A=\lambda_l(\mathrm{Id}+N_l)$, where $N_l\in \mathrm{End}_A(V_l(A))$ is nilpotent. Fixing a complex number $\tau_l$ satisfying $e^{2\pi \sqrt{-1}\tau_l}=\lambda_l$ with $\mathrm{Re}(\tau_l)\in[0,1)$, the operator
\begin{align*}
  R_l&=\tau_l\mathrm{Id}+\frac{1}{2\pi\sqrt{-1}}\log(\mathrm{Id}+N_l)\\
 & =\tau_l\mathrm{Id}+\frac{1}{2\pi\sqrt{-1}}\sum\limits_{k=1}^{L_l}\frac{(-1)^{k+1}N_l^k}{k}
\end{align*}
makes sense.
In particular, we have $e^{2\pi\sqrt{-1}R_l}={_l}T_A$.
 \end{proof}

We introduce the following groupoids:
\begin{itemize}
  \item the groupoid $\mathbb{PFB}(r,\mathbf{d})(A)$ is defined as the disjoint union,  over  the  set $ \mathrm{Iso}\mathfrak{PB}(r,\mathbf{d})(A)$ of isomorphism classes of  parabolic  $A$-bundles
 of rank $r$ and type $\mathbf{d}$, of the transformation groupoid $\mathbb{LC}_{(E_A,{\mathbf{F}_A})}$ corresponding  to $(\mathfrak{LC}_{(E_A,{\mathbf{F}_A})}, \Aut(E_A,\mathbf{F}_A))$;
 \item the  groupoid  $\mathbb{PFB}^{[0,1)}(r,\mathbf{d})(A)$ is defined to have the same objects as those in $\mathbb{PFB}(r,\mathbf{d})(A)$ and have morphisms given by those in $\mathbb{PFB}(r,\mathbf{d})(A)$ together with the RHD-relatedness\footnote{Similarly, two  parabolic logarithmic flat  $A$-bundles $( E_A,  {\mathbf{F}_A},\nabla_A)$ and $( E'_A,  {\mathbf{F}'_A},\nabla'_A)$ are RHD-related  if they are non-isomorphic and  $\mathrm{RHD}_A^{-1}( E_A,  {\mathbf{F}_A},\nabla_A)=\mathrm{RHD}_A^{-1}( E'_A,  {\mathbf{F}'_A},\nabla'_A)$.};
       \item the groupoid $\mathbb{PFB}_{( E, {\mathbf{F}}, \nabla)}(A)$ is defined as the disjoint union, over the set $ \mathrm{Iso}\mathfrak{PB}_{(E,\mathbf{F})}(A)$ of isomorphism classes of parabolic  $A$-bundles $(E_A,\mathbf{F}_A)$ of rank $r$ and type $\mathbf{d}$ with $q^\sharp(E_A)=E$ and $q({\mathbf{F}_A})=\mathbf{F}$, of the groupoid $\mathbb{LC}_{(E, \mathbf{F},\nabla, E_A,{\mathbf{F}_A})}$ corresponding  to $(\mathfrak{LC}_{(E, \mathbf{F},\nabla, E_A,{\mathbf{F}_A})}, \Aut_0(E_A,\mathbf{F}_A))$;
               \item the groupoid $\mathbb{PFB}'_{( E, {\mathbf{F}}, \nabla)}(A)$ is defined to have    the same objects as those in $\mathbb{PFB}_{( E, {\mathbf{F}}, \nabla)}(A)$ and have morphisms given by  those in the kernels of $\varepsilon_A$.
\end{itemize}

  We can also discuss parabolic logarithmic flat vector bundles in $C^\infty$-categories, as in Goldman--Millson's work \cite[Sections 5,6]{gm}. More precisely, denoting by $C^\infty(E,\mathbf{F})$ the underlying $C^\infty$-parabolic vector bundle of $(E,\mathbf{F})\in \mathfrak{PB}(r,\mathbf{d})$, let $\mathfrak{HCLC}_{C^\infty(E,\mathbf{F})}$ be the set/category consisting of holomorphic structures $ \bar\partial_E$ on $E$ and holomorphic logarithmic connections $\nabla$ compatible with the parabolic structure $\mathbf{F}$. There is a   group homomorphism $\varepsilon:\Aut(C^\infty(E,\mathbf{F})))\rightarrow G$ such that the following commutative diagram
\begin{align*}
   \CD
     \mathfrak{HCLC}_{C^\infty(E,\mathbf{F})}@>\mathrm{RHD}^{-1}>>R{(\Gamma,G;\mathbf{d})} \\
     @V \Phi VV @V \varepsilon(\Phi) VV  \\
     \mathfrak{HCLC}_{C^\infty(E,\mathbf{F})} @>\mathrm{RHD}^{-1}>> R{(\Gamma,G;\mathbf{d})}
   \endCD
  \end{align*}
for any $\Phi\in\Aut(C^\infty(E,\mathbf{F}))$. The transformation groupoid corresponding to $(\mathfrak{HCLC}_{C^\infty(E,\mathbf{F})}, \Aut(C^\infty(E,\mathbf{F})))$ is denoted by $\mathbb{HCLC}_{C^\infty(E,\mathbf{F})}$. The disjoint union, over the set $\mathrm{Iso}C^\infty\mathfrak{PB}(r,\mathbf{d})$ of isomorphism classes of $C^\infty$-parabolic vector bundles of rank $r$ and type $\mathbf{d}$ over $\bar X$, of the groupoid $\mathbb{HCLC}_{C^\infty(E,\mathbf{F})}$ forms a groupoid denoted by $C^\infty\mathbb{PFB}(r,\mathbf{d})$. The  groupoid  $C^\infty\mathbb{PFB}^{[0,1)}(r,\mathbf{d})$ is defined to have the same objects as those in $C^\infty\mathbb{PFB}(r,\mathbf{d})(A)$ and have morphisms given by those in $\mathbb{PFB}(r,\mathbf{d})$ together with the RHD-relatedness.

  Fixing a local Artin $\mathbb{C}$-algebra $A\in\mathbf{Art}_\mathbb{C}$,
   we similarly define the groupoids $C^\infty\mathbb{PFB}(r,\mathbf{d})(A)$ and $C^\infty\mathbb{PFB}^{[0,1)}(r,\mathbf{d})(A)$.
   Let $\mathfrak{HCLC}_{C^\infty(E,\mathbf{F})}(A)$ denote the set/category  consisting of holomorphic structures $\bar \partial_{E_A}$ on $C^\infty(E)\otimes_\mathbb{C} A$ and holomorphic logarithmic connections $\nabla_A$ compatible with the parabolic structure $\mathbf{F}(A)=\{\mathcal{F}^{(i)}(A)\}_{i=1,\cdots,n}$.
Chosen a point $(\bar\partial_E,\nabla)\in\mathfrak{HCLC}_{C^\infty(E,\mathbf{F})}$, let $\mathfrak{HCLC}_{(C^\infty(E,\mathbf{F}), \bar\partial_E,\nabla)}(A)$ denote the subset/subcategory of $\mathfrak{HCLC}_{C^\infty(E,\mathbf{F})}(A)$ satisfying $q^\sharp(\bar \partial_{E_A})=\bar\partial_E, q^\sharp(\nabla_A)=\nabla$, and $\mathbb{HCLC}_{(C^\infty(E,\mathbf{F}), \bar\partial_E,\nabla)}(A)$ be the groupoid corresponding to $(\mathfrak{HCLC}_{(C^\infty(E,\mathbf{F}), \bar\partial_E,\nabla)}(A),\exp(\End(C^\infty(E,\mathbf{F}))\otimes_{\mathbb{C}} \mathfrak{m}_A)$.  Denoting by $(\rho, \mathbf{P})$ the parabolic representation pair corresponding to $(C^\infty(E,\mathbf{F}),\bar\partial_E, \nabla)$, then there is a group homomorphism
$$
\varepsilon_A: \exp(\End(C^\infty(E,\mathbf{F}))\otimes_{\mathbb{C}} \mathfrak{m}_A)\longrightarrow G^0(A)
$$
such that the following diagram commutes:
\begin{align*}
   \CD
    \mathfrak{HCLC}_{(C^\infty(E,\mathbf{F}), \bar\partial_E,\nabla)}(A)@>\mathrm{RHD}^{-1}_A>> \mathrm{Def}^R_{\rho}(A) \\
     @V \Phi_A VV @V \varepsilon_A(\Phi_A) VV  \\
 \mathfrak{HCLC}_{(C^\infty(E,\mathbf{F}), \bar\partial_E,\nabla)}(A)@>\mathrm{RHD}^{-1}_A>>  \mathrm{Def}^R_{\rho}(A)
   \endCD
  \end{align*}
for any $\Phi_A\in \exp(\End(C^\infty(E,\mathbf{F}))\otimes_{\mathbb{C}} \mathfrak{m}_A)$. We define the groupoid $\mathbb{HCLC}'_{(C^\infty(E,\mathbf{F}), \bar\partial_E,\nabla)}(A)$  to have the same objects as those in ${\mathbb{HCLC}}_{(C^\infty(E,\mathbf{F}), \bar\partial_E,\nabla)}(A)$ and have  morphisms given by those in the kernels of $\varepsilon_A$.

The various parabolic-groupoid versions of the Riemann--Hilbert--Deligne correspondence are summarized as follows.

\begin{proposition}
With the above notation, the Deligne canonical extension provides equivalences of the following groupoids:
\begin{enumerate}
       \item  ${\mathbb{PRP}}(r,\mathbf{d})$ and $\mathbb{PFB}^{[0,1)}(r,\mathbf{d})$;
       \item  ${\mathbb{PRP}}(r,\mathbf{d})(A)$ and  $\mathbb{PFB}^{[0,1)}(r,\mathbf{d})(A)$;
       \item $\mathbb{R}(\Gamma,G;\mathbf{d})$ and $C^\infty\mathbb{PFB}^{[0,1)}(r,\mathbf{d})$;
       \item $\mathbb{R}(\Gamma,G;\mathbf{d})(A)$ and $C^\infty\mathbb{PFB}^{[0,1)}(r,\mathbf{d})(A)$;
       \item  ${\mathbb{DEF}}_{(\rho,\mathbf{P})}(A)$ and $\mathbb{PFB}_{( E, {\mathbf{F}}, \nabla)}(A)$;
       \item  $\mathbb{DEF}^{R}_{\rho}(A)$ and $\mathbb{HCLC}_{(C^\infty(E,\mathbf{F}), \bar\partial_E,\nabla)}(A)$ .
\end{enumerate}
\end{proposition}

\begin{corollary}\label{zaa}
Let $(E,\mathbf{F},\nabla)$ be a parabolic logarithmic flat bundle and let $(\rho, \mathbf{P})$ be the corresponding parabolic representation pair.
\begin{enumerate}
  \item The analytic germ of $\mathfrak{PRP}(r,\mathbf{d})$ at $(\rho,\mathbf{P})$ pro-represents the functor $\mathbb{F}_{(E,\mathbf{F},\nabla)}: \mathbf{Art}_\mathbb{C}\rightarrow \mathbf{Set}$ defined by
      \begin{align*}
        A\longmapsto \mathrm{Iso}\mathbb{PFB}'_{( E, {\mathbf{F}}, \nabla)}(A),
      \end{align*}
      where $\mathrm{Iso}\mathbb{PFB}'_{( E, {\mathbf{F}}, \nabla)}(A)$ denotes the set of isomorphism classes in $\mathbb{PFB}'_{( E, {\mathbf{F}}, \nabla)}(A)$.

  \item The analytic germ of $\mathbb{R}(\Gamma,G;\mathbf{d})$ at $\rho$ pro-represents the functor $\mathbb{G}_{(E,\mathbf{F},\nabla)}: \mathbf{Art}_\mathbb{C}\rightarrow \mathbf{Set}$ defined by
      \begin{align*}
        A\longmapsto \mathrm{Iso}\mathbb{HCLC}'_{(C^\infty(E,\mathbf{F}), \bar\partial_E,\nabla)}(A),
      \end{align*}
      where $\mathrm{Iso}\mathbb{HCLC}'_{(C^\infty(E,\mathbf{F}), \bar\partial_E,\nabla)}(A)$ denotes the set of isomorphism classes in $\mathbb{HCLC}'_{(C^\infty(E,\mathbf{F}), \bar\partial_E,\nabla)}(A)$.
\end{enumerate}
\end{corollary}

\subsection{Zariski tangent spaces and tangent quadratic cones}

  Let $\mathcal{A}^{0,q}_{\bar X}$ be the sheaf of $(0,q)$-forms on $\bar X$, and set
   \begin{align*}
     \mathcal{A}^{p,q}_{(\bar X,D)}(\End(E,\mathbf{F}))=\mathcal{A}^{0,q}_{\bar X}\otimes_{\mathcal{O}_{\bar X}}(\End(E,\mathbf{F})\otimes_{\mathcal{O}_{\bar X}}\Omega^p_{\bar X}(\mathcal{D})).
   \end{align*}
The double complex  $({\mathcal{A}^{\bullet,\bullet}_{(\bar X,D)}(\End(E,\mathbf{F}))},   \bar\partial_E, \nabla)$ provides  resolutions of the parabolic de Rham complex
\begin{align*}
  \mathcal{C}_{(E,\mathbf{F},\nabla)}: \End(E,\mathbf{F})\xrightarrow{\nabla}\End(E,\mathbf{F})\otimes_{\mathcal{O}_{\bar X}}\Omega^1_{\bar X}(\mathcal{D}),
\end{align*}
 where  the induced holomorphic structure and logarithmic connection on $\End(E)$ are  also denoted by $\bar\partial_E$ and $\nabla$, respectively. Let
 \begin{align*}
 \Gamma^{\bullet}_{(\bar X,D)}(\End(E,\mathbf{F}))=\bigoplus\limits_{\bullet=p+q}\Gamma(\bar X, \mathcal{A}^{p,q}_{(\bar X,D)}(\End(E,\mathbf{F})))
 \end{align*}
  be the spaces of
global sections, then we   define the following complexes
:
\begin{align*}
\mathcal{R}_{(E,\mathbf{F},\nabla)}:& \ \Gamma^{0}_{(\bar X,D)}(\End(E,\mathbf{F}))\xrightarrow{d_\nabla} \Gamma^{1}_{(\bar X,D)}(\End(E,\mathbf{F}))\xrightarrow{d_\nabla} \Gamma^{2}_{(\bar X,D)}(\End(E,\mathbf{F})),\\
\tilde{\mathcal{R}}_{(E,\mathbf{F},\nabla)}:& \  \tilde \Gamma^{0}_{(\bar X,D)}(\End(E,\mathbf{F}))\xrightarrow{d_\nabla} \tilde\Gamma^{1}_{(\bar X,D)}(\End(E,\mathbf{F}))\xrightarrow{d_\nabla} \tilde\Gamma^{2}_{(\bar X,D)}(\End(E,\mathbf{F}))
\end{align*}
where \begin{align*}
        \tilde \Gamma^{i}_{(\bar X,D)}(\End(E,\mathbf{F}))= \left\{
                                                              \begin{array}{ll}
                                                               \Gamma^{i}_{(\bar X,D)}(\End(E,\mathbf{F})), & \hbox{$i=0$;} \\
                                                               \Gamma(\bar X,\mathcal{A}^{0,1}_{(\bar X,D)}(\End(E)))\bigoplus\Gamma(\bar X,\mathcal{A}^{1,0}_{(\bar X,D)}(\End(E,\mathbf{F}))), & \hbox{$i=1$;}\\
                                                               \Gamma^{i}_{(\bar X,D)}(\End(E)), & \hbox{$i=2$,}
                                                              \end{array}
                                                            \right.
      \end{align*}
and $d_\nabla=\bar\partial_E+\nabla$.
When the parabolic structure $\mathbf{F}$
is trivial, we will omit $\mathbf{F}$ in the above notations.

\begin{theorem}\
\begin{enumerate}
    \item Fix a parabolic logarithmic flat bundle   $(E,\mathbf{F},\nabla)\in\mathfrak{PFB}(r,\mathbf{d})$.
          \begin{enumerate}
             \item The Zariski tangent space of $\mathfrak{PFB}(r,\mathbf{d})$ at    $(E,\mathbf{F},\nabla)$  is given by
 \begin{align*}
   &ZT_{(E,\mathbf{F},\nabla)}\mathfrak{PFB}(r,\mathbf{d})=\mathfrak{PFB}_{( E, {\mathbf{F}}, \nabla)}(A_1)\\
   =& \ \{(\Theta, Y_1,\cdots, Y_n)\in \Gamma^{1}_{(\bar X,D)}(\End(E))\bigoplus(\bigoplus_{i=1}^n\mathfrak{g}/\mathfrak{p}_{\mathcal{F}^{(i)}}):\\
   & \ \ \ \ \ \ \ \ \ \ \    d_\nabla\Theta=0, [\mathrm{Res}(\Theta,x_i)]+[\mathrm{Res}(\nabla,x_i), Y_i]= 0, i=1,\cdots,n\}\\
   \simeq&\ \{\Theta\in \Gamma^{1}_{(\bar X,D)}(\End(E)): d_\nabla\Theta=0\}, 
   \end{align*}where $[\bullet]_i$ denotes the projection of $\mathfrak{g}$ onto $\mathfrak{g}/\mathfrak{p}_{\mathcal{F}^{(i)}}$;
               \item The tangent quadratic cone of $\mathfrak{PFB}(r,\mathbf{d})$ at    $(E,\mathbf{F},\nabla)$ is given by
  \begin{align*}
  & \mathrm{QC}_{(E,\mathbf{F},\nabla)}\mathfrak{PFB}(r,\mathbf{d})\\
  =&\ \{ (\Theta, Y_1,\cdots, Y_n)\in ZT_{(E,\mathbf{F},\nabla)}\mathfrak{PFB}(r,\mathbf{d}):  d_\nabla\Xi+\frac{1}{2}  [\Theta,\Theta]=0 \ \textrm{ for } \Xi\in \Gamma^{1}_{(\bar X,D)}(\End(E)),\\
  & \ \ \ \ \ \ \ \ \ \ \ \  [\mathrm{Res}(\Xi,x_i)]_i+[\mathrm{Res}(\nabla,x_i), Z_i]+\frac{1}{2}[ Y_i,[\mathrm{Res}(\nabla,x_i), Y_i]]=0\\
  &\ \ \ \ \ \ \ \ \ \ \ \   \ \ \ \ \ \ \  \ \ \ \ \ \ \ \ \ \ \ \ \ \ \ \ \ \ \ \   \textrm{ for }  Z_i\in \mathfrak{g}/( \mathfrak{p}_{\mathcal{F}^{(i)}}+[ Y_i,\mathfrak{p}_{\mathcal{F}^{(i)}}]), i=1,\cdots,n\}.
  \end{align*}
            \end{enumerate}
    \item Fix a parabolic bundle  $(E,\mathbf{F})\in \mathfrak{PB}(r,\mathbf{d})$ and a point $(\bar\partial_E,\nabla)\in\mathfrak{HCLC}_{C^\infty(E,\mathbf{F})}$.
             \begin{enumerate}
                \item The Zariski tangent space of $\mathfrak{HCLC}_{C^\infty(E,\mathbf{F})}$ at    $(\bar\partial_E,\nabla)$ is given by
                      \begin{align*}
   ZT_{(\bar\partial_E,\nabla)}\mathfrak{HCLC}_{C^\infty(E,\mathbf{F})}&=\mathfrak{HCLC}_{(C^\infty(E,\mathbf{F}), \bar\partial_E,\nabla)}(A_1)\\
   & = \{\Theta\in \tilde\Gamma^{1}_{(\bar X,D)}(\End(E,\mathbf{F})): d_\nabla\Theta=0\}.\end{align*}
                \item The tangent quadratic cone of $\mathfrak{HCLC}_{C^\infty(E,\mathbf{F})}$ at   $(\bar\partial_E,\nabla)$ is given by
                        \begin{align*}
                         &\mathrm{QC}_{(\bar\partial_E,\nabla)}\mathfrak{HCLC}_{C^\infty(E,\mathbf{F})}\\
                         =&\ \{\Theta\in ZT_{(\bar\partial_E,\nabla)}\mathfrak{HCLC}_{C^\infty(E,\mathbf{F})}: d_\nabla\Xi+\frac{1}{2}  [\Theta,\Theta]=0 \ \textrm{ for } \Xi\in \tilde\Gamma^{1}_{(\bar X,D)}(\End(E,\mathbf{F}))\}.
                        \end{align*}
              \end{enumerate}
\end{enumerate}

 \end{theorem}

\begin{proof}
The above theorem can be regarded as the linear version of Theorem \ref{o}, and the  treatment of infinitesimal deformations proceeds in parallel.
\end{proof}

\section{DGLA-descriptions of deformation functors}

\subsection{Review of Deligne--Goldman--Millson deformation theory}

  A differential graded Lie algebra (DGLA) is a graded $\mathbb{C}$-vector space $L=\bigoplus\limits_{p\in \mathbb{Z}}L^p$ with a family of differentials
  $d: L^p\rightarrow L^{p+1}$ and a family of bilinear  homogeneous graded skew-symmetric  brackets $[\bullet,\bullet]: L^p\otimes L^{q}\rightarrow L^{p+q}$ such  that $d$ is a derivation for the Lie bracket and the bracket satisfies the graded Jacobi identity. The quadratic cone $\mathrm{QC}_L$ of $L$ consists of all $u\in L^1$ such that $[u,u]=0$. Two DGLAs $L_1, L_2$ are quasi-isomorphic if there exists a finite zigzag $L_1\rightarrow X_1\leftarrow X_2\rightarrow\cdots \rightarrow X_n \leftarrow L_2$ of quasi-isomorphisms between DGLAs.

   For a local Artin $\mathbb{C}$-algebra $A\in \mathbf{Art}_\mathbb{C}$, let $\mathfrak{m}_A$ be its maximal ideal. Then $L\otimes_{\mathbb{C}} \mathfrak{m}_A$
  is a nilpotent DGLA with differential and bracket given respectively by
  \begin{align*}
    d(u\otimes a)&=du\otimes a,\\
    [u\otimes a,v\otimes b]&=[u,v]\otimes ab
  \end{align*}
  for $u,v\in L, a,b\in\mathfrak{m}_A $. In particular, $L^0\otimes_{\mathbb{C}}\mathfrak{m}_A $ is a nilpotent Lie algebra, and we have the corresponding Lie group $\exp(L^0\otimes\mathfrak{m}_A)$ given by the Baker--Campbell--Hausdorff formula, namely the group multiplication is given by
  \begin{align*}
    \exp(\alpha)\exp(\beta)=\exp(\alpha+\beta+\frac{1}{2}[\alpha,\beta]+\frac{1}{12}[\alpha,[\alpha,\beta]]-\frac{1}{12}[\beta,[\alpha,\beta]]+\cdots)
  \end{align*}
  for $\alpha,\beta\in L^0\otimes_{\mathbb{C}}\mathfrak{m}_A$.
   We have the following three functors:
  \begin{itemize}
    \item the exponential functor $\exp_L$ from $\mathbf{Art}_\mathbb{C}$ to $\mathbf{Grp}$, the category of groups, defined by
    \begin{align*}
     \exp_L:  A\mapsto \exp(L^0\otimes_{\mathbb{C}}\mathfrak{m}_A);
    \end{align*}
    \item the Maurer--Cartan functor $\mathrm{MC}_L$ from $\mathbf{Art}_\mathbb{C}$ to $\mathbf{Set}$, defined by
        \begin{align*}
        \mathrm{MC}_L: A\mapsto  \textrm{MC}(L,A)=\{\alpha\in L^1\otimes_{\mathbb{C}}\mathfrak{m}_A :d\alpha+\frac{1}{2}[\alpha,\alpha]=0\};
        \end{align*}
    \item the deformation functor $\mathrm{Def}_L$, which is the quotient functor $ \mathrm{MC}_L/\exp_L$, where $\exp(L^0\otimes_{\mathbb{C}}\mathfrak{m}_A)$ acts on $L^1\otimes_{\mathbb{C}}\mathfrak{m}_A$
  as
  \begin{align*}
       \exp(\alpha_0)*\alpha_1=\exp(\mathrm{ad}\alpha_0)\alpha_1 +\frac{\mathrm{Id}-\exp(\mathrm{ad} \alpha_0)}{\mathrm{ad} \alpha_0}d\alpha_0.
     \end{align*}
  \end{itemize}

 It is known that the cohomology spaces $H^0(L)$, $H^1(L)$, and $H^2(L)$ characterize the space of infinitesimal automorphisms, the tangent space, and the obstruction space associated to $\mathrm{Def}_L$, respectively.
  Every morphism $f:L\rightarrow M$ of DGLAs induces a natural transformation of the associated deformation functors $f_*:\mathrm{Def}_L\rightarrow \mathrm{Def}_M$.
Goldman and Millson showed that if the induced morphism on 0-th cohomology is surjective, on 1-th cohomology is bijective, and on 2-th cohomology is injective, then $f_*$ is an isomorphism of functors \cite[Theorem 2.4]{gm}. In other words, the deformation functor is a 1-homotopy invariant.

Let $\mathfrak{h}$ be a Lie algebra over $\mathbb{C}$, which is viewed as a DGLA in degree zero.  An $\mathfrak{h}$-augmented DGLA is a pair $(L,\varepsilon)$, where $L$ is a DGLA and the augmentation $\varepsilon: L\rightarrow\mathfrak{h}$ is a homomorphism of DALAs. Given a local Artin $\mathbb{C}$-algebra $A$, the  $\mathfrak{h}$-augmented Deligne--Goldman--Millson groupoid $\mathbb{DGM}_{\mathfrak{h}}(L,\varepsilon,A)$ is defined by
\begin{itemize}
  \item the set of objects being
  \begin{align*}
   \mathrm{Obj}(\mathbb{DGM}_{\mathfrak{h}}(L,\varepsilon,A))= \{(\alpha,\varphi)\in (L^1\otimes_{\mathbb{C}}\mathfrak{m}_A)\times \exp(\mathfrak{h}\otimes_{\mathbb{C}}\mathfrak{m}_A): d\alpha+\frac{1}{2}[\alpha,\alpha]=0\};
  \end{align*}
  \item the set of morphisms between $(\alpha,\varphi)$ and $(\beta,\psi)$ being
  \begin{align*}
   \Hom((\alpha,\varphi),(\beta,\psi)) =\left\{\lambda\in L^0\otimes_{\mathbb{C}}\mathfrak{m}_A: \exp(\lambda)*\alpha=\beta,\exp(\varepsilon(\lambda))\varphi=\psi\right\}.
  \end{align*}
\end{itemize}
In particular, if $\mathfrak{h}$ is a trivial Lie algebra, $\mathbb{DGM}(L,\varepsilon,A)$ reduces to the usual Deligne--Goldman--Millson groupoid $\mathbb{DGM}(L,A)$. In terms of the notations of \cite{gm}, one writes
$$
\mathbb{DGM}_{\mathfrak{h}}(L,\varepsilon,A)=\mathbb{DGM}(L,A)\bowtie \exp( \mathfrak{h}\otimes\mathfrak{m}_A)_\varepsilon.
$$
The deformation functor
$\mathrm{Def}_{L,\varepsilon}$ is defined as
\begin{align*}
 \mathrm{Def}_{L,\varepsilon}: \mathbf{Art}_\mathbb{C}&\longrightarrow \mathbf{Set}\\
 A&\longmapsto\mathrm{Iso} {\mathbb{DGM}}_{\mathfrak{h}}(L,\varepsilon,A),
\end{align*}
where $\mathrm{Iso} {\mathbb{DGM}}_{\mathfrak{h}}(L,\varepsilon,A)$ denotes the set of isomorphism classes in $\mathbb{DGM}_{\mathfrak{h}}(L,\varepsilon,A)$.
An $\mathfrak{h}$-augmented DGLA  $(L,\varepsilon)$ is called formal if there is a finite zigzag
$$
(L,\varepsilon)\rightarrow(L_1,\varepsilon_1)\leftarrow (L_2,\varepsilon_2)\rightarrow\cdots\rightarrow(L_m,\varepsilon_m)\leftarrow (H^\bullet(L),\varepsilon_H)
$$
of quasi-isomorphisms between $\mathfrak{h}$-augmented DGLA's, where $H^\bullet(L)$ is treated as an $\mathfrak{h}$-augmented DGLA equipped with trivial differentials and the augmentation $\varepsilon_H$ is the restriction of $\varepsilon$ to $H^0(L)\subset L^0$.
Let $(L,\varepsilon)$ be a formal $\mathfrak{h}$-augmented DGLA, where the augmentations $\varepsilon$ is surjective and $\varepsilon|_{H}$ is injective. Then Goldman and Millson showed that the analytic germ of $\mathrm{QC}_{H^\bullet(L)}\times \mathfrak{h}/\varepsilon(H^0(L))$ pro-represents the functor $A\mapsto \mathrm{Iso} {\mathbb{DGM}}(L',A)$, where $L'=\Ker(\varepsilon)$ is the augmentation ideal \cite[Theorem 3.5]{gm}.

As one of the main applications of the Goldman--Millson theorem, together with the formality due to Deligne, Griffiths, Morgan, and Sullivan \cite{d}, they showed that if $\rho\in\mathrm{Hom}(\pi_1(Y), H)$, where $Y$ is a compact K\"{a}hler manifold, and $H$ is a real  algebraic Lie group, such that one of the following conditions is satisfied
\begin{itemize}
  \item the image of $\rho$ lies in a
compact subgroup of $H$;
  \item $\rho$ is semisimple;
  \item $H/V$ is a classifying space for polarized real Hodge structures, and the principal $H$-bundle associated to $\rho$ admits a horizontal holomorphic $V$-reduction;
      \item  $H/V$ is  a Hermitian symmetric space with automorphism group $H$ and the principal $H$-bundle associated to $\rho$ admits a holomorphic $V$-reduction,
\end{itemize}
then the representation variety $\mathrm{Hom}(\pi_1(Y), H)$ is quadratic at $\rho$. Hence the tangent cone to $\mathrm{Hom}(\pi_1(Y), H)$  at $\rho$ is a quadratic cone,
and $\mathrm{Hom}(\pi_1(Y), H)$ is locally analytically isomorphic to its tangent cone at $\rho$.

\subsection{DGLAs associated to parabolic logarithmic flat bundles}

Let $(E,\mathbf{F},\nabla)$ be a parabolic logarithmic flat bundle over $(\bar X,D)$. The following two  DGLAs are taken into account:
\begin{itemize}
  \item  the DGLA $M_{(E,\mathbf{F},\nabla)}$ is defined by \begin{itemize}
  \item $M_{(E,\mathbf{F},\nabla)}^i=\Gamma^{i}_{(\bar X,D)}(\End(E,\mathbf{F}))$, $i=0,1,2$,
  \item the bracket is naturally induced  from the Lie bracket on $\End(E)$,
  \item the differential is given by $d_\nabla$;
\end{itemize}
\item  the DGLA $\tilde M_{(E,\mathbf{F},\nabla)}$ is defined by $\tilde M_{(E,\mathbf{F},\nabla)}^i=\tilde\Gamma^{i}_{(\bar X,D)}(\End(E,\mathbf{F}))$ for $i=0,1,2$, with the brackets and differentials the same as those in $M_{(E,\mathbf{F},\nabla)}$.
\end{itemize}

\begin{example}
    We define
\begin{itemize}
 \item the graded  vector space $N^\bullet_{(E,\mathbf{F},\nabla)}$   by
 \begin{align*}
   N_{(E,\mathbf{F},\nabla)}^0&=\Gamma^{0}_{(\bar X,D)}(\End(E,\mathbf{F})),\\
   N_{(E,\mathbf{F},\nabla)}^1&=\Gamma^{1}_{(\bar X,D)}(\End(E,\mathbf{F}))\bigoplus(\bigoplus_{i=1}^n \mathfrak{n}_i),\\
   N_{(E,\mathbf{F},\nabla)}^2&=\Gamma^{2}_{(\bar X,D)}(\End(E,\mathbf{F}))\bigoplus(\bigoplus_{i=1}^n \mathfrak{n}_i),
 \end{align*}
 where $\mathfrak{n}_i=\{a\in \mathfrak{p}_{\mathcal{F}^{(i)}}: a(\mathcal{F}^{(i)}_\ell)\subset \mathcal{F}^{(i)}_{\ell+1}, \ell=1,\cdots,s_i\}$;
  \item the differentials by
  \begin{align*}
    d^0:N_{(E,\mathbf{F},\nabla)}^0&\longrightarrow N_{(E,\mathbf{F},\nabla)}^1\\
    \alpha&\longmapsto (d_\nabla\alpha,0,\cdots,0),\\
    d^1:N_{(E,\mathbf{F},\nabla)}^1&\longrightarrow N_{(E,\mathbf{F},\nabla)}^2 \\
    (\Theta,Y_1,\cdots, Y_n)&\longmapsto (d_\nabla\Theta,[\mathrm{Res}(\nabla,x_1), Y_1],\cdots,[\mathrm{Res}(\nabla,x_n), Y_n]);
  \end{align*}
  \item the brackets  by
  \begin{align*}
   [\bullet,\bullet]: N_{(E,\mathbf{F},\nabla)}^0\times N_{(E,\mathbf{F},\nabla)}^0&\longrightarrow N_{(E,\mathbf{F},\nabla)}^0 \\
    (\alpha,\beta)&\longmapsto[\alpha,\beta],\\
    [\bullet,\bullet]: N_{(E,\mathbf{F},\nabla)}^0\times N_{(E,\mathbf{F},\nabla)}^1&\longrightarrow N_{(E,\mathbf{F},\nabla)}^1 \\
    (\alpha,(\Theta,Y_1,\cdots,Y_n))&\longmapsto([\alpha,\Theta],[\alpha|_{x_1}, Y_1],\cdots,[\alpha|_{x_n}, Y_n]),\\
    [\bullet,\bullet]: N_{(E,\mathbf{F},\nabla)}^0\times N_{(E,\mathbf{F},\nabla)}^2&\longrightarrow N_{(E,\mathbf{F},\nabla)}^2\\
    (\alpha,(\Lambda,Z_1,\cdots,Z_n))&\longmapsto([\alpha,\Lambda],[\alpha|_{x_1}, Z_1],\cdots,[\alpha|_{x_n}, Z_n]),\\
    [\bullet,\bullet]: N_{(E,\mathbf{F},\nabla)}^1\times N_{(E,\mathbf{F},\nabla)}^1&\longrightarrow N_{(E,\mathbf{F},\nabla)}^2 \\
        ( (\Theta,Y_1,\cdots,Y_n),(\Theta',Y'_1,\cdots,Y'_n)) &\longmapsto([\Theta,\Theta'],
        [\mathrm{Res}(\Theta,x_1), {Y'_1}]+[\mathrm{Res}(\Theta',x_1),   Y_1],\\
        &\ \ \ \ \ \ \ \ \ \ \ \cdots,[\mathrm{Res}(\Theta,x_n), {Y'_n}]+[\mathrm{Res}(\Theta',x_n),   Y_n]).
  \end{align*}
\end{itemize}
Then the above data form a DGLA, which is denoted by  $N_{(E,\mathbf{F},\nabla)}$.
\end{example}

\begin{definition}
A parabolic logarithmic flat bundle $(E,\mathbf{F},\nabla)$ over $(\bar X,D)$ is called \emph{very generic} if the residue $\mathrm{Res}(\nabla,x_i)$ of $\nabla$ at each $x_i$ is a regular semisimple element in $\mathfrak{g}$.
\end{definition}

\begin{proposition}\label{zaq}
If $(E,\mathbf{F},\nabla)$ is a very generic parabolic logarithmic flat bundle, then the deformation functors $\mathrm{Def}_{M_{(E,\mathbf{F},\nabla)}}$, $\mathrm{Def}_{\tilde M_{(E,\mathbf{F},\nabla)}}$, $\mathrm{Def}_{\tilde M_{(E,\nabla)}}$ and $\mathrm{Def}_{N_{(E,\mathbf{F},\nabla)}}$ are isomorphic.

\end{proposition}

\begin{proof}(1) We have the following short exact sequence of complexes
\begin{align*}
0\rightarrow \mathcal{R}_{(E,\mathbf{F},\nabla)}\rightarrow \tilde{\mathcal{R}}_{(E,\mathbf{F},\nabla)}\rightarrow \mathcal{S}\rightarrow 0,
\end{align*}
where the complex $\mathcal{S}$ is given by
\begin{align*}
  0\rightarrow \bigoplus_{i=1}^n S_i\xrightarrow{\bigoplus_{i=1}^n\mathrm{Res}(\nabla,x_i)} \bigoplus_{i=1}^nS_i.
\end{align*}
If the residue $\mathrm{Res}(\nabla,x_i)$  is a regular semisimple element, namely it is a  diagonalizable matrix  with distinct eigenvalues,  we have that the adjoint action of $\mathrm{Res}(\nabla,x_i)$ on $\mathfrak{g}/\mathfrak{p}_{\mathcal{F}^{(i)}}$ gives rise to an isomorphism. Therefore the condition of very genericity implies that
  $H^1(\mathcal{S})=H^2(\mathcal{S})=0$.  It follows from the induced  long exact sequence that the DGLA $ M_{(E,\mathbf{F},\nabla)}$ is quasi-isomorphic to the DGLA $\tilde{M}_{(E,\mathbf{F},\nabla)}$.
 By  Goldman--Millson theorem, the deformation functor $\mathrm{Def}_{M_{(E,\mathbf{F},\nabla)}}$ is isomorphic to the deformation functor $\mathrm{Def}_{\tilde M_{(E,\mathbf{F},\nabla)}}$.

(2) We have following short exact sequence of complexes
\begin{align*}
 0\rightarrow\mathcal{C}_{(E,\mathbf{F},\nabla)}\rightarrow\mathcal{C}_{(E,\nabla)}\rightarrow\mathcal{S}'\rightarrow 0,
\end{align*}
where the complex $\mathcal{S}'=\mathcal{S}[1]$ is given by
\begin{align*}
   \bigoplus_{i=1}^n S_i\xrightarrow{\bigoplus_{i=1}^n\mathrm{Res}(\nabla,x_i)} \bigoplus_{i=1}^nS_i.
\end{align*}
Again by  the condition of very genericity and the induced  long exact sequence, we have that the DGLA $ M_{(E,\mathbf{F},\nabla)}$ is quasi-isomorphic to the DGLA $\tilde{M}_{(E,\mathbf{F},\nabla)}$, thus the deformation functor $\mathrm{Def}_{M_{(E,\mathbf{F},\nabla)}}$ is isomorphic to the deformation functor $\mathrm{Def}_{\tilde M_{(E,\mathbf{F},\nabla)}}$.

(3)
 There is a natural embedding  $\iota: M_{(E,\mathbf{F},\nabla)}\rightarrow N_{(E,\mathbf{F},\nabla)}$. It is clear that $H^0(M_{(E,\mathbf{F},\nabla)})=H^0(N_{(E,\mathbf{F},\nabla)})$, and $\iota_*:H^2(M_{(E,\mathbf{F},\nabla)})\rightarrow H^2(N_{(E,\mathbf{F},\nabla)}) $ is injective. Let $C(\mathrm{Res}(\nabla,x_i),\mathfrak{g})$ be the centralizer of $\mathrm{Res}(\nabla,x_i)$ in $\mathfrak{g}$. Since $\mathrm{Res}(\nabla,x_i)$ is a regular semisimple element, we have $C(\mathrm{Res}(\nabla,x_i),\mathfrak{g})\bigcap \mathfrak{n}_i=0$. Hence $\iota_*:H^1(M_{(E,\mathbf{F},\nabla)})\rightarrow H^1(L_{(E,\mathbf{F},\nabla)}) $ is bijective.  which implies that the deformation functor $\mathrm{Def}_{M_{(E,\mathbf{F},\nabla)}}$ is isomorphic to the deformation functor $\mathrm{Def}_{N_{(E,\mathbf{F},\nabla)}}$.
\end{proof}

For a parabolic logarithmic flat bundle $(E,\mathbf{F},\nabla)$ over $(\bar X,D)$, we regard the parabolic de Rham complex
 \begin{align*}
   \mathcal{C}_{(E,\mathbf{F},\nabla)}: \End(E,\mathbf{F})\xrightarrow{\nabla}\End(E,\mathbf{F})\otimes_{\mathcal{O}_{\bar X}}\Omega^1_{\bar X}(\mathcal{D})
 \end{align*}
as a sheaf of DGLAs. Let $\mathcal{A}_{(E,\mathbf{F},\nabla)}$ be a sheaf of DGLAs providing  a resolution $\mathcal{C}_{(E,\mathbf{F},\nabla)}\rightarrow{\mathcal{A}_{(E,\mathbf{F},\nabla)}}$ of $\mathcal{C}_{(E,\mathbf{F},\nabla)}$ as a  morphism of sheaves of DGLAs which induces isomorphisms between hypercohomologies. Moreover, assume that this resolution is acyclic, then we have a DGLA $\Gamma(\bar X,{\mathcal{A}_{(E,\mathbf{F},\nabla)}})$.

\begin{theorem}\label{za}
Let $(E,\mathbf{F},\nabla)$ be a parabolic logarithmic flat bundle over $(\bar X,D)$ and $(\rho, \mathbf{P})$ be the corresponding parabolic representation pair. The deformation functor $ \mathrm{Def}_{\Gamma(\bar X, {\mathcal{A}_{(E,\mathbf{F},\nabla)}})}$ is isomorphic to the functor $\overline{\mathrm{ Def}}_{(\rho,\mathbf{P})}$, which is defined by
 \begin{align*}
 \overline{\mathrm{ Def}}_{(\rho,\mathbf{P})}: \mathbf{Art}_\mathbb{C}&\longrightarrow \mathbf{Set}\\
 A&\longmapsto\mathrm{Iso} \mathbb{DEF}_{(\rho,\mathbf{P})}(A)
 \end{align*}
for $\mathrm{Iso} \mathbb{DEF}_{(\rho,\mathbf{P})}(A)$ denoting the set of isomorphism classes in $\mathbb{DEF}_{(\rho,\mathbf{P})}(A)$.
 \end{theorem}

\begin{proof}
\textbf{\emph{Step 1:}} Given a parabolic logarithmic flat bundle $(E,\mathbf{F},\nabla)$ and a local Artin $\mathbb{C}$-algebra $A\in \mathbf{Art}_\mathbb{C}$, let $\mathcal{U} = \{U_i\}_{i\in I}$ be an affine open cover of $\bar X$. Picking
$$
\alpha=\{\alpha_i\}_{i\in I}\in\prod\limits_{i\in I}\Gamma(U_i,(\End(E,\mathbf{F})\otimes_{\mathcal{O}_{\bar X}}\Omega^1_{\bar X}(\mathcal{D}))\otimes_\mathbb{C}\mathfrak{m}_A),
$$
we obtain a deformed parabolic de Rham complex
 \begin{align*}
   E|_{U_i}\otimes_{\mathbb{C}} A\xlongrightarrow{\iota^\sharp(\nabla)+\alpha_i} (E\otimes_{\mathcal{O}_{\bar X}}\Omega^1_{\bar X}(D))|_{U_i}\otimes_{\mathbb{C}} A,
 \end{align*}where $\iota^\sharp(\nabla)=\iota^*\circ q^*(\nabla)$.
 We need to specify isomorphisms between the deformed complexes
on the double intersections $U_{ij}$ of the cover $\mathcal{U}$. Namely, for $m=\{m_{ij}\}_{i<j}\in \prod\limits_{i<j}\Gamma(U_{ij}, \End(E,\mathbf{F})\otimes_{\mathbb{C}}\mathfrak{m}_A)$, we have the following commutative diagram
 \begin{align*}
   \CD
  E|_{U_{ij}}\otimes_{\mathbb{C}} A @>(\iota^\sharp(\nabla)+\alpha_j)|_{U_{ij}}>> (E\otimes_{\mathcal{O}_{\bar X}}\Omega^1_{\bar X}(D))|_{U_{ij}}\otimes_{\mathbb{C}} A \\
  @V e^{m_{ij}} VV @V e^{m_{ij}} VV  \\
  E|_{U_{ij}}\otimes_{\mathbb{C}} A @>(\iota^\sharp(\nabla)+\alpha_i)|_{U_{ij}}>> (E\otimes_{\mathcal{O}_{\bar X}}\Omega^1_{\bar X}(D))|_{U_{ij}}\otimes_{\mathbb{C}} A
\endCD,
 \end{align*}
 which implies that
\begin{align}\label{a}
  \alpha_i|_{U_{ij}}=e^{m_{ij}}* \alpha_j|_{U_{ij}}\quad \text{for all}\ i<j.
\end{align}
The above isomorphisms must satisfy the cocycle condition
   \begin{align}\label{b}
     m_{jk}|_{U_{ijk}}\star(-m_{ik}|_{U_{ijk}})\star m_{ij}|_{U_{ijk}}=0\quad \text{for all}\ i<j<k,
   \end{align}
where the product $\star$ is defined by
 \begin{align*}
  X\star Y=\log(\exp X\exp Y)
 \end{align*}
 for $X,Y$ lying in a certain Lie algebra. The pair $(\alpha,m)$ satisfying the above conditions \eqref{a} and \eqref{b} gives rise to a parabolic logarithmic flat $A$-bundle $(E_A, F(A),\nabla_A)\in\mathfrak{PFB}_{( E, {\mathbf{F}}, \nabla)}(A)$. Conversely, if $(\alpha,m)$ and $(\alpha',m')$ produce the isomorphic parabolic flat $A$-bundles, then there exists
 $$
 a=\{a_i\}_{i\in I}\in\prod\limits_{i\in I}\Gamma(U_i,\End(E,\mathbf{F})\otimes_{\mathbb{C}}\mathfrak{m}_A)
 $$
 such that the diagrams
\begin{align*}
\CD
  E|_{U_i}\otimes_{\mathbb{C}} A @>\iota^\sharp(\nabla)+\alpha_i>>(E\otimes_{\mathcal{O}_{\bar X}}\Omega^1_{\bar X}(D))|_{U_i}\otimes_{\mathbb{C}} A  \\
  @V e^{a_i} VV @V e^{a_i} VV  \\
  E|_{U_i}\otimes_{\mathbb{C}} A @>\iota^\sharp(\nabla)+\alpha'_i>>(E\otimes_{\mathcal{O}_{\bar X}}\Omega^1_{\bar X}(D))|_{U_i}\otimes_{\mathbb{C}} A
\endCD
\end{align*}
and
 \begin{align*}
 \CD
 ( E|_{U_{ij}}\otimes_{\mathbb{C}} A, (\iota^*(\nabla)+\alpha_i)|_{U_{ij}}) @> e^{m_{ij}}>> ( E|_{U_{ij}}\otimes_{\mathbb{C}} A, (\iota^*(\nabla)+\alpha_i)|_{U_{ij}}) \\
  @V (e^{a_j})|_{U_{ij}}VV @V e^{a_i}|_{U_{ij}} VV  \\
  ( E|_{U_{ij}}\otimes_{\mathbb{C}} A, (\iota^*(\nabla)+\alpha_i)|_{U_{ij}}) @> e^{m'_{ij}}>> ( E|_{U_{ij}}\otimes_{\mathbb{C}} A, (\iota^*(\nabla)+\alpha_i)|_{U_{ij}})
\endCD
\end{align*}
commute, hence
\begin{align*}
 \alpha_i'&=e^{a_i}*\alpha_i,\\
 m'_{ij}&=a_i|_{U_{ij}}\star m_{ij}\star (-a_j|_{U_{ij}}).
\end{align*}

Define a DGLA $\mathfrak{g}_\kappa(\mathcal{U})$ by
\begin{align*}
 \mathfrak{g}^0_\kappa(\mathcal{U})&=\prod\limits_{i_1<\cdots< i_\kappa}\Gamma(U_{i_1\cdots i_\kappa}, \End(E,\mathbf{F})),\\
 \mathfrak{g}^1_\kappa(\mathcal{U})&=\prod\limits_{i_1<\cdots< i_\kappa}\Gamma(U_{i_1\cdots i_\kappa}, \End(E,\mathbf{F})\otimes_{\mathcal{O}_{\bar X}}\Omega^1_{\bar X}(D))
\end{align*}
with differential $\nabla: \mathfrak{g}^0_\kappa(\mathcal{U})\rightarrow \mathfrak{g}^1_\kappa(\mathcal{U})$ and Lie bracket induced from $\End(E,\mathbf{F})$. This yields a \v{C}ech cosimplicial DGLA  \cite{f,f1}
\begin{align*}
\mathfrak{g}_\mathcal{U}:  \xymatrix{
 \mathfrak{g}_1(\mathcal{U})\ar@<2pt>[r]\ar@<-2pt>[r]& \mathfrak{g}_2(\mathcal{U})\ar@<4pt>[r] \ar[r] \ar@<-4pt>[r]& \mathfrak{g}_3(\mathcal{U})\ar@<6pt>[r] \ar@<2pt>[r] \ar@<-2pt>[r] \ar@<-6pt>[r]& \cdots},
\end{align*}
where, for each $\kappa$ there are $\kappa+1$ morphisms of DGLAs
  \begin{align*}
   \delta_{\lambda,\kappa}: \mathfrak{g}_\kappa(\mathcal{U})\rightarrow \mathfrak{g}_{\kappa+1}(\mathcal{U}),\ \lambda=1,\cdots,\kappa+1,
  \end{align*}
 defined by sending $\theta\in \Gamma(U_{i_1\cdots \hat{i_\lambda}\cdots i_\kappa}, \End(E,\mathbf{F})\otimes_{\mathcal{O}_{\bar X}}\Omega^\bullet_{\bar X}(D))$ to $\delta_{\lambda,\kappa}(\theta)=\theta|_{U_{i_1\cdots i_\lambda\cdots i_\kappa}}$, with relations
 \begin{align*}
    \delta_{\lambda'+1,\kappa+1}\circ  \delta_{\lambda,\kappa}= \delta_{\lambda,\kappa+1} \circ\delta_{\lambda',\kappa}
 \end{align*}
  for all $\lambda'\geq \lambda$.

For a local Artin $\mathbb{C}$-algebra $A\in \mathbf{Art}_\mathbb{C}$, one defines the groupoid  $\mathfrak{H}^1(\exp(\mathfrak{g}_\mathcal{U}),A)$ as follows:
\begin{itemize}
  \item the set of objects is
\begin{align*}
\mathrm{Obj}(\mathfrak{H}^1(\exp(\mathfrak{g}_\mathcal{U}),A))=\left\{(\alpha,m)\in(\mathfrak{g}^1_1(\mathcal{U})\oplus \mathfrak{g}^0_2(\mathcal{U}) )\otimes_{\mathbb{C}} \mathfrak{m}_A:\
\begin{aligned}
&\ \delta_{2,1}\alpha=e^m*\delta_{2,1}\alpha,\\
&\ \delta_{1,2}m\star(-\delta_{2,2}m)*\delta_{3,2}m=0
\end{aligned}
\right\};
\end{align*}
  \item the set of morphisms between $(\alpha,m)$ and $(\beta,n)$ is
\begin{align*}
\Hom((\alpha,m),(\beta,n))=\left\{\mathfrak{g}_1^0(\mathcal{U})\otimes_{\mathbb{C}}\mathfrak{m}_A:\
\begin{aligned}
&\ \exp(a)*\alpha=\beta,\\
&\ \delta_{2,1}a\star m\star (-\delta_{1,1}a) =n
\end{aligned}
\right\}.
\end{align*}
\end{itemize}
Take another  affine open cover $\mathcal{U}' = \{U'_{i'}\}_{i'\in I'}$ of $\bar X$, which is a refinement of $\mathcal{U} $. Let $\rho_1, \rho_2: I'\rightarrow I$ be two refinement maps, then they induce functors $\tilde \rho_1, \tilde \rho_2$ from $\mathfrak{H}^1(\exp(\mathfrak{g}_\mathcal{U}),A)$ to $\mathfrak{H}^1(\exp(\mathfrak{g}_{\mathcal{U}'}),A)$, defined by
\begin{align*}
  \tilde \rho_{1,2}(\{\alpha_i\},\{m_{ij}\})=(\{\alpha^{1,2}_{i'}\},\{m^{1,2}_{i'j'}\}),
\end{align*}
where
 \begin{align*}
  \alpha^{1,2}_{i'}&=\alpha_{\rho_{1,2}(i')}|_{U'_{i'}},\\
   m^{1,2}_{i'j'}&=m_{\rho_{1,2}(i')\rho_{1,2}(j')}|_{U'_{i'j'}}.
 \end{align*}
 One can show that  $\tilde \rho_1, \tilde \rho_2$ are isomorphic.
Indeed, by \eqref{a} and \eqref{b}, we have
\begin{align*}
 e^{m_{\rho_{2}(i')\rho_{1}(i')}|_{U'_{i'}}}*\alpha^{1}_{i'}&=\alpha^{2}_{i'},\\
  m_{\rho_2(i')\rho_1(i')}|_{U'_{i'}}\star m^1_{i'j'}\star (-m_{\rho_2(j')\rho_1(j')}|_{U'_{j'}})&= m^2_{i'j'},
\end{align*}
thus the natural morphism  $\xi_{(\alpha,m)}: \tilde \rho_{1}(\alpha,m)\mapsto \tilde \rho_{2} (\alpha,m)$ is determined  by $\{a_{i'}= m_{\rho_2(i')\rho_1(i')}|_{U'_{i'}}\}_{i'\in I'}\in \mathfrak{g}_1^0(\mathcal{U}')\otimes_{\mathbb{C}}\mathfrak{m}_A$.
Consequently, the groupoid
 \begin{align*}
   \mathfrak{H}^1(\exp(\mathcal{C}_{(E,\mathbf{F},\nabla)}),A)=\lim\limits_{\substack{\longrightarrow \\ U}}\mathfrak{H}^1(\exp(\mathfrak{g}_\mathcal{U}),A)
 \end{align*}
  is well-defined. Then the groupoid $\mathfrak{H}^1(\exp(\mathcal{C}_{(E,\mathbf{F},\nabla)}),A)$ is equivalent to the groupoid $\mathbb{PFB}_{( E, {\mathbf{F}}, \nabla)}(A)$.

\textbf{\emph{Step 2:}} Following \cite{f,f1}, we introduce the Thom--Whitney--Sullivan  DGLA (TWS DGLA).
Define the standard $\kappa$-simplex as the affine space
\begin{align*}
  \Delta^\kappa=\{(t_1,\cdots,t_{\kappa+1})\in \mathbb{C}^{\kappa+1}: t_1+\cdots+t_{\kappa+1}=1\},
\end{align*}
then let $\Lambda_\kappa^p$ ($0\leq p\leq \kappa$) be the space of polynomial differential $p$-forms on $\Delta^\kappa$, and set $\Lambda_\kappa=\bigoplus\limits_{p=0}^\kappa\Lambda_\kappa^p$ as a differential graded algebra.
The face map $\partial_{\lambda,\kappa}^*: \Lambda_{\kappa+1}\rightarrow  \Lambda_{\kappa}$ ($1\leq \lambda\leq \kappa+2$) is given by pullback of the forms under the inclusion
$$
\partial_{\lambda,\kappa}:  \Delta^\kappa\rightarrow  \Delta^{\kappa+1},\quad (t_1,\cdots, t_\lambda,\cdots, t_{\kappa+1})\mapsto( t_1,\cdots,t_{\lambda-1},0, t_\lambda,\cdots, t_{\kappa+1}).
$$
The Thom--Whitney--Sullivan bicomplex $\mathrm{TWS}^{\bullet,\bullet}_\mathcal{U}(\mathcal{C}_{(E,\mathbf{F},\nabla)})$ associated to the \v{C}ech cosimplicial DGLA $\mathfrak{g}_{\mathcal{U}}$ is defined by
\begin{align*}
( \mathrm{TWS}^{p,q}_\mathcal{U}(\mathcal{C}_{(E,\mathbf{F},\nabla)})=\{(x_\kappa)\in\prod\limits_{\kappa\geq1}\Lambda^p_\kappa\otimes \mathfrak{g}^q_\kappa(\mathcal{U}): (\partial^*_{\lambda,\kappa-1}\otimes \mathrm{Id})x_\kappa=( \mathrm{Id}\otimes \delta_{\lambda,\kappa-1})x_{\kappa-1}, 1\leq \lambda\leq \kappa+1\}, d, \nabla).
\end{align*}
The TWS DGLA $\mathrm{TWS}_\mathcal{U}(\mathcal{C}_{(E,\mathbf{F},\nabla)})$ is then defined by the corresponding total complex
$$
\mathrm{TWS}_\mathcal{U}^i(\mathcal{C}_{(E,\mathbf{F},\nabla)})=\bigoplus\limits_{p+q=i}\mathrm{TWS}^{p,q}_\mathcal{U}(\mathcal{C}_{(E,\mathbf{F},\nabla)}),
$$
with differential
$$
d^i_{\mathrm{TWS}}=\sum\limits_{p+q=i}d+(-1)^p\nabla,
$$
and the Lie bracket induced from $\mathfrak{g}_\kappa(\mathcal{U})$.

Let $\mathrm{Def}_{\mathrm{TWS}_\mathcal{U}(\mathcal{C}_{(E,\mathbf{F},\nabla)})}$ be the deformation functor associated to the TWS DGLA $\mathrm{TWS}_\mathcal{U}(\mathcal{C}_{(E,\mathbf{F},\nabla)})$.  By \cite[Theorem 4.11]{f},  every refinement $\mathcal{U}'$  of $\mathcal{U}$ induces a natural morphism of deformation functors $\mathrm{Def}_{\mathrm{TWS}_\mathcal{U}(\mathcal{C}_{(E,\mathbf{F},\nabla)})}\rightarrow\mathrm{Def}_{\mathrm{TWS}_{\mathcal{U}'}(\mathcal{C}_{(E,\mathbf{F},\nabla)})}$. Hence,
the direct limit
 \begin{align*}
  \mathrm{Def}_{\mathrm{TWS}(\mathcal{C}_{(E,\mathbf{F},\nabla)})}=\lim\limits_{\substack{\longrightarrow \\ U}}\mathrm{Def}_{\mathrm{TWS}_\mathcal{U}(\mathcal{C}_{(E,\mathbf{F},\nabla)})}
 \end{align*}
is well-defined. Moreover, there is a natural isomorphism of functors
\begin{align*}
  \mathrm{Def}_{\mathrm{TWS}(\mathcal{C}_{(E,\mathbf{F},\nabla)})}\simeq \mathfrak{H}^1_{\exp(\mathcal{C}_{(E,\mathbf{F},\nabla)})},
\end{align*}
where the functor $\mathfrak{H}^1_{\exp(\mathcal{C}_{(E,\mathbf{F},\nabla)})}$ is defined by
\begin{align*}
\mathfrak{H}^1_{\exp(\mathcal{C}_{(E,\mathbf{F},\nabla)})}: \mathbf{Art}_\mathbb{C}&\longrightarrow \mathbf{Set}\\
 A&\longmapsto\mathrm{Iso}\mathfrak{H}^1(\exp(\mathcal{C}_{(E,\mathbf{F},\nabla)}),A),
\end{align*}
where $\mathrm{Iso}\mathfrak{H}^1(\exp(\mathcal{C}_{(E,\mathbf{F},\nabla)}),A)$ denotes the set of isomorphism classes in $\mathfrak{H}^1(\exp(\mathcal{C}_{(E,\mathbf{F},\nabla)}),A)$.

\textbf{\emph{Step 3:}} Assume $\mathcal{U}$ is an  open cover of $\bar X$, which is acyclic with respect to both $\mathcal{C}_{(E,\mathbf{F},\nabla)}$ and $\mathcal{A}_{(E,\mathbf{F},\nabla)}$. For ${\mathcal{A}_{(E,\mathbf{F},\nabla)}}$, one can also define the corresponding  \v{C}ech cosimplicial DGLA and TWS DGLA $\mathrm{TWS}_\mathcal{U}({\mathcal{A}_{(E,\mathbf{F},\nabla)}})$. Then there is a natural morphism of DGLAs (see \cite[Theorem 3.7]{f1})
$$
\Gamma(\bar X, {\mathcal{A}_{(E,\mathbf{F},\nabla)}})\rightarrow \mathrm{TWS}_\mathcal{U}({\mathcal{A}_{(E,\mathbf{F},\nabla)}}),
$$
which is a quasi-isomorphism since their $i$-th cohomologies are both isomorphic to the hypercohomology $\mathbb{H}^i(\mathcal{C}_{(E,\mathbf{F},\nabla)})$. On the other hand, the natural morphism of DGLA's
$$
\mathrm{TWS}_\mathcal{U}(\mathcal{C}_{(E,\mathbf{F},\nabla)})\rightarrow \mathrm{TWS}_\mathcal{U}({\mathcal{A}_{(E,\mathbf{F},\nabla)}})
$$
is also a quasi-isomorphism. Therefore, the DGLA $\Gamma(\bar X, {\mathcal{A}_{(E,\mathbf{F},\nabla)}})$ is quasi-isomorphic to the DGLA $\mathrm{TWS}_\mathcal{U}(\mathcal{C}_{(E,\mathbf{F},\nabla)})$. By the Goldman--Millson theorem, we have following isomorphism of deformation functors:
$$
\mathrm{Def}_{\Gamma(\bar X, {\mathcal{A}_{(E,\mathbf{F},\nabla)}})}\simeq\mathrm{Def}_{\mathrm{TWS}_\mathcal{U}(\mathcal{C}_{(E,\mathbf{F},\nabla)})}.
$$
Then combining Step 1 and Step 2, we obtain the final conclusion.
\end{proof}

\begin{corollary}The deformation functor $ \mathrm{Def}_{M_{(E,\mathbf{F},\nabla)}}$ is isomorphic to the functor $\overline{\mathrm{ Def}}_{(\rho,\mathbf{P})}$.
\end{corollary}
\begin{remark}
It is known that, for the local Artin $\mathbb{C}$-algebra $A_1$, $\mathrm{Iso}\mathfrak{H}^1(\exp(\mathcal{C}_{(E,\mathbf{F},\nabla)}),A_1)$ is isomorphic to the hypercohomology $\mathbb{H}^1(\mathcal{C}_{(E,\mathbf{F},\nabla)})\simeq \mathbb{H}^1(\mathcal{C}_{(E,\nabla)})$ as sets \cite{y,m,i,hh,bi,b,si2,bm}.
\end{remark}

\begin{corollary}Let $(E,\mathbf{F},\nabla)\in \mathfrak{PFB}(r,\mathbf{d})$ be  a very generic parabolic logarithmic flat bundle and $(\rho,\mathbf{P})\in \mathfrak{PRP}(r,\mathbf{d})$ be the corresponding
parabolic representation pair. Then the parabolic representation pair variety $\mathfrak{PRP}(r,\mathbf{d})$ is non-singular at $(\rho,\mathbf{P})$.
\end{corollary}

\begin{proof}By Proposition \ref{zaq}, the deformation functor $\mathrm{Def}_{ M_{(E,\mathbf{F},\nabla)}}$ is isomorphic to the deformation functor $\mathrm{Def}_{ M_{(E,\nabla)}}$, then from \cite[Lemma 3.8]{gm}, it follows that the $\mathfrak{g}$-augmented Deligne-Goldman-Millson groupoid $\mathbb{DGM}_{\mathfrak{g}}(M_{(E,\mathbf{F},\nabla)},\varepsilon,A)$ for the evaluation map $\varepsilon:M^0_{(E,\mathbf{F},\nabla)}\rightarrow\mathfrak{g} $ is equivalent the $\mathfrak{g}$-augmented Deligne-Goldman-Millson groupoid $\mathbb{DGM}_{\mathfrak{g}}(M_{(E,\mathbf{F},\nabla)},\varepsilon,A)$ for the evaluation map $\varepsilon:M^0_{(E,\nabla)}\rightarrow\mathfrak{g} $. By \cite[Lemma 3.9]{gm}, the groupoids $\mathbb{DGM}_{\mathfrak{g}}(M_{(E,\mathbf{F},\nabla)},\varepsilon,A)$ and $\mathbb{DGM}_{\mathfrak{g}}(M_{(E,\nabla)},\varepsilon,A)$ are equivalent to the groupoids $\mathbb{PFB}'_{( E, {\mathbf{F}}, \nabla)}(A)$ and $\mathbb{PFB}'_{( E,  \nabla)}(A)$, respectively. Due to Corollary \ref{zaa}, the analytic germ of $\mathfrak{PRP}(r,\mathbf{d})$ at $(\rho,\mathbf{P})$ pro-represents the functor $\mathbb{F}_{(E,\mathbf{F},\nabla)}: A\longmapsto \mathrm{Iso}\mathbb{PFB}'_{( E, {\mathbf{F}}, \nabla)}(A) $, in particular,  the analytic germ of the non-singular affine variety $R(\Gamma,G)$ at $\rho$ pro-represents the functor $\mathbb{F}_{(E,,\nabla)}: A\longmapsto \mathrm{Iso}\mathbb{PFB}'_{( E,  \nabla)}(A) $. The conclusion immediately follows.
\end{proof}

\subsection{Partial Formality and quadraticity}

\begin{definition}A DGLA $L$ is called \emph{partially formal} if it contains a sub-DGLA $\tilde L$ satisfying
\begin{itemize}
  \item $\tilde L$ is a formal DGLA;
  \item $\tilde L^0=L^0$;
  \item $\tilde L^p$ is nonzero if $L^p$ is nonzero.
\end{itemize}
\end{definition}

\begin{definition}\
\begin{enumerate}
  \item Given a logarithmic flat bundle $(E,\nabla)$ of rank $r$ over $(\bar X,D)$, there is a canonical parabolic structure $\mathbf{J}=\{\mathcal{J}^{(i)}\}_{i=1,\cdots,n}$ on $E$ compatible with $\nabla$, called the \emph{Jordan parabolic structure}.
Define the subspace $E^{(i)}_\lambda$ of $E|_{x_i}$ by
\begin{align*}
 E^{(i)}_\lambda=\{v\in E|_{x_i}: (\mathrm{Res}(\widehat{\nabla},x_i)-\lambda)^{r}v=0\}.
\end{align*}
There are finitely many  $\lambda$  such that $E^{(i)}_\lambda\neq \{0\}$, then take the distinct real parts of these $\lambda$ and arrange them as $ \alpha^{(i)}_0< \alpha^{(i)}_1<\cdots<\alpha^{(i)}_{s_i}$. Define
\begin{align*}
 E^{(i)}_{\alpha^{(i)}_\ell}=\bigoplus_{\mathrm{Re}(\lambda)=\alpha^{(i)}_\ell}E^{(i)}_\lambda, \ell=0, 1,\cdots,s_i,
\end{align*}
then the Jardon parabolic structure $\mathbf{J}$ is given by
\begin{align*}
 \mathcal{J}^{(i)}: 0\varsubsetneq  E^{(i)}_{\alpha^{(i)}_{s_i}}\subset  E^{(i)}_{\alpha^{(i)}_{s_i}}\bigoplus E^{(i)}_{\alpha^{(i)}_{s_i-1}}\subset\cdots\subset E^{(i)}_{\alpha^{(i)}_{s_i}}\bigoplus \cdots\bigoplus E^{(i)}_{\alpha^{(i)}_{1}}\subset E|_{x_i}.
\end{align*}
We say $(E,\nabla)$ is \emph{Jordan stable} if the parabolic logarithmic flat bundle $(E,\mathbf{J},\nabla)$ is stable with respect to the weight system $\mathbf{w}=\{\overrightarrow{w^{(i)}}\}_{i=1,\cdots,n}$, where $\overrightarrow{w^{(i)}}=(\alpha^{(i)}_0,\cdots, \alpha^{(i)}_{s_i})\in \mathbb{R}^{s_i+1}$.
  \item A logarithmic flat bundle $(E,\nabla)$  over $(\bar X,D)$ is called \emph{perfect} if  it satisfies the following conditions
\begin{itemize}
  \item  the residues of $\nabla$ at the punctures are semisimple and all their eigenvalues lie in $[0,1)$;
    \item  $(E,\nabla)$ is Jordan stable.
    \end{itemize}
\end{enumerate}

\end{definition}

\begin{theorem}If $(E,\nabla)$ is a perfect logarithmic flat bundle over $(\bar X,D)$,
then the DGLA $M_{(E,\mathbf{J},\nabla)}$ is partially formal.
\end{theorem}

\begin{proof}
Let $\mathbb{V}$ be the local system over $X$ corresponding to the flat bundle $E|_X$,
then   $\End(\mathbb{V})$ is the local system over $X$ corresponding to the flat bundle $\End(E)|_X$,
 Let $\rho(\gamma_i)$ be the monodromy action around $z_i$ on $\mathbb{V}$, and chooses a stalk $V^{(i)}$ of  $ \mathbb{V}$ around $x_i$. Since $\rho(\gamma_i)$ is semisimple, there is an eigenvalue decomposition $V^{(i)}=\bigoplus\limits_{\alpha}{V}^{(i)}_\alpha$, where ${V}^{(i)}_\alpha$  is an eigenspace of  $\rho(\gamma_i)$-action on $V^{(i)}$ with the eigenvalue $\alpha$ (by assumption, $|\alpha|=1$).
Each ${V}^{(i)}_\alpha$ produces  a local system $ \mathbb{V}^{(i)}_\alpha$  over a  small deleted neighborhood
 $ U_i$ around $x_i$. We define
 \begin{align*}
    \mathbb{H}^{(i)}&=\bigoplus_{\alpha}\mathbb{V}^{(i)}_\alpha\otimes_{\mathbb{C}}(\mathbb{V}^{(i)}_\alpha)^*,\\
   \mathbb{K}^{(i)}&=\bigoplus_{\frac{\log \alpha}{2\pi\sqrt{-1}}\geqslant\frac{\log \beta}{2\pi\sqrt{-1}}}\mathbb{V}^{(i)}_\alpha\otimes_{\mathbb{C}}(\mathbb{V}^{(i)}_\beta)^*,\\
    \mathbb{N}^{(i)}&=\bigoplus_{\frac{\log \alpha}{2\pi\sqrt{-1}}>\frac{\log \beta}{2\pi\sqrt{-1}}}\mathbb{V}^{(i)}_\alpha\otimes_{\mathbb{C}}(\mathbb{V}^{(i)}_\beta)^*,
 \end{align*}where we choose $\frac{\log \alpha}{2\pi\sqrt{-1}}$ and $ \frac{\log \beta}{2\pi\sqrt{-1}}$ to lie in $[0,1)$.
 Over $\bar U_i=U_i\bigcup\{x_i\}$ with the local coordinate $z^{(i)}$, $\End(E,\mathbf{J})|_{\bar U_i}$ is generated by $z^{(i)}\End(\mathbb{V})$ and $ \mathbb{K}^{(i)}$. Define a subsheaf $\mathrm{SEnd}(E,\mathbf{J})$ of $\End(E,\mathbf{J})$ by $\mathrm{SEnd}(E,\mathbf{J})|_X=\End(E)|_X$ and $\mathrm{SEnd}(E,\mathbf{J})|_{\bar U_i}$ being  generated by $z^{(i)}\End(\mathbb{V})$ and $ \mathbb{N}^{(i)}$ for each $i$.
 Then we have a subcomplex $\mathcal{SC}_{(E,\mathbf{J},\nabla)}$ of $\mathcal{C}_{(E,\mathbf{J},\nabla)}$ given by
 \begin{align*}
 \mathcal{SC}_{(E,\mathbf{J},\nabla)}:  \mathrm{End}(E,\mathbf{J})\xrightarrow{\nabla}\mathrm{SEnd}(E,\mathbf{J})\otimes_{\mathcal{O}_{\bar X}}\Omega^1_{\bar X}(\mathcal{D}),
 \end{align*}
and a sub-DGLA $\tilde M_{(E,\mathbf{J},\nabla)}$ of $M_{(E,\mathbf{J},\nabla)}$ defined by
\begin{itemize}
  \item $\tilde M^0_{(E,\mathbf{J},\nabla)}=M^0_{(E,\mathbf{J},\nabla)}$;
  \item $\tilde M^p_{(E,\mathbf{J},\nabla)}=\Gamma(\bar X, \bigoplus\limits_{a+b=p} \mathcal{A}^{0,a}_{\bar X}\otimes_{\mathcal{O}_{\bar X}}(\mathrm{SEnd}(E,\mathbf{J})\otimes_{\mathcal{O}_{\bar X}}\Omega^b_{\bar X}(\mathcal{D})))$, $p=1,2$.
\end{itemize}
By Esnault--Viehweg--Timmerscheidt's result (see \cite[Lemma 1.6]{e} and \cite[Proposition D.2]{e}), $\iota_*\End(\mathbb{V})$ for the embedding  $\iota: X\rightarrow \bar X$ is quasi-isomorphic to the complex $ \mathcal{SC}_{(E,\mathbf{J},\nabla)}$.

The assumption on eigenvalues of the residues of $\nabla$ at each puncture $x_i$ guarantees that the eigenvalues of the residues of the induced connection $\nabla$ on $\mathrm{End}(E,\mathbf{J})$ lie in $[0,1)$.
 Since $(E,\nabla)$ is Jordan stable, by Simpson's famous work \cite{si},  there exists a harmonic metric $h$ on $\End(E)|_X$ such that, for $v\in \mathbb{V}^{(i)}_\alpha\otimes_{\mathbb{C}}(\mathbb{V}^{(i)}_\beta)^*$,  we have
\begin{align*}
 ||v||_h=\mathcal{O}(|z^{(i)}|^{|\frac{\log\alpha-\log\beta}{2\pi \sqrt{-1}}|-\varepsilon})
\end{align*}
for any $\varepsilon>0$.
On the other hand, endow $X$ with a complete Poincar\'{e}-type K\"{a}hler metric $\omega_{\textrm{PK}}$ along $D$, namely, $\omega_{\mathrm{PK}}$ has an asymptotic form
\begin{align*}
\frac{\sqrt{-1}}{2}\frac{dz^{(i)}\wedge d\overline{z^{(i)}}}{|z^{(i)}|^2\log^2|z^{(i)}|^2}
\end{align*}near the puncture $x_i$.
We define $\mathcal{L}^{p}_{d_\nabla,(2)}(\End(\mathbb{V}))$ to be the sheaf over $\bar X$ of measurable $p$-forms $\sigma$ valued in $\End(\mathbb{V})$ such that both $\sigma$ and $d_\nabla\sigma$ are locally square-integrable with respect to the harmonic metric $h$ and the K\"{a}hler metric $\omega_{\mathrm{PK}}$.  By the result of Zucker--Timmerscheidt (see \cite[Theorem 6.2]{z} and \cite[Proposition D.4]{e}), $\iota_*\End(\mathbb{V})$ is also quasi-isomorphic to the complex
\begin{align*}
 {\mathcal{L}}_{(E,{\nabla})}: \mathcal{L}^{0}_{d_\nabla,(2)}(\End(\mathbb{V}))\xrightarrow{d_\nabla}\mathcal{L}^{1}_{d_\nabla,(2)}(\End(\mathbb{V}))\xrightarrow{d_\nabla}\mathcal{L}^{2}_{d_\nabla,(2)}(\End(\mathbb{V})).
 \end{align*}
\begin{lemma}The local  sections $ \mathrm{End}(E,\mathbf{J})$ and $\mathrm{SEnd}(E,\mathbf{J})\otimes_{\mathcal{O}_{\bar X}}\Omega^1_{\bar X}(\mathcal{D})$ over $\bar U_i$ are precisely the square-integrable local  sections of $\iota_*({\End(E)}|_{U_i})$ and $\iota_*(\mathrm{End}(E)|_{U_i}\otimes_{\mathcal{O}_{U_i}}\Omega^1_{U_i})$ with respect to $h$ and  $\omega_{\mathrm{PK}}$, respectively.
\end{lemma}
\begin{proof}This lemma follows local estimations by means of the following elementary facts:
\begin{align*}
|\frac{dz^{(i)}}{z^{(i)}}|^2_{\omega_{\mathrm{PK}}}&\sim \log^2t,\\
  \int_{0}^{2\pi}\int_0^A t^n|\log t|^m\frac{tdtd\theta }{t^2\log^2t}<+\infty &\textrm{ when }m\leq 0\textrm{ and } n\geq 0,\textrm{ or }m>0\textrm{ and }n>0,
\end{align*}
where one writes $z^{(i)}=te^{\sqrt{-1}\theta}$.
\end{proof}

One introduces
\begin{align*}
  D''&=\frac{1}{2}(d_\nabla+\delta''_h-\delta'_h),\\
  D'&=\frac{1}{2}(d_\nabla+\delta'_h-\delta''_h)
\end{align*}
with $\delta'_h, \delta''_h$ given by, respectively,
\begin{align*}
  \bar\partial h(u,v)&=h(\bar\partial_Eu,v)+h(u,\delta'_h v),\\
  \partial h(u,v)&=h(\nabla u,v)+h(u,\delta''_h v)
\end{align*}
for $u,v\in C^\infty(\End(E))$.
Then we similarly define the sheaves $\mathcal{L}^{p,q}_{D',(2)}(\End(\mathbb{V}))$ and $\mathcal{L}^{p,q}_{D'',(2)}(\End(\mathbb{V}))$.
Let $L^{p}_{d_\nabla,(2)}(\End(\mathbb{V}))$ be the space of square-integrable global $\End(\mathbb{V})$-valued $p$-forms with measurable coefficients on $X$. Similarly, we can define the spaces $L^{p,q}_{D',(2)}(\End(\mathbb{V}))$ and $L^{p,q}_{D'',(2)}(\End(\mathbb{V}))$. We treat the operators $\bullet$ as densely-defined closed operators acting on the space $L^{*/(*,*)}_{\bullet,(2)}(\End(\mathbb{V}))$ with maximal domain $\mathrm{Dom}(\bullet^{*/(*,*)})$. Let $\bullet^\dagger$ be the formal adjoint of $\bullet$ with respect to $h$ and $\omega_{\mathrm{PK}}$, and let the Laplacian operator of $\bullet$ be $\square _\bullet=\bullet\circ\bullet^\dagger+\bullet^\dagger\circ\bullet$. Then
we define the following spaces
  \begin{align*}
    H^p_{d_\nabla,(2)}&=\{\varphi\in \mathrm{Dom}(\square ^p_{d_\nabla}): \square_{d_\nabla}\varphi=0\}\\
    &=\{\varphi\in \mathrm{Dom}({d_\nabla^p})\bigcap \mathrm{Dom}({(d_\nabla^\dagger)^p}): d_\nabla\varphi={d_\nabla^\dagger}\varphi=0\},\\
    H^{p,q}_{D', (2)}&=\{\varphi\in \mathrm{Dom}(\square ^{p,q}_{D'}): \square_{D'}\varphi=0\}\\
    &=\{\varphi\in \mathrm{Dom}({D'^p})\bigcap \mathrm{Dom}({(D'^\dagger)^p}): D'\varphi={D'^\dagger}\varphi=0\},\\
    H^{p,q}_{D'', (2)}&=\{\varphi\in \mathrm{Dom}(\square ^{p,q}_{D''}): \square_{D''}\varphi=0\}\\
    &=\{\varphi\in \mathrm{Dom}({D''^p})\bigcap \mathrm{Dom}({(D''^\dagger)^p}): D''\varphi={D''^\dagger}\varphi=0\},
  \end{align*}
Since $h$ is a harmonic metric, we have the K\"{a}hler identities
   \begin{align*}
     \square _{d_\nabla}=2\square _{D'}=2\square _{D''},
   \end{align*}
   which lead to the decomposition \cite{e,z}
\begin{align*}
  \Gamma(\bar X, \mathcal{L}^{p}_{(2)}(\mathbb{\End(\mathbb{V})}))&\simeq H^p_{d_\nabla,(2)}\bigoplus \square_{d_\nabla}\big(\Gamma(\bar X, \mathcal{L}^{p}_{d_\nabla,(2)}(\End(\mathbb{V})))\big)\\
  & \simeq \bigoplus_{a+b=p}\Big(H^{a,b}_{D',(2)}\bigoplus \square_{D'}\big(\Gamma(\bar X, \mathcal{L}^{a,b}_{D',(2)}(\mathbb{V}))\big)\Big)\\
   & \simeq \bigoplus_{a+b=p}\Big(H^{a,b}_{D'',(2)}\bigoplus \square_{D''}\big(\Gamma(\bar X, \mathcal{L}^{a,b}_{D'',(2)}(\mathbb{V}))\big)\Big).
\end{align*}
As a consequence, the $D'D''$-lemma ensures the quasi-isomorphisms of DGLAs
\begin{align*}
 (\tilde M_{(E,\mathbf{J},\nabla)},d_\nabla)\leftarrow (\Ker(D''|_{\tilde M_{(E,\mathbf{J},\nabla)}}),D')\rightarrow (H^*(\tilde M_{(E,\mathbf{J},\nabla)}),0).
\end{align*}
Thus, we complete the proof of the formality of the DGLA $\tilde M_{(E,\mathbf{J},\nabla)}$.
\end{proof}

Fixing $\mathbf{e}=\{\overrightarrow{e^{(i)}}\}_{i=1,\cdots,n}=\{(e^{(i)}_1,\cdots,e^{(i)}_{r})\}_{i=1,\cdots,n}$ for $e^{(i)}_\mu=e^{2\pi\sqrt{-1}a^{(i)}_\mu}$ with $a^{(i)}_\mu\in[0,1)$, let $\mathfrak{PRP}(r,\mathbf{d},\mathbf{e})$ be a subvariety of $\mathfrak{PRP}(r,\mathbf{d})$ consisting of pairs of parabolic representations $(\rho,\mathbf{P})$ such that the eigenvalues of $\rho(\gamma_i)$ are exactly $e^{(i)}_1,\cdots,e^{(i)}_{r}$ for each $i$.
\begin{corollary} If $\mathfrak{PRP}(r,\mathbf{d},\mathbf{e})$ contains a parabolic representation pair $(\rho,\mathbf{P})$ corresponding to a perfect logarithmic flat bundle, then
$\mathfrak{PRP}(r,\mathbf{d},\mathbf{e})$ is quadratic at $(\rho,\mathbf{P})$.
\end{corollary}

\begin{proof}It is known that the deformation of $(\rho,\mathbf{P})$ in $\mathfrak{PRP}(r,\mathbf{d},\mathbf{e})$  is controlled by the DGLA $\tilde M_{(E,\mathbf{J},\nabla)}$ (cf. \cite{bf,foth}). Using the  standard arguments by  Goldman--Millson  in \cite{gm}, we conclude the quadraticity.
\end{proof}

\section{Moduli spaces of weighted parabolic representation pairs}

Consider the flag variety $\mathrm{FL}(r,\overrightarrow{d})$. It is identified with $G/P_{\mathcal{F}}$ for a certain parabolic subgroup $P_{\mathcal{F}}$ of $G$ preserving a flag $\mathcal{F}$ of type $\overrightarrow{d}=(d_0,\cdots,d_s)$. Let $T_{\mathcal{F}}$ and $B_{\mathcal{F}}$ be the maximal torus and  the Borel subgroup contained in $P_{\mathcal{F}}$, respectively, and let $X(T_{\mathcal{F}})$ and $ Y(T_{\mathcal{F}})$ be the character group and cocharacter group of $T_{\mathcal{F}}$, respectively, with the natural pairing $\langle,\rangle:X(T_{\mathcal{F}})\times Y(T_{\mathcal{F}})\rightarrow\mathbb{Z}$. Denote by $S$ the set of simple roots with respect to $T_{\mathcal{F}}$. Then $P_{\mathcal{F}}$ is determined by a subset $I\subset S$; namely, $P_{\mathcal{F}}$ is generated by $B_{\mathcal{F}}$ and the root subgroups $U_{-\upsilon}$ of $G$ corresponding to the root $-\upsilon\in I$. The character group of $P_{\mathcal{F}}$ is given by
\begin{align*}
  X(P_{\mathcal{F}})=\{\delta\in X(T_{\mathcal{F}}): \langle\delta,\upsilon^\vee\rangle=0\textrm{ for all } \upsilon\in I\},
\end{align*}
where  $\upsilon^\vee \in Y(T_{\mathcal{F}})$ denotes the coroot corresponding to $\upsilon$.
Let $\delta: P_{\mathcal{F}}\rightarrow\mathbb{C}^*$  be a character of $P_{\mathcal{F}}$. We define an action of $P_{\mathcal{F}}$ on the product $G\times \mathbb{C}$  by
\begin{align*}
 h\cdot(g, c)=(gh^{-1},\delta(h)c),
\end{align*}
 then we have a quotient
$L_\delta=(G\times \mathbb{C})/P_{\mathcal{F}}$, which maps to $G/P_{\mathcal{F}}$ via $P_{\mathcal{F}}\cdot(g,c)\mapsto gP_{\mathcal{F}}$. Thus, we obtain a line bundle
$$
\pi:L_\delta\longrightarrow\mathrm{FL}(r,\overrightarrow{d})\simeq G/P_{\mathcal{F}}.
$$
The group $G$ acts on $L_\delta$ via $g'\cdot (g,c)=(g'g,c)$. One easily checks that the projection $\pi$ is $G$-equivariant and that any $g\in G$ induces a linear map from the fiber $\pi^{-1}(F)$ to the fiber  $\pi^{-1}(gF)$ for $F\in \mathrm{FL}(r,\overrightarrow{d})$. In other words, $L_\delta$ is a $G$-linearized line bundle over $ \mathrm{FL}(r,\overrightarrow{d})$.

The following proposition is well-known \cite{br}.

\begin{proposition}
$L_\delta$ is ample if and only if $\langle\delta,\upsilon^\vee\rangle>0$  for all $\upsilon\in S\backslash I$.
\end{proposition}

For our case $G=\mathrm{GL}(V)$,   the Levi subgroup $L_{\mathcal{F}}$ of $P_{\mathcal{F}}$ is isomorphic to $\mathrm{GL}(d_0,\mathbb{C})\times\cdots\times \mathrm{GL}(d_{s},\mathbb{C})$ via a suitable choice of basis of $V$ adapted to the flag $\mathcal{F}$. Then the composition of projections $P_{\mathcal{F}}\rightarrow L_{\mathcal{F}}$ and $L_{\mathcal{F}}\rightarrow \mathrm{GL}(d_j,\mathbb{C})$ is denoted by $\mathfrak{l}_{j}: P_{\mathcal{F}}\rightarrow \mathrm{GL}(d_j,\mathbb{C})$. Given a weight $\overrightarrow{w}=(w_0,\cdots,w_s)\in \mathbb{Z}^{s+1}$ with $w_0<\cdots<w_s$, we have a character of $P_{\mathcal{F}}$ defined by
\begin{align*}
 \delta_{\overrightarrow{w}}: P_{\mathcal{F}}&\longrightarrow\mathbb{C}^*\\
 h&\longmapsto \det(\mathfrak{l}_{0}(h))^{-w_0}\cdots\det(\mathfrak{l}_{s}(h))^{-w_s}.
\end{align*}
The corresponding line bundle is denoted by $L_{\overrightarrow{w}}$, which is ample.
For the multiple flag variety $\mathrm{FL}(r,\mathbf{d})$ and a weight system $\mathbf{w}=\{\overrightarrow{w^{(i)}}\}_{i=1,\cdots,n}$, where $\overrightarrow{w^{(i)}}=(w^{(i)}_0,\cdots, w^{(i)}_{s_i})\in \mathbb{Z}^{s_i+1}$ satisfying $w^{(i)}_0<\cdots <w^{(i)}_{s_i}$, we define a line bundle
\begin{align*}
  L_{\mathbf{w}}=L_{\overrightarrow{w^{(1)}}}\otimes\cdots\otimes L_{\overrightarrow{w^{(n)}}}
\end{align*}
over $\mathrm{FL}(r,\mathbf{d})$.

On the other hand, since $R(\Gamma,G)$ is an affine variety, the trivial line bundle $\mathcal{O}_{R(\Gamma,G)}$ is an ample line bundle over $R(\Gamma,G)$. Then, by considering the fibrations
$$
f_1: R(\Gamma,G)\times\mathrm{FL}(r,\mathbf{d}) \rightarrow R(\Gamma,G),\quad f_2: R(\Gamma,G)\times\mathrm{FL}(r,\mathbf{d})\rightarrow \mathrm{FL}(r,\mathbf{d}),
$$
we get a $G$-linearized ample line bundle
$$
\mathcal{L}_{\mathbf{w}}=f_1^*\mathcal{O}_{R(\Gamma,G)}\otimes f_2^*L_{\mathbf{w}}
$$
over $ R(\Gamma,G)\times\mathrm{FL}(r,\mathbf{d}) $. Since $\mathfrak{PRP}(r,\mathbf{d})$ is a closed subvariety of $R(\Gamma,G)\times\mathrm{FL}(r,\mathbf{d})$, the restriction of $\mathcal{L}_{\mathbf{w}}$ on $\mathfrak{PRP}(r,\mathbf{d})$ is again a $G$-linearized ample line bundle over $\mathfrak{PRP}(r,\mathbf{d})$, which we also denote by $\mathcal{L}_{\mathbf{w}}$.

 \begin{proposition}
 Given a weight system $\mathbf{w}=\{\overrightarrow{w^{(i)}}\}_{i=1,\cdots,n}$, where $\overrightarrow{w^{(i)}}=(w^{(i)}_0,\cdots, w^{(i)}_{s_i})\in \mathbb{Z}^{s_i+1}$ satisfies $w^{(i)}_0<\cdots <w^{(i)}_{s_i}$. A weighted parabolic representation pair $((\rho,\mathbf{P}),\mathbf{w})$ of degree zero is semistable if and only if it is GIT-semistable in $\mathfrak{PRP}(r,\mathbf{d})$ with respect to the $G$-linearized ample line bundle $ \mathcal{L}_{\mathbf{w}}$.
 \end{proposition}

\begin{proof}
Consider a one-parameter subgroup $\lambda: \mathbb{C}^*\rightarrow G$, then the vector space $V$ decomposes into the direct sum of weight spaces for the action of $\lambda$, namely
\begin{align*}
 V=\bigoplus\limits_{N\in\mathbb{N}} V^N,
\end{align*}
where
\begin{align*}
 \lambda(c)\cdot v=c^Nv
\end{align*}
  for every $v\in V^N$. We write $V^{\geq N}=\bigoplus\limits_{M\geq N}V^M$. For any  $\omega\in \Gamma$, the map $\rho(\omega)$ also decomposes into the direct sum of linear maps between weight spaces, namely
  \begin{align*}
   \rho(\omega)=\bigoplus\limits_{M,N\in\mathbb{N}} \rho(\omega)^{MN},
  \end{align*}
  where $\rho(\omega)^{MN}$ is the compositions of the embedding $\iota^M: V^M\hookrightarrow V$, the linear map $\rho(\omega): V\rightarrow V$, and the projection $\rho^N:V\rightarrow V^N$.

  Thus  we have
  \begin{align*}
 \lambda(c)\cdot \rho(\omega)=\mathrm{  Ad}_{\lambda(c)}\rho(\omega)=\bigoplus\limits_{M,N\in\mathbb{N}}c^{N-M}\rho(\omega)^{MN}.
  \end{align*}
  If the limit $\lim\limits_{c\rightarrow 0} \lambda(c)\cdot \rho(\omega)$ exists, we must have $\rho(\omega)^{MN}=0$ whenever $N<M$. It follows that each $V^{\geq N}$ is a subrepresentation of $\rho$. On the other hand, if the  limit $\lim\limits_{c\rightarrow 0} \lambda(c)\cdot \mathbf{P}$ exists, we must have $V^N=\{0\}$ for $N<0$. Then we calculate the Mumford weight $\mu^{\mathcal{L}_{\mathbf{w}}}_\lambda(\rho,\mathbf{P})$ of the $\lambda$-action on $\mathcal{L}_{\mathbf{w}}$ as follows (cf. \cite[Proposition 3.7]{t})
  \begin{align*}
   \mu^{\mathcal{L}_{\mathbf{w}}}_\lambda(\rho,\mathbf{P})=-\sum_{N\geq 0}N\sum_{i=1}^n\sum_{\ell=0}^{s_i}w^{(i)}_\ell\dim_\mathbb{C}(V^{(i)}_\ell\bigcap V^N)/(V^{(i)}_{\ell+1}\bigcap V^N).
  \end{align*}
  Since $((\rho,\mathbf{P}),\mathbf{w})$ is of zero degree and semistable, each summand with respect to the index $N$ is non-positive. Consequently,
  \begin{align*}
    \mu^{\mathcal{L}_{\mathbf{w}}}_\lambda(\rho,\mathbf{P})\geq 0,
  \end{align*}
   so $(\rho,\mathbf{P})$  is GIT-semistable in $\mathfrak{PRP}(r,\mathbf{d})$ with respect to  $ \mathcal{L}_{\mathbf{w}}$. Conversely,  the non-negativity of the Mumford weight $\mu^{\mathcal{L}_{\mathbf{w}}}_\lambda(\rho,\mathbf{P})$ for a one-parameter subgroup $\lambda: \mathbb{C}^*\rightarrow G$ corresponding to a proper subrepresentation immediately implies that $((\rho,\mathbf{P}),\mathbf{w})$ is semistable.
\end{proof}

By classical geometric invariant theory, we have the following corollary.

\begin{corollary}
Given a rational weight system $\mathbf{w}=\{\overrightarrow{w^{(i)}}\}_{i=1,\cdots,n}$, where $\overrightarrow{w^{(i)}}=(w^{(i)}_0,\cdots, w^{(i)}_{s_i})\in \mathbb{Q}^{s_i+1}$ satisfies $w^{(i)}_0<\cdots <w^{(i)}_{s_i}$. Let $\mathfrak{PRP}^{\mathbf{w}\text{-}\mathrm{ss}}(r,\mathbf{d})$ denote the open subset of $\mathfrak{PRP}(r,\mathbf{d})$ consisting of semistable parabolic representation pairs of degree zero with respect to $\mathbf{w}$, then $\mathfrak{PRP}^{\mathbf{w}\text{-}\mathrm{ss}}(r,\mathbf{d})$ admits a good categorical GIT quotient
\begin{align*}
  q: \mathfrak{PRP}^{\mathbf{w}\text{-}\mathrm{ss}}(r,\mathbf{d})\longrightarrow \mathfrak{PRP}(r,\mathbf{d})/\!\!/_{ \mathcal{L}_{\mathbf{w}}} G.
\end{align*}
 Moreover, $\mathfrak{PRP}(r,\mathbf{d})/\!\!/_{ \mathcal{L}_{\mathbf{w}}} G$ has the following properties:
\begin{itemize}
  \item it is a quasi-projective variety;
  \item it co-represents the moduli functor $\mathcal{M}:\mathbf{Sch}^{\mathrm{op}}\rightarrow \mathbf{Set}$ from the opposite category $\mathbf{Sch}^{\mathrm{op}}$ of $\mathbb{C}$-schemes to the category   $\mathbf{Set}$, which sends a scheme $S$ to the set of isomorphism classes of flat families of semistable  parabolic representation pairs of rank $r$ and type $\mathbf{d}$ of degree zero with respect to $\mathbf{w}$, parametrized by $S$;
  \item its closed points correspond to the  isomorphism classes of polystable parabolic representation pairs  of rank $r$ and type $\mathbf{d}$ of degree zero with respect to $\mathbf{w}$.
\end{itemize}
\end{corollary}

Another approach to constructing the moduli space of weighted parabolic representation pairs is to translate it into the moduli space of representations of a certain quiver. More precisely, we construct the following star-shaped quiver $Q = (Q_0,Q_1,s,t)$ with loops, where
\begin{itemize}
  \item $Q_0=\{u\}\cup\{u_{\ell}^{(i)}\}_{1\le\ell\le s_i,1\le i\le n}$ is the set of vertices;
  \item $Q_1 = \{a_j,b_j\}_{1\le j\le g}\cup\{c^{(i)}\}_{1\le i\le n}\cup\{e_{\ell}^{(i)}\}_{1\le\ell\le s_i,1\le i\le n}\cup\{c_{\ell}^{(i)}\}_{1\le\ell\le s_i,1\le i\le n}$ is the set of arrows;
  \item the source and target maps $s,t: Q_1 \to Q_0$ are given by
  \begin{align*}
    s(a_j)=u, \qquad & t(a_j)=u, \\
    s(b_j)=u, \qquad & t(b_j)=u, \\
    s(c^{(i)})=u, \qquad & t(c^{(i)})=u, \\
    s(e_{\ell}^{(i)})=u_{\ell}^{(i)}, \qquad & t(e_{\ell}^{(i)})=u_{\ell-1}^{(i)} (\textrm{one assigns }u_{0}^{(i)}=u), \\
    s(c_{\ell}^{(i)})=u_{\ell}^{(i)}, \qquad & t(c_{\ell}^{(i)})=u_{\ell}^{(i)}.
\end{align*}
\end{itemize}
That is, the quiver $Q$ is depicted as follows:
\begin{figure}[H] 
\centering 
\includegraphics[width=0.7\textwidth]{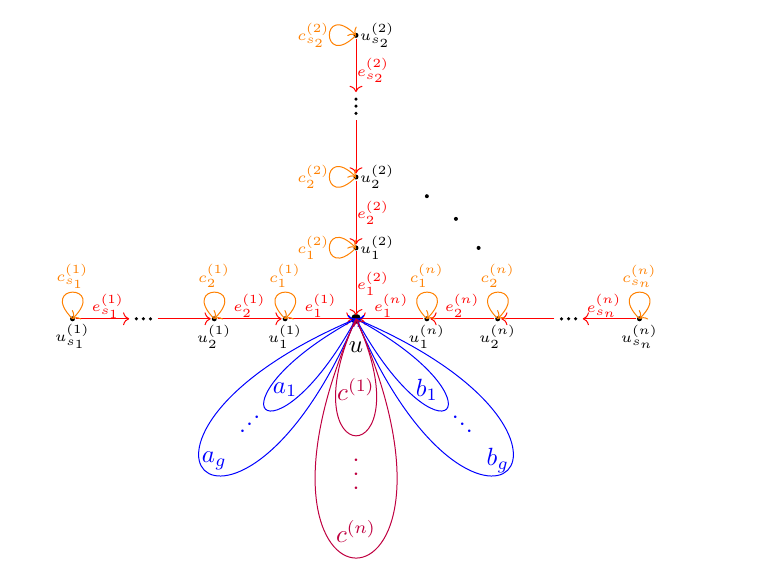} 
\end{figure}
Let us consider representations of $Q$ with dimension vector $\overrightarrow{d_Q}= (d_v)_{v \in Q_0}$ induced from the type $\mathbf{d}$, where
 \begin{itemize}
   \item $d_u = r$;
   \item $d_{u_{\ell}^{(i)}} := \sum\limits_{\jmath\geq \ell} d_\jmath^{(i)}$ for $1 \le \ell \le s_i$ and $1 \le i \le n$.
 \end{itemize}
 The representation space of $Q$ with dimension vector $\overrightarrow{d_Q}$ is the affine space given by
 \begin{align*}
   \mathfrak{REP}(Q,\overrightarrow{d_Q}) = \bigoplus\limits_{a \in Q_1} \Hom_{\mathbb{C}}(\mathbb{C}^{d_{s(a)}},\mathbb{C}^{d_{t(a)}}).
 \end{align*}
 Let $\mathrm{GL}(\overrightarrow{d_Q})$ denote the algebraic group $\prod\limits_{v \in Q_0} \mathrm{GL}(d_v,\mathbb{C})$. There is a canonical algebraic action of $\mathrm{GL}(\overrightarrow{d_Q})$ on $\mathfrak{REP}(Q,\overrightarrow{d_Q})$, given by
  \begin{align*}
    (g \cdot x)_a := g_{t(a)} \circ x_a \circ g_{s(a)}^{-1}
  \end{align*}
 for $x \in \mathfrak{REP}(Q,\overrightarrow{d_Q})$, $g \in \mathrm{GL}(\overrightarrow{d_Q})$, and $a \in Q_1$.

Define $\mathfrak{REP}_{\mathrm{prp}}(Q,\overrightarrow{d_Q})$ to be the subset of $\mathfrak{REP}(Q,\overrightarrow{d_Q})$ consisting of tuples $(x_a)_{a \in Q_1}$ satisfying the following conditions:
\begin{itemize}
\item $x_a$ is injective for all $a \in Q_1$;
\item $x_{c_{\ell-1}^{(i)}} \circ x_{e_{\ell}^{(i)}} = x_{e_{\ell}^{(i)}} \circ x_{c_{\ell}^{(i)}}$ for all $1 \le \ell \le s_i$ and $1 \le i \le n$;
\item $[x_{a_1},x_{b_1}] \circ \cdots \circ [x_{a_g},x_{b_g}] \circ x_{c^{(1)}} \circ \cdots \circ x_{c^{(n)}} = \mathrm{Id}$.
\end{itemize}
Clearly, $\mathfrak{REP}_{\mathrm{prp}}(Q,\overrightarrow{d_Q})$ is a locally closed subset of $\mathfrak{REP}(Q,\overrightarrow{d_Q})$ and is invariant under the $\mathrm{GL}(\overrightarrow{d_Q})$-action.

\begin{proposition}
There is a natural $\mathrm{GL}(\overrightarrow{d_Q})$-equivariant morphism $\varphi$ from $\mathfrak{REP}_{\mathrm{prp}}(Q,\overrightarrow{d_Q})$ to $\mathfrak{PRP}(r,\mathbf{d})$. Moreover, $\varphi$ is a geometric quotient.
\end{proposition}

\begin{proof}
The morphism $\varphi$ is defined by $(x_a)_{a \in Q_1} \mapsto (\rho,\mathbf{P})$, where:
\begin{itemize}
\item $\rho : \Pi \to\mathrm{GL}(V)$ is given by
$$
\rho(\alpha_j) =x_{a_j},\ \ \rho(\beta_j )= x_{b_j},\  \ \rho(\gamma_i) = x_{c^{(i)}}\qquad (1 \le j \le g,\ 1 \le i \le n);
$$
\item $\mathbf{P} = (P_{\mathcal{F}^{(1)}},\cdots,P_{\mathcal{F}^{(n)}})$ with $V_{\ell}^{(i)}$ being the subspace of $V$ given by the image of $ x_{e_{\ell}^{(i)}}\circ\cdots x_{e_{1}^{(i)}}$.
\end{itemize}
The identity
\begin{align*}
x_{c^{(i)}} \circ x_{e_{1}^{(i)}} \circ \cdots \circ x_{e_{\ell}^{(i)}} = x_{e_{1}^{(i)}} \circ \cdots \circ x_{e_{\ell}^{(i)}} \circ x_{c_{\ell}^{(i)}}
\end{align*}
implies that $ x_{c^{(i)}} = \rho(\gamma_i)\in P_{\mathcal{F}^{(i)}} $. One easily checks that
\begin{align*}
 \varphi(g \cdot x) = g_u \cdot \varphi(x)
\end{align*}
for all $g \in\mathrm{GL}(\overrightarrow{d_Q})$ and $x \in \mathfrak{REP}_{\mathrm{prp}}(Q,\overrightarrow{d_Q})$, where $g_u \in\mathrm{GL}(V)$ denotes the projection of $g$ onto $\mathrm{GL}(d_u,\mathbb{C})=\mathrm{GL}(V)$.

Let $G'=\prod\limits_{v \ne u} \mathrm{GL}(d_v,\mathbb{C})$. The canonical embedding $G' \hookrightarrow \mathrm{GL}(\overrightarrow{d_Q})$ induces a natural $G'$-action on $\mathfrak{REP}_{\mathrm{prp}}(Q,\overrightarrow{d_Q})$ and a trivial $G'$-action on $\mathfrak{PRP}(r,\mathbf{d})$. Then the morphism
$$
\varphi : \mathfrak{REP}_{\mathrm{prp}}(Q,\overrightarrow{d_Q}) \longrightarrow \mathfrak{PRP}(r,\mathbf{d})
$$
is a geometric quotient for the action of $G'$ on $\mathfrak{REP}_{\mathrm{prp}}(Q,\overrightarrow{d_Q})$. This claim is confirmed as follows. Let $k$ and $l$ be positive integers with $k < l$, and let $\mathrm{Hom}'(k,l)$ be the subset of $\mathbb{C}^{k \times l}$ consisting of the full row-rank matrices. There is a natural morphism $\mathrm{GL}(l,\mathbb{C}) \to \mathrm{Hom}'(k,l)$  given by
$$
\begin{pmatrix}X & Y \\
Z & W\end{pmatrix} \longmapsto (X,Y),
$$
where $X,Y,Z$ and $W$ are matrices of sizes $k \times k, k \times (l-k), (l-k) \times k$, and $(l - k)\times(l- k)$, respectively.
Let $P\subset\mathrm{GL}(l,\mathbb{C})$ be the parabolic subgroup consisting of block matrices of the form  $\left(\begin{smallmatrix}X & 0 \\ Z & W\end{smallmatrix}\right)$.  Hence we have
 \begin{align*}
  \mathrm{Hom}'(k,l) / \mathrm{GL}(k,\mathbb{C})\simeq \mathrm{GL}(l,\mathbb{C}) / P \simeq \mathrm{Gr}(k,l).
 \end{align*}

We thus complete the proof.
\end{proof}

\begin{corollary}
Let $\mathbb{REP}_{\mathrm{prp}}(Q,\overrightarrow{d_Q})$ be the transformation groupoid corresponding to $(\mathfrak{REP}_{\mathrm{prp}}(Q,\overrightarrow{d_Q}),\mathrm{GL}(\overrightarrow{d_Q}))$, then $\mathbb{REP}_{\mathrm{prp}}(Q,\overrightarrow{d_Q})$ is equivalent to $\mathbb{PRP}(r,\mathbf{d})$.
\end{corollary}

We now introduce a weight vector $\overrightarrow{w_Q}=(w_v)_{v\in Q_0}\in \mathbb{Q}^{|Q_0|}$  on the quiver $Q$ induced from a rational weight system $\mathbf{w}$, defined as follows:
\begin{itemize}
  \item $w_u=-\frac{1}{r} \sum\limits_{i=1}^n\sum\limits_{\ell=0}^{s_i}w_{\ell}^{(i)}d_{\ell}^{(i)}+\sum\limits_{i=1}^nw_{0}^{(i)}$;
  \item $w_{u_\ell^{(i)}}=w^{(i)}_\ell-w^{(i)}_{\ell-1}$ for $1 \le \ell \le s_i$ and $1\le i\le n$.
\end{itemize}
The representations in $\mathfrak{REP}(Q,\overrightarrow{d_Q})$ that are semistable with respect to $\overrightarrow{w_Q}$ (i.e. for any proper subrepresentation with dimension vector $\overrightarrow{d'_Q}$ we have $\sum\limits_{v\in Q_0}w_v d'_v\leq 0$) form an open $\mathrm{GL}(\overrightarrow{d_Q})$-invariant subvariety
$$
\mathfrak{REP}^{\overrightarrow{w_Q}\text{-}\mathrm{ss}}(Q,\overrightarrow{d_Q})\subset\mathfrak{REP}(Q,\overrightarrow{d}).
$$
By geometric invariant theory for quiver representations, there exists a good categorical GIT quotient \cite{k}
$$
q: \mathfrak{REP}^{\overrightarrow{w_Q}\text{-}\mathrm{ss}}(Q,\overrightarrow{d_Q})\longrightarrow \mathfrak{REP}(Q,\overrightarrow{d})/\!\!/\mathrm{GL}(\overrightarrow{d_Q}).
$$

\begin{proposition}
Let $(\rho,\mathbf{P}) \in \mathfrak{PRP}(r,\mathbf{d})$ and $x = (x_a)_{a\in Q_1} \in \varphi^{-1}(\rho,\mathbf{P})$. Then $(\rho,\mathbf{P})$ is semistable with respect to $\mathbf{w}$ if and only if $x$ is semistable with respect to $\overrightarrow{w_Q}$. In particular,
$$
\varphi^{-1}(\mathfrak{PRP}^{\mathbf{w}\text{-}\mathrm{ss}}(r,\mathbf{d}))=\mathfrak{REP}^{\overrightarrow{w_Q}\text{-}\mathrm{ss}}_{\mathrm{prp}}(Q,\overrightarrow{d_Q}):=\mathfrak{REP}_{\mathrm{prp}}(Q,\overrightarrow{d_Q})\bigcap\mathfrak{REP}^{\overrightarrow{w_Q}\text{-}\mathrm{ss}}(Q,\overrightarrow{d_Q}).
$$
\end{proposition}

\begin{proof}
It is clear that the semistability of $x$ implies the semistability of $(\rho,\mathbf{P})$. Conversely, we only need to note that  for a proper subrepresentation with dimension vector $\overrightarrow{d'_Q}$, we have the inequalities
 \begin{align*}
  d'_{u^{(i)}_{\ell+1}}\leq \dim_{\mathbb{C}}(\iota^{(i)}_{\ell+1}(\mathbb{C}^{  d_{u^{(i)}_{\ell+1}}})\bigcap   V'_{u^{(i)}_{\ell}})
 \end{align*}
for $1 \le \ell \le s_i$ and $1\le i\le n$, where  $\iota^{(i)}_{\ell}:\mathbb{C}^{  d_{u^{(i)}_{\ell}}}\rightarrow \mathbb{C}^{  d_{u^{(i)}_{\ell-1}}}$ denotes the chosen embedding, and $ V'_{u^{(i)}_{\ell}}$ is a subspace of $\mathbb{C}^{  d_{u^{(i)}_{\ell-1}}}$ arising from the subrepresentation. Since $w_{u_\ell^{(i)}}>0$ for all $1 \le \ell \le s_i$ and $1\le i\le n$, the semistability of $(\rho,\mathbf{P})$ guarantees the  semistability of $x$.
\end{proof}

\begin{corollary}
Assume $\sum\limits_{i=1}^n\sum\limits_{\ell=0}^{s_i}w_{\ell}^{(i)}d_{\ell}^{(i)}=0$. There is a good categorical GIT quotient
 \begin{align*}
   q: \mathfrak{REP}^{\overrightarrow{w_Q}\text{-}\mathrm{ss}}_{\mathrm{prp}}(Q,\overrightarrow{d_Q})\longrightarrow \mathfrak{REP}_{\mathrm{prp}}(Q,\overrightarrow{d})/\!\!/\mathrm{GL}(\overrightarrow{d_Q}),
 \end{align*}
 such that the following diagram commutes:
$$\CD
  \mathfrak{REP}^{\overrightarrow{w_Q}\text{-}\mathrm{ss}}_{\mathrm{prp}}(Q,\overrightarrow{d_Q}) @>q>> \mathfrak{REP}_{\mathrm{prp}}(Q,\overrightarrow{d})/\!\!/\mathrm{GL}(\overrightarrow{d_Q}) \\
  @V \varphi VV @V \simeq VV  \\
  \mathfrak{PRP}^{\mathbf{w}\text{-}\mathrm{ss}}(r,\mathbf{d}) @>q>> \mathfrak{PRP}(r,\mathbf{d})/\!\!/_{ \mathcal{L}_{\mathbf{w}}} G
\endCD.
$$
\end{corollary}

As an application of the correspondence between semistable parabolic representation pairs and semistable representations of the quiver $Q$, we obtain the following Kobayashi--Hitchin-type theorem.

\begin{theorem}
Given a rational weight system $\mathbf{w}$, a  parabolic representation pair $(\rho,\mathbf{P})\in \mathfrak{PRP}(r,\mathbf{d})$ is polystable with respect to $\mathbf{w}$ if and only if there exists a $\mathbf{w}$-good Hermitian metric $h$ on $V$, namely the following equations are satisfied
\begin{align*}
    \left\{
    \begin{aligned}
        &\sum_{j=1}^g ([\rho(\alpha_j),\rho^{*_h}(\alpha_j)]+[\rho(\beta_j),\rho^{*_h}(\beta_j)])+\sum_{i=1}^n[\rho(\gamma_i),\rho^{*_h}(\gamma_i)]+\sum\limits_{i=1}^n{p_h}^{(i)}_1\\
 &\hspace{138pt}=\ (\frac{1}{r} \sum\limits_{i=1}^n\sum\limits_{\ell=0}^{s_i}w_{\ell}^{(i)}d_{\ell}^{(i)}-\sum\limits_{i=1}^nw_{0}^{(i)})\mathrm{Id},\\[3pt]
 &[\rho(\gamma_i)|_{V^{(i)}_\ell},(\rho(\gamma_i)|_{V^{(i)}_\ell})^{*_h}]+{p_h}^{(i)}_\ell
  =\ (w^{(i)}_{\ell-1}-w^{(i)}_\ell+1)\mathrm{Id},\quad 1\leq \ell\leq s_i, 1\leq i\leq n,
    \end{aligned}
    \right.
\end{align*}
where $\rho^{*_h}(\omega)$ and ${p_h}^{(i)}_\ell$ denote the conjugate of $\rho(\omega)$ and the orthogonal projection of $V^{(i)}_{\ell}$ onto $V^{(i)}_{\ell+1}$, respectively, with respect to the metric $h$. Moreover, if $h'$ is another $\mathbf{w}$-good Hermitian metric  on $V$, then there exists $g\in G$ such that
\begin{align*}
  h(u,v)=h'(g\cdot u,g\cdot v)
\end{align*}
for any $u,v\in V$.
\end{theorem}

\begin{proof}
Applying King's result (cf. \cite[Section 6]{k} or \cite{da}) to  polystable representations of the quiver $Q$, there exists a Hermitian metric $h$ on the vector space $\mathbb{C}^r$ and Hermitian metrics $h^{(i)}_\ell$ on the vector spaces  $\mathbb{C}^{d_{u^{(i)}_{\ell}}}$  for  $1\leq \ell\leq s_i$ and $1\leq i\leq n$ such that the following equations hold:
\begin{align*}
\left\{
\begin{aligned}
& \sum_{j=1}^g ([\rho(\alpha_j),\rho^{*_h}(\alpha_j)]+[\rho(\beta_j),\rho^{*_h}(\beta_j)])+\sum_{i=1}^n[\rho(\gamma_i),\rho^{*_h}(\gamma_i)]+\sum\limits_{i=1}^n\iota^{(i)}_1(\iota^{(i)}_1)^{*_{h,h^{(i)}_1}}\\
 &\hspace{173pt}=\ (\frac{1}{r} \sum\limits_{i=1}^n\sum\limits_{\ell=0}^{s_i}w_{\ell}^{(i)}d_{\ell}^{(i)}-\sum\limits_{i=1}^nw_{0}^{(i)})\mathrm{Id},\\
 &[\rho(\gamma_i)|_{V^{(i)}_\ell},(\rho(\gamma_i)|_{V^{(i)}_\ell})^{*_{h^{(i)}_\ell}}]+\iota^{(i)}_{\ell+1}(\iota^{(i)}_{\ell+1})^{*_{h^{(i)}_\ell, h^{(i)}_{\ell+1}}}-(\iota^{(i)}_{\ell})^{*_{h^{(i)}_{\ell-1},h^{(i)}_\ell}}\iota^{(i)}_{\ell}\\
  &\hspace{173pt}=\ (w^{(i)}_{\ell-1}-w^{(i)}_\ell)\mathrm{Id},\quad 1\leq \ell\leq s_i, 1\leq i\leq n,
  \end{aligned}
  \right.
\end{align*}
where $(\iota^{(i)}_{\ell})^{*_{h^{(i)}_{\ell-1},h^{(i)}_\ell}}:\mathbb{C}^{  d_{u^{(i)}_{\ell-1}}}\rightarrow \mathbb{C}^{  d_{u^{(i)}_{\ell}}}$ denotes the conjugate of $\iota^{(i)}_{\ell}$ with respect to the metrics $h^{(i)}_{\ell-1}$ and $h^{(i)}_\ell$.

Now, we put $\iota^{(i)}_{\ell}$ as the natural embedding, then the metrics $h^{(i)}_\ell$ are given by the restriction of the metric $h$ on the corresponding subspaces. Hence, the above equations are simplified accordingly.
\end{proof}

\vspace{3mm}

    \noindent\textbf{Acknowledgements}.
 The authors would like to thank Ya Deng, Yikun Qiao, Junchao Shentu and Carlos Simpson  for their useful discussions, and thank Jonathan Pridham for his helpful communications.   The authors are grateful to the anonymous reviewer for pointing out a  mistake in an earlier version.
  Z. Hu is supported by the National Natural Science Foundation of China (No. 12271253).  R. Zong is supported by the National Natural Science Foundation of China (Key Program 12331002).

 \end{document}